\documentclass{article}
\usepackage{amsmath}
\usepackage{latexsym}
\usepackage{amssymb}
\newcommand{\qed}{\Box}
\newtheorem{thm}{Theorem}[section]
\newtheorem{la}[thm]{Lemma}
\newtheorem{Defn}[thm]{Definition}
\newtheorem{Exam}[thm]{Example}
\newtheorem{Remark}[thm]{Remark}
\newtheorem{prop}[thm]{Proposition}
\newtheorem{cor}[thm]{Corollary}
\newtheorem{Number}[thm]{\!\!}
\newenvironment{defn}{\begin{Defn}\rm}{\end{Defn}}
\newenvironment{ex}{\begin{Exam}\rm}{\end{Exam}}
\newenvironment{rem}{\begin{Remark}\rm}{\end{Remark}}
\newenvironment{numba}{\begin{Number}\rm}{\end{Number}}
\newenvironment{proof}{{\noindent\bf Proof.}}%
                  {\nopagebreak\hspace*{\fill}$\qed$\medskip\par}   
\newcommand{\Punkt}{\nopagebreak\hspace*{\fill}$\qed$}
\newcommand{\wb}{\overline}

\newcommand{\impl}{\Rightarrow}
\newcommand{\mto}{\mapsto}
\newcommand{\ve}{\varepsilon}

\newcommand{\N}{{\mathbb N}}
\newcommand{\R}{{\mathbb R}}

\newcommand{\Z}{{\mathbb Z}}
\newcommand{\C}{{\mathbb C}}
\newcommand{\bL}{{\mathbb L}}

\newcommand{\K}{{\mathbb K}}

\newcommand{\cS}{{\cal S}}

\newcommand{\cL}{{\cal L}}
\newcommand{\cW}{{\cal W}}

\newcommand{\sub}{\subseteq}

\DeclareMathOperator{\pr}{pr}
\DeclareMathOperator{\id}{id}
\DeclareMathOperator{\Lip}{Lip}
\newcommand{\sbull}{{\scriptscriptstyle \bullet}}
\DeclareMathOperator{\Diff}{Diff}

\DeclareMathOperator{\graph}{graph}
\newcommand{\pl}{{\displaystyle\lim_{\longleftarrow}}}

\DeclareMathOperator{\ev}{ev}
\begin{document}
\begin{center}
{\bf\Large Implicit Functions from\\[.9mm]
Topological Vector
Spaces to Fr\'{e}chet Spaces\\[2.1mm]
in the Presence of Metric Estimates}\vspace{3.7mm}\\
{\bf Helge Gl\"{o}ckner}\vspace{2.3mm}
\end{center}
\begin{abstract}\vspace{.4mm}
{\protect\noindent}We prove an implicit function theorem
for Keller $C^k_c$-maps from arbitrary real or complex topological
vector spaces to Fr\'{e}chet
spaces, imposing only a certain metric
estimate on the partial differentials.
As a tool, we show the $C^k$-dependence
of fixed points on parameters
for suitable families of contractions
of a Fr\'{e}chet space.
The investigations were stimulated by
a recent metric approach to differentiability
in Fr\'{e}chet spaces by Olaf M\"{u}ller.
Our results also subsume generalizations
of M\"{u}ller's Inverse Function Theorem
for mappings between Fr\'{e}chet spaces.
As an application, we study
existence, uniqueness and parameter-dependence
of solutions to suitable ordinary differential
equations in Fr\'{e}chet spaces.
\vspace{3mm}
\end{abstract}
{\footnotesize
% = differentiation theory (Gateaux, Frechet,...)
}\vspace{2.8mm}
\begin{center}
{\Large\bf Introduction}\vspace{.5mm}
\end{center}
One of the most famous and useful results of
infinite-dimensional
differential\linebreak
calculus beyond Banach spaces
is the Nash-Moser Inverse Function Theorem
(see \cite{Ham}, \cite{KaM}; cf.\
\cite{Nas}, \cite{Mos}),
which provides a smooth local inverse
under restrictive conditions in terms
of a given fundamental sequence
of seminorms (a ``grading'') on the space.
A variant of the Nash-Moser Theorem for
implicit functions is also available~\cite{Ser}.
These theorems are difficult to prove,
and also their hypotheses are usually
difficult to check in applications.
Besides these results (and some variants),
inverse and implicit function
theorems are available for mappings between
bornological spaces in the framework of ``bounded
differential calculus'' by Colombeau
(see \cite[Chapter~13]{Hog} for a survey).
Implicit functions from topological
vector spaces to Banach spaces
have been studied in various settings
of infinite-dimensional calculus
and in varying generality
(see \cite{Hil}, \cite{Tei}, \cite{IMP},
\cite{IM2}).
Furthermore, \cite{Hi3} provides results
concerning the solutions~$\phi$ to
equations $f(x,\phi(x))=0$,
where $F=\pl\,F_j$\vspace{-.3mm} is a projective limit of Banach spaces
and $f\colon E\times F\to F$
of the form $f=\pl\, f_j$\vspace{-.9mm}
for suitable maps $f_j\colon E\times F_j\to F_j$.\\[2.5mm]
Recently, Olaf M\"{u}ller
formulated a metric approach to differential
calculus for mappings between Fr\'{e}chet spaces
and provided an Inverse Function Theorem
for certain ``bounded differentiable'' maps~\cite{Mul}
(which we call ``$MC^1$-maps''
to avoid confusion with Colombeau's
venerable ``bounded differential calculus'').
M\"{u}ller does not need to introduce gradings
on~$F$ and work with
``tame'' smooth maps as in the case of the
Nash-Moser Theorem.
Rather, he equips~$F$ with a
translation invariant metric~$d$
defining its topology (in which case
$(F,d)$ is called a ``metric Fr\'{e}chet space''),
and then introduces metric concepts
which strongly depend on the choice
of~$d$.
Using~$d$ systematically,
he succeeds in adapting many familiar arguments
and results from the Banach case to
the Fr\'{e}chet case,
and obtains a simple and natural proof
of his inverse function theorem.\\[2.5mm]
M\"{u}ller asserts (see \cite[\S\,5]{Mul})
that the metric approach is general enough
to cover some of the standard applications
of the Nash-Moser Theorem
(like those given by Hamilton~\cite{Ham}).
But this claim is too optimistic,
as shown by Hiltunen~\cite{HKi}
(cf.\ also Section~\ref{nontriv} below).
Rather, the metric approach
(and its variants described
in the current article)
should be seen as a method
which yields very strong conclusions,
but only in quite restrictive situations.\footnote{See \cite{HKi}
for a class
of examples to which the method
does not (and was not intended to)
apply.
In the case of bounded metrics,
this is also clear
from our Propositions~\ref{propshape}
and~\ref{quasiisom},}
Despite its natural limitations, the method produces
valuable new results (and we shall even
encounter novel aspects of differential calculus
in Banach spaces in this article).\\[2.5mm]
One of the essential ideas of M\"{u}ller is to
replace
the (unwieldy) space $\cL(E,F)$ of all continuous
linear operators between metric
Fr\'{e}chet spaces $(E,d)$
and $(F,d')$
by the space\footnote{We use the notational
conventions
of the present article here, rather than those from~\cite{Mul}.
Also, we tacitly assume that $d$ and $d'$
have absolutely convex balls.}
\[
\cL_{d,d'}(E,F)
\]
of all linear maps
from~$E$ to~$F$ which
are (globally) Lipschitz continuous
as mappings between the metric
spaces $(E,d)$ and $(F,d')$.
Then $\cL_{d,d'}(E,F)$ is a vector
space, and also a
topological group under addition
with respect to the topology
defined by the complete metric
$(A,B)\mto \|A-B \|_{d,d'}$,
where
\[
\|A \|_{d,d'}\;:=\; \sup_{x\in E\setminus\{0\}}\,
\frac{d'(A(x),0)}{d(x,0)}
\]
is the (minimal) Lipschitz constant $\Lip(A)$
of $A \in \cL_{d,d'}(E,F)$.
The spaces $\cL_{d,d'}(E,F)$
have good properties which
would be impossible for $\cL(E,F)$ (cf.\ \cite{Mai}):
For example, the evaluation map
$\cL_{d,d'}(E,F)\times E\to F$ is continuous,
and $\cL_d(E):=\cL_{d,d}(E,E)$ is a topological
ring with open unit group~$\cL_d(E)^\times$
and continuous inversion
$\cL_d(E)^\times \to \cL_d(E)^\times$, $A \mto A^{-1}$.\\[2.5mm]
In the present article, we
combine M\"{u}ller's ideas with
the approach to implicit functions
from topological vector spaces to Banach spaces
developed in~\cite{IM2}.
In contrast to M\"{u}ller,
we formulate all of our results
in a standard setting of differential
calculus:
Our $C^k_\K$-maps are $C^k$-maps over $\K\in \{\R,\C\}$
in the sense of Michal and Bastiani
(also known as Keller's $C^k_c$-maps).\footnote{See \cite{Bas}, \cite{RES},
\cite{GaN}, \cite{Ham}, \cite{Mil}
for discussions of such maps,
in varying generality.}
These are the maps
widely used as the basis of
infinite-dimensional Lie theory
(except for the literature
based on the convenient differential calculus
as in~\cite{KaM}).
By contrast, M\"{u}ller uses
``bounded differentiable'' maps ($MC^1$-maps):
these are $C^1$-maps $f\colon U\to F$
from an open subset $U\sub E$
of a metric Fr\'{e}chet space $(E,d)$
to a metric Fr\'{e}chet space $(F,d')$
such that $f'(U)\sub \cL_{d,d'}(E,F)$
and $f'\colon U\to\cL_{d,d'}(E,F)$
is continuous (where $f'(x)\colon E\to F$
is the differential of~$f$ at~$x$).
For our results, this continuity property is
not required, and this is a real
advantage because the class of $MC^1$-maps
can be quite small in some cases
(see Remark~\ref{proverpr}).\\[3mm]
Among our main results is the following
Implicit Function Theorem for Keller $C^k_c$-maps
from arbitrary topological vector spaces to Fr\'{e}chet
spaces.\\[3.5mm]
%
%
%
%\ma{genimp}
{\bf Theorem~A (Generalized Implicit Function Theorem).}
\emph{Let $\K\in \{\R,\C\}$,
$E$ be a topological $\K$-vector space,
$F$ be a Fr\'{e}chet space over~$\K$,
and $f\colon U\times V\to F$ be a $C^k_\K$-map,
where $U\sub E$
and $V\sub F$ are open sets.
Given $x\in U$, abbreviate $f_x:=f(x,\sbull)
\colon V\to F$.
Assume that $f(x_0,y_0)=0$ for some $(x_0,y_0)\in U\times V$
and that $f_{x_0}'(y_0)\colon F\to F$
is invertible.
Furthermore, assume that there
exists a translation-invariant metric~$d$
on~$F$ defining its topology
such that all $d$-balls are absolutely convex
and}
%
%\ma{ensurenice}
\begin{equation}\label{ensurenice}
\sup_{(x,y)\in U\times V}\, \|\id_F-f'_{x_0}(y_0)^{-1}f_x'(y)\|_{d,d}
\; <\; 1\,.
\end{equation}
\emph{Then there exist open neighborhoods
$U_0\sub U$ of~$x_0$ and $V_0\sub V$ of~$y_0$
such that}
\[
\{(x,y)\in U_0\times V_0\colon f(x,y)=0\}
\; =\; \graph\lambda
\]
\emph{for a $C^k_\K$-map $\lambda \colon U_0\to V_0$.}\\[4mm]
Note that Theorem~A also covers the case of
complex analytic maps (in the usual sense,
as in \cite{BaS})
because a map from an open subset
of a complex topological vector space
to a complex locally convex space
is $C^\infty_\C$
if and only if it is complex analytic
(see \cite[Propositions~7.4 and~7.7]{Ber}).
We can also deal with local inverses.\\[4mm]
%
%
%\ma{genimp}
{\bf Theorem~B (Local inverses for {\boldmath$C^k$}-maps
between Fr\'{e}chet spaces).}
\emph{Let $F$ be a Fr\'{e}chet space
over~$\K \in \{\R,\C\}$
and $f\colon U \to F$ be a $C^k_\K$-map
on an open subset
$U\sub F$, where $k\in \N\cup\{\infty\}$.
Let $x_0\in U$.
If $f'(x_0)\colon F\to F$
is invertible
and there exists a translation-invariant metric~$d$
on~$F$ defining its topology
such that all $d$-balls are absolutely convex
and}
\begin{equation}\label{latermoregen}
\sup_{x\in U} \|\id_F-f'(x_0)^{-1}f'(x)\|_{d,d}
\; <\; 1\,,
\end{equation}
\emph{then there exists an open neighborhood
$U_0\sub U$ of~$x_0$ such that
$f(U_0)$ is open in~$F$ and
$f|_{U_0}\colon U_0\to f(U_0)$
is a $C^k_\K$-diffeomorphism.}\\[4mm]
We remark that, slightly more generally,
(\ref{latermoregen})
can be replaced by the following condition:
There exist isomorphisms
$S,A,T\colon F\to F$
of topological vector spaces
such that $S\circ A \circ T\in \cL_d(F)^\times$
and
%\ma{nowmoregen}
\begin{equation}\label{nowmoregen}
\sup_{x\in U} \|S\circ (A-f'(x))\circ T\|_{d,d}
\; <\; \frac{1}{\|(S\circ A \circ T)^{-1}\|_{d,d}}\,.
\end{equation}
Likewise, (\ref{ensurenice})
can be replaced by the condition:
There exist isomorphisms
of topological vector spaces $S,A,T\colon F\to F$
such that $S\circ A \circ T\in \cL_d(F)^\times$
and
%\ma{nowmoregen2}
\begin{equation}\label{nowmoregen2}
\sup_{(x,y)\in U\times V}
\|S\circ (A-f'_x(y))\circ T\|_{d,d}
\; <\; \frac{1}{\|(S\circ A \circ T)^{-1}\|_{d,d}}\,.
\end{equation}
Both Theorem~A and~B will be deduced from
a suitable ``Inverse Function Theorem with Parameters''
(Theorem~\ref{advif}),
dealing with families of local diffeomorphisms.
This theorem is our
main result
(whence we should count it as Theorem~C,
although we shall not restate it here
in the introduction).
As a technical tool, in Section~\ref{secfp}
we prove $C^k$-dependence
of fixed points on
parameters,
for certain ``uniform families of special contractions''
(as in Definition~\ref{defunicon}
below):\\[4mm]
{\bf Theorem D (Dependence of Fixed
Points on Parameters).}
\emph{Let
$(F,d)$ be a metric Fr\'{e}chet space over~$\K\in \{\R,\C\}$
with absolutely convex balls,
and
$E$ be a topological $\K$-vector space.
Let $P\sub E$
and $U\sub F$ be open sets,
and
$f\colon P \times U\to F$
be a continuous map such that
$f_p:=f(p,\sbull)\colon U\to F$
defines a uniform family
$(f_p)_{p\in P}$
of contractions.
Then the following holds:}
\begin{itemize}
\item[\rm (a)]
\emph{The set $Q$ of all $p\in P$ such that
$f_p$ has a fixed point $x_p$
is open in~$P$.
Furthermore, the map $\phi\colon Q \to U$, $\phi(p):=x_p$
is continuous.}
\item[\rm (b)]
\emph{If $f$ is $C^k_\K$
for some $k\in \N \cup\{\infty\}$
and $(f_p)_{p\in P}$
is a uniform family of special contractions,
then also~$\phi$
is $C^k_\K$.}\vspace{1.3mm}
\end{itemize}
{\bf Variants for non-open domains.}
We mention that, if~$E$ is locally convex,
then Theorem~\hspace*{-.1mm}A\hspace*{.3mm}
and Theorem~\hspace*{-.1mm}D\hspace*{.3mm} hold just as
well if $U\sub E$ (resp., $P\sub E$)
is a locally convex subset
with dense interior. Our proofs
also cover these variants.\\[3mm]
{\bf The case of mappings into Banach spaces.}
In the case of mappings into Banach spaces,
we recover the inverse function theorem
with parameters and the theorem
on implicit functions from
topological vector spaces
to real or complex Banach spaces
from~\cite{IM2}.
%(see Section~\ref{secBanach}).
The proofs of Theorem~A
and Theorem~D are direct adaptations
of the proofs in~\cite{IM2}.\\[3mm]
{\bf Applications to ODEs in Fr\'{e}chet spaces.}
In Section~\ref{secODE},
we prove existence, uniqueness
and $C^k$-dependence on parameters
for $C^k$-solutions to
suitable ordinary differential equations
in Fr\'{e}chet spaces.
Our results (recorded as Theorem~\ref{exandun}) are
slightly more general than the following.\\[3mm]
{\bf Theorem~E (Existence and Uniqueness Theorem for ODEs
in Fr\'{e}chet Spaces).}
\emph{Let
$(F,d)$ be a metric Fr\'{e}chet space
over~$\R$, with absolutely convex balls.
Let $k\in \N_0\cup\{\infty\}$,
$J \sub \R$ be an interval,
$E$ be a locally convex space
and $P\sub E$ as well as $U\sub F$ be open subsets.
Let $f\colon J \times U\times P \to F$ be a $C^k_\R$-map
which satisfies a local special
contraction condition {\rm(SCC)} in its second argument
$($as in Definition~{\rm \ref{deflocSCC}}$)$.
Let $t_0\in J$, $x_0\in U$ and $p_0\in P$.
Then there exist
open neighborhoods
$U_1\sub U$ of~$x_0$,
$P_1\sub P$ of~$p_0$
and $r>0$
such that
for all $(x_1,p_1)\in U_1\times P_1$
and $t_1\in J_1:= \;]t_0-r,t_0+r[\;\cap \, J$,
the initial value problem}
\begin{equation}\label{initval}
x'(t)\;=\; f(t,x(t),p_1)\,,\qquad
x'(t_1)\;=\; x_1
\end{equation}
\emph{has a unique $C^k_\R$-solution
$\phi_{t_1,x_1,p_1}\colon J_1 \to U$
and also the following map is~$C_\R^k$}:
\[
\Psi \colon J_1\times J_1 \times U_1\times P_1 \to U\,,
\quad \Psi(t_1,t,x_1,p_1):=\phi_{t_1,x_1,p_1}(t)\,.\vspace{.7mm}
\]
To prove Theorem~E,
we use Theorem~A and a Lipschitz version thereof
(Corollary~\ref{lipimpl}),
combined with some
preparatory results concerning
differentiability properties
of pushforwards depending on
parameters provided in Section~\ref{sec-push-it}.
We remark that, if $(F,\|.\|)$
is a Banach space and $d(x,y)=\|x-y\|$,
then the local SCC can be replaced with the ordinary
local Lipschitz condition (in the middle
argument) for~$f$.
If also the space of parameters~$E$
is a Banach space, then
an analogue of
Theorem~E for $k$ times continuously
Fr\'{e}chet differentiable maps
($FC^k$-maps)
is known
(see the classical literature
or also \cite[\S3.1, Theorem~1.1]{CaH},
where $F$ is assumed finite-dimensional
and $k\geq 1$).
But for Keller $C^k_c$-maps,
the result is new even in the
Banach case.\\[3mm]
{\bf Global Inverse Function Theorems for Fr\'{e}chet Spaces.}
Beyond the standard theorems on local
inverses,
there is Hadamard's Global
Inverse Function Theorem
for continuously Fr\'{e}chet-differentiable
self-maps of a Banach space
(see \cite[Chapter II.C, \S\,4, Theorem~1]{CWD},
\cite[Chapter~2, Theorem~3.9]{CaH},
or \cite[Theorem~6.2.4]{KaP} for a
more restricted version).
In Section~\ref{secglob},
we prove analogous global inverse
function theorems for self-maps
of a Fr\'{e}chet space,
both for $C^k$-maps
(Theorem~\ref{globCk})
and $MC^k$-maps (Theorem~\ref{globMC}).\\[3mm]
{\bf Further variations.}
In~\cite{Mul}, one also finds
a discussion of left and right inverses.
Along the lines of Theorem~C and its proof,
one could use Theorem~D
also to prove parameter-dependent
versions of these one-sided
inverse function theorems,
providing left (resp., right) inverses
depending on a parameter in a general topological
vector space. However, we refrain
from doing so here and prefer to concentrate
on the central results.\\[3mm]
{\bf Prospects.} As the next stage,
it would be interesting to study examples
and to explore the scope of the approach.
For example, it might go along well
with certain topologically nilpotent
Fr\'{e}chet
Lie algebras and corresponding Lie groups.\\[3mm]
{\bf Complications of metric differential calculus.}
Let us mention in closing
that a problem has been overlooked in~\cite{Mul}:
Contrary to claims made there (after
\cite[Definition~3.13]{Mul}),
$\cL_{d,d'}(E,F)$
is not always
a Fr\'{e}chet space.
In fact,
examples show that
$\cL_{d,d'}(E,F)$ is, in general,
not a topological vector space
because balls around~$0$ are not absorbing
(see Proposition~\ref{vectgp}).
It merely is a locally convex vector group
in the sense of Ra\u{\i}kov
(as in~\cite{Rai}, also~\cite{Aus}).
Fortunately, this does not endanger the use
of higher order differentiability properties
in \cite{Mul} (as clarified in Remark~\ref{commt4}).
It implies, however, that the
class of $MC^1$-maps is quite small
in many typical situations
(see Remark~\ref{proverpr}).\\[2.5mm]
Another comment concerns the
type of metrics used by M\"{u}ller: These
are somewhat problematic, because they
need not have convex balls
(see Remark~\ref{commt1}).
By contrast, we prefer to use
metrics with absolutely convex balls.
%
%
%
%
%
%
%
%
%
%
%
%\ma{secprel}
\section{Preliminaries and basic facts}\label{secprel}
In this section, we set up our notation
and terminology concerning
differential calculus in
infinite-dimensional spaces
and mappings between Fr\'{e}chet spaces.\\[2.5mm]
Throughout the article, $\K\in \{\R,\C\}$.
All topological vector spaces
and all topological groups are assumed Hausdorff.
Our basic terminology concerning locally
convex spaces follows~\cite{Rud}.
We write $\N:=\{1,2,\ldots\}$
and $\N_0:=\N\cup \{0\}$.
\subsection*{{\normalsize Prerequisites concerning {\boldmath$C^k$}-maps}}
Naturally,
we are mainly interested
in results concerning mappings
from open subsets of real or complex
locally convex spaces to Fr\'{e}chet
spaces.
%and the reader
%is invited to concentrate on this case
%if he or she so desires.
However, most of the results
(and their proofs)
apply just as well
to mappings on open subsets
of non-locally convex spaces,
and also to mappings on
suitable
subsets with dense interior.
Since mappings on non-open
sets are useful and frequently
encountered
in infinite-dimensional
analysis and Lie theory,
we present our results in full
generality.
%Also in the current article,
This is also vital
for our main application:
The approach to
ODEs in Fr\'{e}chet spaces
in Section~\ref{secODE}
hinges on the consideration of $C^k$-maps
on sets of the form $[0,1]\times U$, with $U$
an open subset of a locally
convex space.\\[2.5mm]
The exact framework
of differential calculus
will be described now.
\begin{defn}\label{C11}
Given $\K\in \{\R,\C\}$,
let $E$ be a topological $\K$-vector
space and $F$ be a locally convex
topological $\K$-vector space.
If~$E$ is not locally convex, let
$U\sub E$ be an open set.
If $E$ is a locally convex space,
then more generally let $U\sub E$ be a subset with dense
interior which is locally convex
in the sense that each $x\in U$ has
a convex neighborhood $V\sub U$
(and hence arbitrarily small convex neighborhoods).
Let $f\colon U\to F$ be a map.
The map~$f$ is called $C^0_\K$ if it is continuous.
The map~$f$ is called $C^1_\K$ if it is continuous
and there exists a (necessarily unique)
continuous map $df\colon U\times E\to F$
such that
\[
df(x,y)\;=\;\lim_{t\to 0}\frac{f(x+ty)-f(x)}{t}
\]
for all $x$ in the interior~$U^0$ of~$U$
and all $y\in E$
(with $0\not=t\in \K$ sufficiently small).
Since $U\times E$ is open in $E\times E$
(resp., a locally convex subset with
dense interior),
we can proceed by induction:
Given $k\in \N$, we say that
$f$ is $C^{k+1}_\K$ if $f$ is
$C^1_\K$ and $df\colon U\times E\to F$
is~$C^k_\K$.
We say that~$f$ is $C^\infty_\K$
if~$f$ is~$C^k_\K$ for each $k\in \N_0$.
If~$\K$ is understood,
we simply write $C^k$ instead of $C^k_\K$,
for $k\in \N_0\cup\{\infty\}$.
\end{defn}
If $f\colon E\supseteq U\to F$
is $C^1_\K$, then $f'(x):=df(x,\sbull)\colon E\to F$
is a continuous $\K$-linear map
(cf.\
\cite[Chapter~1]{GaN}
and \cite[Lemma~1.9]{RES}).\\[2.5mm]
If $U\sub \K$,
we shall occasionally write $f'(x)$ also
for $f'(x)(1)=\frac{d}{dx}f(x)$,
in particular when dealing with solutions
to differential equations.
It will always be clear from the context
which meaning of~$f'(x)$ is intended.\\[2.5mm]
At some places, we use an alternative
approach to~$C^1_\K$-maps
based on
continuous extensions~$f^{[1]}$ of directional
difference quotients,
which even remains meaningful
for mappings
into non-locally convex
spaces.\footnote{As introduced in~\cite{Ber}
for maps on open sets and in~\cite{IM2}
for maps on sets with dense interior.}
This
alternative approach
is not an unnecessary ballast,
but invaluable for
our purposes, because the proof
of our main technical result
(Lemma~\ref{prepdep2})
is best formulated
in terms of the maps~$f^{[1]}$.
\begin{defn}\label{C12}
Let $E$ and~$F$ be topological $\K$-vector
spaces over $\K\in \{\R,\C\}$
and $U\sub E$ be a subset with dense
interior.
Given a map
$f\colon U\to F$, its directional
difference quotients
%\ma{defopf}
\begin{equation}\label{defopf}
f^{]1[}(x,y,t)\; :=\; \frac{f(x+ty)-f(x)}{t}
\end{equation}
make sense for all
$(x,y,t)\in U\times E\times \K^\times$
such that $x+ty\in U$.
Allowing now also the value $t=0$, we define
\[
U^{[1]}\; :=\; \{
(x,y,t)\in U\times E\times \K\colon
x+ty\in U\}
\]
and say that $f\colon U\to F$ is $C^1_\K$
if $f$ is continuous and
there exists a (necessarily unique)
continuous map
\[
f^{[1]}\colon U^{[1]}\to F
\]
which extends the difference quotient map,
i.e., $f^{[1]}(x,y,t)=f^{]1[}(x,y,t)$
for all $(x,y,t)\in U^{[1]}$ such that $t\not=0$.
\end{defn}
\begin{rem}
If $f$ is $C^1_\K$ in the sense
of Definition~\ref{C12},
then the mapping\linebreak
$df\colon U\times E\to F$, $df(x,y):=f^{[1]}(x,y,t)$
is continuous
and the differential
$f'(x):=df(x,\sbull)\colon E\to F$
is continuous linear,
for each $x\in U$
(cf.\ \cite[Proposition~2.2]{Ber}).
\end{rem}
We mention that no ambiguity occurs
because if $E,F$ and $U\sub E$
happen to satisfy the hypotheses
of both Definition~\ref{C11} and Definition~\ref{C12},
then a map
$f\colon E\supseteq U\to F$
is $C^1_\K$ in the sense of
Definition~\ref{C11}
if and only if it is $C^1_\K$
in the sense of Definition~\ref{C12}
(cf. \cite[Proposition~7.4]{Ber}
or \cite[Chapter~1]{GaN}).
\begin{numba}\label{chainr}
We need two versions of the Chain Rule
(cf.\ \cite[Proposition~1.15]{RES},
\cite[Chapter~1]{GaN}
and \cite[Proposition~3.1]{Ber}):
\begin{itemize}
\item[\rm (a)]
If $E$, $F$ and $H$ are topological
$\K$-vector spaces, $U\sub E$ and $V\sub F$ are
subsets with dense interior, and $f\colon  U\to V\sub F$,
$g\colon  V\to H$ are $C^1_\K$-maps,
then also the composition
$g\circ f\colon  U\to H$ is~$C^1_\K$,
and $(g\circ f)'(x)=g'(f(x))\circ f'(x)$
for all $x\in U$.
\item[\rm (b)]
Let $E$ be a topological
$\K$-vector space and $F$ as well
as $H$ be locally convex topological
$\K$-vector spaces.
Let $U\sub E$ be open
(if~$E$ is not locally convex)
or a locally convex subset with
dense interior
(if~$E$ is locally convex).
Let $V\sub F$ be a locally convex subset
with dense interior.
If $k\in \N_0\cup\{\infty\}$
and both $f\colon  U\to V\sub F$
and $g\colon  V\to H$ are $C^k_\K$-maps,
then also their composition
$g\circ f\colon  U\to H$ is~$C^k_\K$.
\end{itemize}
\end{numba}
Given a linear map $A\colon E\to F$ between
vector spaces, we shall frequently write~$A.x$
instead of~$A(x)$.
\subsection*{{\normalsize Metric Fr\'{e}chet spaces
and linear, Lipschitz maps}}
Given a metric space $(X,d)$,
we write $\wb{B}_r^d(x):=\{y\in X\colon d(x,y)\leq r\}$
for $x\in X$ and $r\in [0,\infty[$
and $B_r^d(x):=\{y\in X\colon d(x,y)<r\}$
if $r>0$. If $d$ or~$X$ is understood,
we also write $B_r^X(x)$ for $B_r^d(x)$,
or simply $B_r(x)$.
Likewise for $\wb{B}_r^d(x)$.
\begin{numba}
A \emph{metric Fr\'{e}chet space}
is a Fr\'{e}chet space~$F$,
equipped with a metric\linebreak
$d\colon F\times F\to [0,\infty[$
defining its topology which is translation
invariant, i.e., $d(x+z,y+z)=d(x,y)$
for all $x,y,z\in F$.
In this case, we define
$\|x\|_d:=d(x,0)$ for $x\in F$
and note that $d$ can be recovered from
$\|.\|_d\colon F\to [0,\infty[$ via
$d(x,y)=\|x-y\|_d$.
Recall that a $0$-neighborhood
$U\sub F$ is \emph{absolutely convex}
if it is convex and \emph{balanced},
i.e., $\wb{B}^\K_1(0)U\sub U$.
We say that \emph{$d$ has symmetric} (\emph{resp., balanced,
resp., convex,
resp., absolutely convex})
\emph{balls}
if $\wb{B}_r^d(0)=-\wb{B}_r^d(0)$
(resp., $\wb{B}_r^d(0)$
is balanced, resp., it is convex, resp., absolutely
convex) for each $r\geq 0$.
Then $B_r^d(0)$
has analogous properties,
for each $r>0$.
\end{numba}
\begin{ex}\label{exam1}
Every Fr\'{e}chet space~$F$
admits a translation invariant
metric~$d$ which has absolutely
convex balls and defines the topology of~$F$.
In fact, pick
any sequence $w=(w_n)_{n\in \N}$ of
real numbers~$w_n>0$
such that $\lim_{n\to\infty}\,w_n=0$,
and any sequence
$p=(p_n)_{n\in \N}$
of continuous seminorms
$p_n\colon F\to[0,\infty[$
which define the topology of~$F$ in the sense
that finite intersections of
sets of the form $p_n^{-1}([0,\ve[)$
(with $n\in \N$, $\ve>0$)
form a basis of $0$-neighborhoods
in~$F$.
Then
\[
d_{w,p}\colon F\times F\to[0,\infty[\,,\quad
d_{w,p}(x,y)\, :=\, \sup_{n\in \N}w_n\,\frac{p_n(x-y)}{1+p_n(x-y)}
\]
is a metric with the desired properties.
\end{ex}
Metrics of the form $d_{w,p}$ (as just defined)
will occasionally be called \emph{standard metrics}
in the following.
%
%
%\ma{splobsv}
\begin{la}\label{splobsv}
Let $(F,d)$ be a metric Fr\'{e}chet
space with balanced balls,
$t\in \K$ and $x\in F$.
Then $\|tx\|_d\leq \|x\|_d$ if $|t|\leq 1$;
$\|tx\|_d=\|x\|_d$ if $|t|=1$;
and
$\|tx\|_d\leq 2|t|\cdot \|x\|_d$
if $|t|\geq 1$.
In any case,
%\ma{sumup}
\begin{equation}\label{sumup}
\|tx\|_d\; \leq \; \max\{1,2|t|\}\, \|x\|_d\,.
\end{equation}
\end{la}
\begin{proof}
The first assertion is clear
since $\wb{B}^d_{\|x\|_d}(0)$ is a
balanced $0$-neighborhood.
The second assertion follows
from the first and the observation that also
$\|x\|_d=\|t^{-1}(tx)\|_d\leq \|tx\|_d$
by the first assertion,
if $|t|=1$.
If $|t|\geq 1$,
set $n:=[\, |t|\, ]+1 \geq |t|$, using the Gau\ss{} bracket
(integer part).
Then $\|tx\|_d\leq \|nx\|_d\leq n\|x\|_d\leq 2|t|\cdot \|x\|_d$.
\end{proof}
%
%
%
%\ma{defnspac}
\begin{defn}\label{defnspac}
Given metric Fr\'{e}chet spaces
$(E,d)$ and $(F,d')$,
we let $\cL_{d,d'}(E,F)$
be the set of all linear maps
$A\colon E\to F$ such that
%
%\ma{condlip}
\begin{equation}\label{condlip}
\|A\|_{d,d'}\; :=\;
\sup_{x\in E\setminus\{0\}} \frac{\|A.x\|_{d'}}{\|x\|_d}
\; <\, \infty\,.
\end{equation}
We abbreviate $\cL_d(E):=\cL_{d,d}(E,E)$;
occasionally, we write
$\|A\|_d:=\|A\|_{d,d}$ for $A\in\cL_d(E)$
(as there is little risk of confusion with $\|x\|_d:=d(x,0)$
for $x\in E$).
\end{defn}
Condition~(\ref{condlip})
means that~$A$ is
Lipschitz continuous
as a map $(E,d)\to (F,d')$.
To prevent misunderstandings,
let us mention that although $f'(x)$
denotes the differential of~$f$,
we shall frequently use the symbol~$d'$
in a different meaning
(it denotes a metric on the range space of a map).
%
%
%\ma{simplestpr}
\begin{rem}\label{simplestpr}
The following simple properties
of $\cL_{d,d'}(E,F)$ and the functions
$\|.\|_{d,d'}$
will be used later.
\begin{itemize}
\item[(a)]
In the situation of Definition~\ref{defnspac},
%
%\ma{standest1}
\begin{equation}\label{standest1}
\|A.x\|_{d'}\;\leq\; \|A\|_{d,d'}\|x\|_d\quad
\mbox{for all $x\in E$,}
\end{equation}
as is clear from the definition of $\|.\|_{d,d'}$.
Furthermore, $0\in \cL_{d,d'}(E,F)$ with $\|0\|_{d,d'}=0$
and
%
%\ma{thusnom}
\begin{equation}\label{thusnom}
\|A\|_{d,d'}\;>\; 0\quad \mbox{if $\,A\in\cL_{d,d'}(E,F)\setminus\{0\}$,}
\end{equation}
because there is a $x\in E$ with $A.x\not=0$
and thus $\|A\|_{d,d'}\geq \frac{\|A.x\|_{d'}}{\|x\|_d}>0$.
\item[(b)]
If also $(G,d'')$ is a metric
Fr\'{e}chet space, then
%
%\ma{standest2}
\begin{equation}\label{standest2}
\|B \circ A\|_{d,d''} \leq
\|B\|_{d',d''} \|A\|_{d,d'}\;
\mbox{for $A\!\in \!\cL_{d,d'}(E,F)$, $B\!\in \!\cL_{d',d''}(F,G)$,}
\end{equation}
as an immediate consequence of~(\ref{standest1}).
\item[(c)]
If $A,B\in \cL_{d,d'}(E,F)$,
then also $A+B\in \cL_{d,d'}(E,F)$
and
%\ma{subaff}
\begin{equation}\label{subaff}
\|A+B\|_{d,d'}\;\leq\; \|A\|_{d,d'}+\|B\|_{d,d'}\;<\;\infty\,,
\end{equation}
because $\frac{\|(A+B).x\|_{d'}}{\|x\|_d}
\leq \frac{\|A.x\|_{d'}}{\|x\|_d}+
\frac{\|B.x\|_{d'}}{\|x\|_d}\leq \|A\|_{d,d'}+\|B\|_{d,d'}$
for all $x\in E\setminus\{0\}$.
Thus $\cL_{d,d'}(E,F)$ is a monoid
under addition.
\item[(d)]
If $d'$ or $d$ has symmetric balls,
then $\frac{\|-Ax\|_{d'}}{\|x\|_d}=\frac{\|A.x\|_{d'}}{\|x\|_d}$
(resp., $\frac{\|-Ax\|_{d'}}{\|x\|_d}=\frac{\|A.(-x)\|_{d'}}{\|-x\|_d}$),
entailing that $-A\in \cL_{d,d'}(E,F)$ and
%
%\ma{trivvv}
\begin{equation}\label{trivvv}
\|{-A}\|_{d,d'}\;=\; \|A\|_{d,d'}\,,\quad\mbox{for each
$A\in \cL_{d,d'} (E,F)$.}
\end{equation}
Hence $\cL_{d,d'}(E,F)$ is a
subgroup of $F^E$ in this case,
and it follows from~(a) and~(c)
that
%
%\ma{transmetr}
\begin{equation}\label{transmetr}
D_{d,d'}\colon \cL_{d,d'}(E,F)\times \cL_{d,d'}(E,F)\to [0,\infty[\,,\quad
(A,B)\mto \|A-B\|_{d,d'}
\end{equation}
is a translation invariant
metric on the abelian group $\cL_{d,d'}(E,F)$
which defines a topology on
$\cL_{d,d'}(E,F)$
turning the latter into a
topological group.
We shall always equip $\cL_{d,d'}(E,F)$
with the metric $D_{d,d'}$
and the corresponding topology.
\end{itemize}
\end{rem}
It is essential to have estimates on
the size of integrals.
%
%
%
%
%
%\ma{estinteg}
\begin{la}\label{estinteg}
Let $(F,d)$ be a metric Fr\'{e}chet
space with convex balls,
and\linebreak
$\gamma\colon [0,1]\to F$
be a continuous curve.
Then
\begin{equation}\label{heregood}
\left\|\int_0^1 \gamma(t)\;dt\right\|_d\;\leq\;
\max_{t\in [0,1]}\,\|\gamma(t)\|_d\,.
\end{equation}
\end{la}
\begin{proof}
Set $r:=\max_{t\in [0,1]}\,\|\gamma(t)\|_d$.
The ball $B:=\wb{B}_r^d(0)$
is convex
and contains $\gamma(t)$ for each $t\in [0,1]$.
Each Riemann sum of~$\gamma$
is a convex combination of values
of~$\gamma$, whence it lies in~$B$.
Since~$B$ is closed,
it follows that also the limit $\int_0^1\gamma(t)\,dt$
of the Riemann sums lies in~$B$.
\end{proof}
We record a variant of
\cite[Proposition~3.18]{Mul}:
%
%\ma{corrfalse}
\begin{la}\label{corrfalse}
Let $(E,d_E)$ and $(F,d_F)$
be metric Fr\'{e}chet spaces
such that~$d_F$ has absolutely convex balls.
Let $U\sub E$ be a convex
subset with non-empty interior
and $f\colon U\to F$ be a $C^1$-map.
Then
%
%\ma{heregood2}
\begin{equation}\label{heregood2}
\|f(y)-f(x)\|_{d_F}\;\leq\;
\|y-x\|_{d_E}\cdot \sup_{t\in [0,1]} \|f'(x+t(y-x))\|_{d_E,d_F}
\;\; \mbox{for all $x,y\in U$.}
\end{equation}
\end{la}
\begin{proof}
Apply (\ref{heregood})
to $\gamma\colon [0,1]\to F$,
$\gamma(t)=f'(x+t(y-x)).(y-x)$
with $\int_0^1\gamma'(t)\, dt=f(y)-f(x)$
and use (\ref{standest1})
to estimate $\|f'(x+t(y-x)).(y-x)\|_{d_F}$.
\end{proof}
\begin{rem}\label{commt1}
We warn the reader that,
in the situation of
Example~\ref{exam1},
the metric $D$ on~$F$ given by $D(x,y):=\sum_{n=1}^\infty
2^{-n}\frac{p_n(x,y)}{1+p_n(x,y)}$
does not have convex balls
in general.
For instance, $\R^\N$ with $D(x,y):=\sum_{n=1}^\infty
2^{-n}\frac{|x_n-y_n|}{1+|x_n-y_n|}$ does
not have convex balls.
To see this, let $e_n:=(0,\ldots, 0,1,0,\ldots)\in \R^\N$
with~$1$ only in the $n$-th slot.
Then $v_1:=10 e_1$ and $v_2:=10 e_2$
are elements\vspace{-.5mm} of the ball $\wb{B}_\frac{1}{2}^D(0)$,
but $\frac{1}{2}v_1+\frac{1}{2}v_2\not\in \wb{B}_{\frac{1}{2}}^D(0)$
because
$\|\frac{1}{2}v_1+\frac{1}{2}v_2\|_D=\frac{5}{8}>\frac{1}{2}$.\\[2.5mm]
Since not all of the $D$-balls
are convex, it is not clear
whether Formula~(\ref{heregood2})
also holds if the metric~$D$ is used.
The contrary is claimed
in \cite[Proposition~3.18]{Mul},
but the author cannot make sense
of its proof.\footnote{No clear explanation
is given for the first inequality in the proof
of \cite[Proposition~3.18]{Mul}.}
This may be a serious problem
for~\cite{Mul},
because sums of metrics
are used in the main results of
that paper.\\[2.5mm]
For many of our main results (outside Section~\ref{nontriv}),
the crucial point is the validity of Lemma~\ref{estinteg},
rather than the absolute convexity of balls.
As long as balls are balanced,
the validity of the lemma should suffice
to carry out the proofs.
\end{rem}
\section{The space
{\boldmath $\cL_{d,d'}(E,F)$} and linear contractions}\label{nontriv}
In our studies,
linear contractions $A\colon F\to F$
of a metric Fr\'{e}chet
space $(F,d)$ will play
an important role,
i.e., mappings $A\in \cL_d(F)$
such that $\|A\|_{d,d}<1$.
It is therefore useful
to know how
linear contractions look like,
and moreover
elements in $\cL_d(F)$ close to~$0$.
With this motivation,
in the current section
we
discuss the groups
$\cL_{d,d'}(E,F)$ and the associated group
norms $\|.\|_{d,d'}$ in more detail.
In particular, we shall see that
although $\cL_{d,d'}(E,F)$
is a vector space if~$d'$ has absolutely
convex balls,
it frequently happens
that $\cL_{d,d'}(E,F)$
is \emph{not} a topological
vector space.
We also discuss various examples
which illustrate
the concept of a linear contraction,
and hint towards the limitations
of the theory.\\[2.5mm]
Although neither the vector space
structure on $\cL_{d,d'}(E,F)$
nor other results of this section
will be used later,\footnote{Except for the definition
of $MC^k$-maps.}
they seem indispensable
for a deeper \mbox{understanding.}\\[3mm]
The following proposition
slightly expands \cite[Theorem~4.2]{Mul}.
We relegate its
simple proof to Appendix~\ref{appA}.
%
%
%\ma{essLspa1}
\begin{prop}\label{essLspa1}
Let $(E,d)$ and $(F,d')$
be metric Fr\'{e}chet spaces
such that all $d'$-balls are absolutely convex.
Then the following holds:
\begin{itemize}
\item[\rm(a)]
$\cL_{d,d'}(E,F)$ is a vector
subspace of $F^E$.
\item[\rm(b)]
The evaluation map
$\cL_{d,d'}(E,F)\times E\to F$,
$(A,x)\mto A.x$ is continuous bilinear.
\item[\rm(c)]
If also $(G,d'')$ is a metric
Fr\'{e}chet space with absolutely convex balls,
then the composition map
\[
\cL_{d',d''}(F,G)\times \cL_{d,d'}(E,F)\to \cL_{d,d''}(E,G),\;\;\;
(A,B)\mto A\circ B
\]
is continuous bilinear.
\item[\rm(d)]
$D_{d,d'}$
from {\rm(\ref{transmetr})}
is a complete metric,
and has absolutely convex balls.
\item[\rm(e)]
$\cL_d(E)$ is a unital
associative $\K$-algebra,
and the topology defined by $D_{d,d}$
turns $\cL_d(E)$ into a
topological ring.
\item[\rm(f)]
The group of units
$\{ A\in \cL_d(E) \colon (\exists B \in \cL_d(E))
\;
B\circ A
=A\circ B = \id_E \}$
$=:\cL_d(E)^\times$
is open in $\cL_d(E)$
and the inversion map
$\iota\colon
\cL_d(E)^\times \to
\cL_d(E)^\times$, $A\mto A^{-1}$
is continuous.\Punkt
\end{itemize}
\end{prop}
We recall that a \emph{locally convex vector
group} is $\K$-vector
space~$E$, equipped with
a topology making $(E,+)$
a topological group
and such that~$0$ has a basis
of absolutely convex neighborhoods
(see \cite{Rai}, also~\cite[\S\,9]{Aus}).
Unlike the case of a topological
vector space, $0$-neighborhoods in~$E$
need not be absorbing.
Quite surprisingly, we have:
%
%
%\ma{vectgp}
\begin{prop}\label{vectgp}
In the situation of Proposition~{\rm\ref{essLspa1}},
$\cL_{d,d'}(E,F)$
is a locally convex vector group.
In some cases,
$\cL_{d,d'}(E,F)$
is not a topological
vector space. It can even happen
that $\cL_{d,d'}(E,F)$ is discrete
$($but $\not=\{0\})$.
\end{prop}
\begin{proof}
We already know that
$\cL_{d,d'}(E,F)$ is a topological group,
a vector space and that all $\|.\|_{d,d'}$-balls
around~$0$ are absolutely convex.
Hence $\cL_{d,d'}(E,F)$ is a locally convex vector
group.\\[2.5mm]
To get an example which is not a topological
vector space,
equip $\R$ with the unusual metric
$d\colon \R\times \R \to [0,\infty[$,
$d(s,t):=\frac{|s-t|}{1+|s-t|}$.
Then $\lambda \id_\R\in \cL_d(\R)$ for
each $\lambda\in \R$, and
\begin{equation}\label{fiendish}
\|\lambda \id_\R\|_{d,d}\;\geq \; 1\quad
\mbox{for all $\lambda\in \R\setminus \{0\}$.}
\end{equation}
In fact, $\|\lambda\id_\R\|_{d,d}\leq \max\{1,2|\lambda|\}<\infty$
by~(\ref{sumup}) and thus $\lambda\id_\R\in \cL_d(\R)$.
If $\lambda\not=0$, we have
\[
\lim_{t\to\infty}\frac{d(\lambda t,0)}{d(t,0)}\;=\;
\lim_{t\to\infty}
\frac{\frac{|\lambda t|}{1+|\lambda t|}}{\frac{|t|}{1+|t|}}\;=\;1
\]
and thus $\|\lambda\id_\R\|_{d,d}\geq 1$.
Hence $0$ is an isolated point in $\cL_d(\R)$
and hence $\cL_d(\R)$ is discrete
(being also a topological group).
\end{proof}
%
%
%\ma{interprnee}
\begin{rem}\label{interprnee}
In \cite[p.\,11]{Mul},
it is claimed that $\cL_{d,d'}(E,F)$
always is a Fr\'{e}chet space
(and hence a topological vector space),
contrary to Proposition~\ref{vectgp}.
If $\cL_{d,d'}(E,F)$ is not a Fr\'{e}chet space,
then the map $f''\colon U\to \cL_{d,d}(E,F)$
used in \cite[Theorem~4.7]{Mul} requires
interpretation,
as well as the use of higher
order differentiability properties
in \cite[Theorem~4.6]{Mul} and its proof.
However, a suitable interpretation
is possible (see Definition~\ref{dfnMCk}).
\end{rem}
Our simple counterexample can be generalized
further.
%
%
%\ma{bdthen}
\begin{prop}\label{bdthen}
Let $(F,d)$ be a metric Fr\'{e}chet
space with absolutely convex balls,
such that $\|\R \, x\|_d\sub \R$
is bounded for some non-zero
vector $x\in F$
$($such $x$ exists, e.g., if $F\not=\{0\}$
and $d$ has bounded image$)$.
Then $\cL_d(F)$ is not a topological
vector space.
\end{prop}
\begin{proof}
It is clear from the definition
that $\|\id_F\|_{d,d}=1$,
whence $\id_F\in \cL_d(F)$. 
We claim that $\|t\id_F\|_d\geq 1$
for all real numbers $t>0$.
If this is so, then
$\|t\id_F\|_d\not\to 0$
as $t\to 0$, whence $t\id_F\not\to 0$ in
$\cL_d(F)$. Thus scalar multiplication
$\K\times \cL_d(F)\to\cL_d(F)$
is discontinuous.
To prove the claim, let~$x$
be as described in the proposition.
Since $d$ has convex balls,
the map
$h\colon [0,\infty[\to [0,\infty[$,
$h(s):=\|sx\|_d$ is monotonically
increasing. Because~$h$ is bounded
by hypothesis and $h(1)=\|x\|_d>0$,
the limit $\lim_{s\to\infty}h(s)$ exists
and coincides with $\sigma:=\sup\, h([0,\infty[)>0$.
For each $s>0$, we have
\begin{equation}\label{stuff}
\|t\id_F\|_{d,d}\; \geq\;
\frac{\|t\id_F(sx)\|_d}{\|sx\|_d}
\; =\; \frac{h(ts)}{h(s)}\,.
\end{equation}
For $s\to\infty$, the right hand
side of~(\ref{stuff}) tends
to $\frac{\sigma}{\sigma}=1$.
Thus $\|t\id_F\|_{d,d}\geq 1$.
\end{proof}
We record a crucial property of
the metrics $d_{w,p}$
from Example~\ref{exam1}.
%
%
%\ma{crucia}
\begin{la}\label{crucia}
Assume that $w=(w_n)_{n\in \N}$
is monotonically decreasing
in the situation of Example~{\rm\ref{exam1}}.
Given non-zero vectors $x,y\in F$,
there exist minimal numbers $n,m\in \N$
such that $p_n(x)>0$
and $p_m(y)>0$, respectively. Then
%
%\ma{twovect}
\begin{equation}\label{twovect}
\sup_{t\in \K^\times} \frac{d_{w,p}(ty,0)}{d_{w,p}(tx,0)}
\;\geq \; \frac{w_m}{w_n}\,.
\end{equation}
\end{la}
\begin{proof}
Abbreviate $d:=d_{w,p}$.
For each $t\in \K^\times$,
we have $p_k(tx)=0$ for $k<n$
and thus $\|tx\|_d\leq w_n$,
entailing that
$\frac{\|ty\|_d}{\|tx\|_d}\geq \frac{\|ty\|_d}{w_n}
\geq \frac{w_m}{w_n}\frac{p_m(ty)}{1+p_m(ty)}$.
Since the right hand side tends
to $\frac{w_m}{w_n}$ as $|t|\to\infty$,
the assertion follows.
\end{proof}
The following example shows that
$\cL_d(F)$ can be quite large.
%
%
%\ma{exam3}
\begin{ex}\label{exam3}
Let
$(F_n,\|.\|_n)_{n\in \N}$
be a sequence of Banach spaces
and $w=(w_n)_{n\in \N}$
be a sequence of
real numbers $w_n>0$ such that
$\lim_{n\to\infty}\,w_n=0$.
We turn the direct product
$F:=\prod_{n\in \N}F_n$ into a metric
Fr\'{e}chet space
using the translation invariant
metric $D\colon F\times F\to [0,\infty[$,
\[
(x,y)\mto \sup_{n\in \N} \, w_n\,\frac{p_n(x_n,y_n)}{1+p_n(x_n,y_n)}
\quad
\mbox{for $x=(x_n)_{n\in \N}$ and $y=(y_n)_{n\in \N}$
in $F$}
\]
with absolutely convex balls.
Then
$\prod_{n\in \N}A_n\in \cL_d(F)$
for all $A_n\in \cL(F_n)$ such that
$\sigma :=\sup_{n\in \N}\, \|A_n\|<\infty$
holds for the operator norm,
exploiting that $\frac{\|A_n.x_n\|_n}{1+\|A_n.x_n\|_n}
\leq \frac{\|A_n\|\cdot \|x_n\|_n}{1+\|A_n\|\cdot \|x_n\|_n}
\leq \frac{\sigma \|x_n\|_n}{1+\sigma \|x_n\|_n}
\leq \max\{1,2 \sigma \} \, \frac{\|x_n\|}{1+\|x_n\|_n}$
for each $n\in \N$ and $x_n\in F_n$
(using Lemma~\ref{splobsv}).
\end{ex}
If $\cL_d(F)\not=\{0\}$, then $\cL_d(F)$
need not contain any linear
contractions except for~$0$,
as we have seen in the proof
of Proposition~\ref{vectgp}.
The next example describes
a situation where contractions
exist in abundance.
%
%
%\ma{exRN}
\begin{ex}\label{exRN}
Let $a\in \;]0,1[$ and equip
$\R^\N$ with the metric
\[
d\colon \R^\N\times\R^\N\to[0,\infty[\,, \;\;\;
d(x,y)\;:=\; \sup_{n\in \N}\, a^n\frac{|x_n-y_n|}{1+|x_n-y_n|}\,.
\]
Let $S\colon \R^\N\to \R^\N$, $x=(x_n)_{n\in \N}\mto
(0,x_1,x_2,\ldots)$ be the right shift.
Then $S$
is a linear contraction of $(\R^\N,d)$,
with $\|S\|_{d,d}=a$
(as is clear from the definition of~$d$).
If $a<\frac{1}{2}$ (and hence $\frac{a}{1-a}<1$),
let $(t_n)_{n\in \N}$ be a sequence
in~$\R$ such that $|t_n|\leq 1$ for each
$n\in \N$.
Since $\sum_{n=1}^\infty \|t_n S^n\|_d\leq
\sum_{n=1}^\infty \|S^n\|_d\leq \sum_{n=1}^\infty a^n=\frac{a}{1-a}
<\infty$,
the series $\sum_{n=1}^\infty t_n S^n$
then converges in $\cL_d(\R^\N)$,
and its limit~$A$ is a contraction
with $\|A\|_{d,d}\leq
\frac{a}{1-a}<1$.
\end{ex}
General contractions of $(\R^\N,d)$
share a property of the shift.
\begin{prop}\label{filtrprop}
Let $d$ be as in Example~{\rm\ref{exRN}}
and $F_n:=\R^{\{n,n+1,n+2,\cdots\}}\sub \R^\N$
for $n\in \N$.
If $A\in \cL_d(\R^\N)$
and $\|A\|_{d,d}<1$, then
there exists $\ell \in \N$ such that
%
%\ma{teles}
\begin{equation}\label{teles}
A.F_k\; \sub \; F_{k+\ell}\quad\mbox{for each $k\in \N$.}
\end{equation}
Moreover, {\rm(\ref{teles})}
holds for each $\ell\in \N$ such that
$\|A\|_{d,d}<a^{\ell-1}$.
\end{prop}
\begin{proof}
If the first assertion is false,
there exists $0\not=x\in F_k$
for some~$k$ such that $y:=A.x\not\in F_{k+1}$.
Let $n$ and~$m$ be as in Lemma~\ref{crucia}.
Then $n\geq k$ and $m\leq k$.
Hence $\|A\|_{d,d}\geq \frac{w_m}{w_n}\geq \frac{w_k}{w_k}=1$,
contradicting the hypothesis that $\|A\|_{d,d}<1$.
If the final assertion is false,
instead we find~$x$ with $y=A.x\not\in F_{k+\ell}$.
Then $m\leq k+\ell-1$
and we conclude as before that
$\|A\|_{d,d}\geq \frac{w_m}{w_n}\geq \frac{w_{k+\ell-1}}{w_k}=a^{\ell-1}$,
contradicting the choice of~$\ell$.
\end{proof}
Each standard metric $d_{w,p}$
(with $w$ monotonically decreasing)
goes along with a
filtration $F=F_0\supseteq F_1\supseteq F_2\supseteq \cdots$
of closed vector subspaces of~$F$,
as we shall presently see.
Each linear contraction $A\colon F\to F$
satisfies $A.F_k\sub F_{k+1}$ for
each~$k$ and hence behaves,
essentially, like the contractions
of~$\R^\N$
just discussed. More generally,
repeating the argument used to prove
Proposition~\ref{filtrprop}, we see:
%
%\ma{propshape}
\begin{prop}\label{propshape}
Let $E$ and $F$ be Fr\'{e}chet spaces.
Let $d:=d_{w,p}$
$d':=d_{v,q}$ be metrics on~$E$,
resp., $F$ of the form described
in Example~{\rm\ref{exam1}},
such that the sequences
$(w_n)_{n\in \N}$ and $(v_n)_{n\in \N}$
are monotonically decreasing.
Set $E_0:=E$ and $E_k:=\bigcap_{j=1}^k p_j^{-1}(0)$
for $k\in \N$.
Then $E=E_0\supseteq E_1\supseteq E_2\supseteq \cdots$
is a descending sequence of closed
vector subspaces of~$E$ such that $\bigcap_{k\in \N}E_k=\{0\}$.
Likewise, set $F_0:=F$ and $F_k:=\bigcap_{j=1}^k q_j^{-1}(0)$.
If $A\in \cL_{d,d'}(E,F)$
and $\|A\|_{d,d'}<1$,
then there exists $\ell\in \N$
such that
\begin{equation}\label{filtratra}
A.E_k\; \sub \; F_{k+\ell}\quad \mbox{for
each $k\in \N$.}
\end{equation}
Moreover, {\rm(\ref{filtratra})}
holds for each $\ell\in \N$ such that
$\|A\|_{d,d'} <
\inf\big\{ \frac{v_{k+\ell-1}}{w_k}\colon k\in \N\big\}$.\Punkt
\end{prop}
Let us sum up our observations
and discuss their relevance
concerning linear contractions.
We have seen that, if $0\not=A\in \cL_d(F)$,
then $tA$ need not be a contraction for
any $t\not= 0$ (no matter how small).
This naturally leads to the
question which elements
$A\in \cL_d(F)$ have the property
that $\lim_{t\to 0}tA=0$.
Since $\|tS\|_d=a$ for each $t\in \R^\times$ such that
$|t|\leq 1$
in the situation of Example~\ref{exRN}
(as a consequence of Lemma~\ref{crucia}),
we see that $\lim_{t\to 0}tA=0$
need not even hold if~$A$ is a contraction.\\[2.5mm]
Because $\cL_{d,d'}(E,F)$
is a locally convex vector group,
the following proposition provides
in particular a characterization of those $A\in
\cL_{d,d'}(E,F)$
such that $tA$ can be made arbitrarily small.
%
%
%\ma{niceelements}
\begin{prop}\label{niceelements}
Let~$E$ be a locally convex vector group
over~$\K$ and $E_0$ be its connected
component of~$0$.
The following conditions
are equivalent for $x\in E$:
\begin{itemize}
\item[\rm(a)]
$tx\to 0$ in~$E$ as $t\to 0$ in~$\K$;
\item[\rm(b)]
The map $\K\to E$, $t\mapsto tx$ is continuous
on $\K$ with the usual
topology.
\item[\rm(c)]
$x\in E_0$.
\end{itemize}
Furthermore, $E_0$ coincides
with the path component of~$0$
and $E_0$ is the largest vector subspace
of~$E$ which is a topological
vector space in the induced topology.
Also, $E_0=\bigcap_U\, \K\, U$,
where $U$ ranges through the set
of all absolutely convex $0$-neighborhoods
in~$E$.
\end{prop}
\begin{proof}
For each absolutely convex $0$-neighborhood~$U$
in~$E$, the set $\K \, U=\bigcup_{n\in \N} nU=:V$
is an open vector subspace of~$E$
and hence also closed. It follows that
$E_0\sub \bigcap_U \K\, U$.
Since~$U\cap V$ is absorbing
in~$V$ for each absolutely convex $0$-neighborhood
$U\sub E$ (by definition of~$V$)
and furthermore $U\cap V$ is absolutely convex,
we deduce that the topology induced by~$E$ on~$V$
is a vector topology. Hence~$V$ is connected
and thus $V\sub E_0$. Hence $V=E_0$.

(a)$\impl$(b): If $tx\to 0$ as $t\to 0$,
then the homomorphism of additive groups
$\K\to E$, $t\mto tx$ is continuous at~$0$
and hence continuous.

(b)$\impl$(c): If the map $f\colon \K\to E$, $f(t):=tx$
is continuous, then $f(\K)$ is path connected
and $0\in f(\K)$, whence $f(\K)\sub E_0$.
But $x=f(1)\in E_0$.

(c)$\impl$(a): If $x\in E_0$, then $\lim_{t\to 0}tx=0$ because
$E_0$ is a topological vector space,
as observed at the beginning.
\end{proof}
The identity component
$\cL_d(F)_0$ of $\cL_d(F)$ is a two-sided
ideal in $\cL_d(F)$ and hence
a (not necessarily unital) subalgebra.
It has beautiful properties.
%
%
%\ma{isCIA}
\begin{prop}\label{isCIA}
Let $(F,d)$ be a metric
Fr\'{e}chet space over~$\K$
with absolutely convex balls.
Then $\cL_d(F)_0$ is a Fr\'{e}chet
space and a so-called
continuous quasi-inverse algebra,
i.e., $\cL_d(F)_0$ is a $($not necessarily unital$)$
locally convex, associative
topological $\K$-algebra
whose group $Q(\cL_d(F)_0)$ of quasi-invertible
elements is open in~$\cL_d(F)_0$
and whose quasi-inversion
map $Q(\cL_d(F)_0)\to Q(\cL_d(F)_0)$ is continuous
$($hence $C^\infty_\K$, and even
$\K$-analytic$)$.
In some cases, $\id_F\not\in \cL_d(F)_0$.
\end{prop}
\begin{proof}
Since $\cL_d(F)_0$ is a topological vector space
and closed in the complete metric abelian group~$\cL_d(F)$,
it is a Fr\'{e}chet space.
Since~$\cL_d(F)$ is a topological ring
with bilinear multiplication, $\cL_d(F)_0$
is a topological algebra.
Let $A\in \cL_d(F)_0$ with $\|A\|_d<1$.
Then $\id_F- A\in \cL_d(F)^\times$
and $(\id_F-A)^{-1}=\sum_{n=0}^\infty A^n
=\id_F-B$ with
$B:=\sum_{n=1}^\infty A^n \in \cL_d(F)_0$
(see \cite[Theorem~4.1]{Mul}).
Here $B$ is the quasi-inverse
of~$A$ in $\cL_d(F)$
and hence also the quasi-inverse of~$A$
in $\cL_d(F)_0$, since $B\in \cL_d(F)_0$
(see \cite[Lemmas~2.3 and~2.5]{ALG}).
Thus $Q(\cL_d(F)_0)$ is a $0$-neighborhood
in $\cL_d(F)_0$ and hence open,
by \cite[Lemma~2.6]{ALG}.
The inversion map $\cL_d(F)^\times\to\cL_d(F)^\times$
is continuous by \cite[Theorem~4.1]{Mul},
whence also the quasi-inversion
$q\colon Q(\cL_d(F))\to Q(\cL_d(F))$
is continuous.
As a consequence,
the quasi-inversion map~$q_0$ of~$\cL_d(F)_0$
is continuous on some $0$-neighborhood
(because we have seen above that it coincides with~$q$
on some $0$-neighborhood).
By \cite[Lemma~2.8]{ALG},
this implies continuity of~$q_0$
on all of~$Q(\cL_d(F)_0)$.
Now $q_0$ is $C^\infty_\K$ and even
$\K$-analytic
automatically
(cf.\ Lemma~3.1, Proposition~3.2 and Proposition~3.4
in~\cite{ALG}).
Here, $\K$-analyticity is understood
as in~\cite{BaS}, or as in \cite{Mil}
and~\cite{RES}.\footnote{Since $\cL_d(F)_0$
is a Fr\'{e}chet space, real analyticity
as in~\cite{BaS} coincides with
real analyticity as in~\cite{Mil} and~\cite{RES}
(cf.\ \cite[Theorem~7.1]{BaS}).}
To complete the proof,
we recall that $t\id_\R\not\to 0$ in $\cL_d(\R)$
as $t\to 0$ for~$d$ as in the
proof of Proposition~\ref{vectgp}.
Hence $\id_\R\not\in \cL_d(\R)_0$.
\end{proof}
Unfortunately,
it frequently happens that
$\cL_d(F)_0=\{0\}$,
as the next example shows.
This can occur
even if the set of contractions
is large and $\cL_d(F)$ is
non-discrete (in which case $\cL_d(F)_0$
is not open in~$\cL_d(F)$),
for instance in the situation
of Example~\ref{exRN}.
%
%
%\ma{exm}
\begin{ex}\label{exm}
Consider a Fr\'{e}chet space~$F$,
equipped with a standard metrics~$d=d_{w,p}$,
where $w=(w_n)_{n\in \N}$
is of the form $w_n=a^n$
for some $a \in \; ]0,1[$.
Set $F_0:=F$ and $F_k:=\bigcap_{j=1}^kp_j^{-1}(0)$.
Then $\cL_d(F)_0=\{0\}$.
To see this, let $A\in \cL_d(F)_0$.
Given $\ell\in \N$, we find
$t\in \K\setminus \{0\}$
such that $\|t A\|_{d,d'}<a^{\ell-1}$.
Since $a^{\ell-1}=\inf\big\{
\frac{w_{k+\ell-1}}{w_k}\colon k\in \N\big\}$,
we deduce with Proposition~\ref{propshape}
that $A.F=tA.F=tA.F_0 \sub F_\ell$.
As~$\ell$ was arbitrary,
it follows that $A.F\sub\bigcap_{\ell\in \N}F_\ell=\{0\}$
and~so~$A=0$.
\end{ex}
The preceding example extends
to much more general situations, due
to the following proposition.
\begin{prop}\label{quasiisom}
Let $(F,d)$ be a metric Fr\'{e}chet
space with absolutely convex balls,
and such that~$d$ is bounded,
say $d(F\times F)\sub [0,M]$ with
some $M\in\;]0,\infty[$. 
Then there exists a sequence
$p=(p_n)_{n\in \N}$
of continuous seminorms $p_1\leq p_2\leq\cdots$
on~$F$ such that
that $\id_F\colon (F,d)\to (F,D)$
is a quasi-isometry for
the standard metric $D=d_{w,p}$ with
$w=(2^{-n})_{n\in \N}$.
More precisely,
\begin{equation}\label{precquasi}
\frac{1}{2}\,\|x\|_D\; \leq\;
\|x\|_d\; \leq\; \max\left\{4, 4M \right\} \|x\|_D
\quad\mbox{for all $\, x\in F$.}
\end{equation}
\end{prop}
\begin{proof}
We define $C_n:=\wb{B}^d_{2^{-n}}(0)$
for $n\in \N$ and let $p_n:=\mu_{C_n}$
be the Minkowski functional
of~$C_n$
(as in \cite[\S\,1.33]{Rud}).
Then $p_1\leq p_2\leq\cdots$
is an ascending sequence of continuous seminorms
on~$F$, with unit balls $\wb{B}^{p_n}_1(0)=C_n$.
Let $0\not= x\in F$.
We first verify the second inequality in~(\ref{precquasi}).\\[2.5mm]
If $\|x\|_d>\frac{1}{2}$, then
$x\not\in C_1$ and hence
$p_1(x)>1$, whence $\|x\|_D
\geq \frac{1}{2}\frac{p_1(x)}{1+p_1(x)}
>\frac{1}{4}\geq \frac{1}{4 M}\|x\|_d$,
as required.\\[2.5mm]
If $\|x\|_d\leq \frac{1}{2}$,
there exists a minimal $n\in \N$
such that $2^{-n}<\|x\|_d$.
Then $p_n(x)>1$ and thus
$\|x\|_D>2^{-n}\frac{p_n(x)}{1+p_n(x)}\geq 2^{-n-1}
\geq  \frac{1}{4} \|x\|_d$,
using in the last step
that $2^{-n+1}\geq \|x\|_d$
by minimality of~$n$.\\[2.5mm]
To check the first inequality,
pick $n\in \N$ minimal
such that $2^{-n}<\|x\|_D$.
Then $n\geq 2$ (since $\|x\|_D<\frac{1}{2}$ by definition of~$D$),
and $2^{1-n}\geq \|x\|_D$.
The definition of~$\|.\|_D$ as a supremum now entails
that there exists $m\in\N$
such that
$2^{-m}\frac{p_m(x)}{1+p_m(x)}>2^{-n}$.
Then $m<n$ and
$p_m(x)\geq \frac{p_m(x)}{1+p_m(x)}>\frac{1}{2^{n-m}}$.
Thus $p_m(2^{n-m}x)>1$ and hence
$2^{n-m}x\not\in \wb{B}^{p_m}_1(0)=C_m=\wb{B}^d_{2^{-m}}(0)$.
Therefore $2^{-m}<\|2^{n-m}x\|_d\leq 2^{n-m}\|x\|_d$
(using the triangle inequality) and hence
$\|x\|_d\geq 2^{-n}\geq \frac{1}{2}\|x\|_D$.
\end{proof}
Note that if $\id_F\colon (F,d)\to (F,D)$
is a quasi-isometry,
then $\cL_d(F)=\cL_D(F)$
and the identity map
$\cL_d(F)\to \cL_D(F)$
is a quasi-isometry for
the metrics on operators
determined by $\|.\|_{d,d}$
and~$\|.\|_{D,D}$.
Combination of Example~\ref{exm}
with Proposition~\ref{quasiisom}
now shows:
\begin{cor}\label{evenworse}
If $(F,d)$ is a metric Fr\'{e}chet
space with a bounded metric
and absolutely convex balls,
then $\cL_d(F)_0=\{0\}$.\,\Punkt
\end{cor}
%
%Prompted by
Varying~\cite{Mul},
we define maps
with certain metric
differentiability properties.
%
%
%\ma{dfnMCk}
\begin{defn}\label{dfnMCk}
Let $(E,d)$ and $(F,d')$
be metric Fr\'{e}chet spaces
over~$\K$, with absolutely convex balls.
Let $U\sub E$ be a locally convex subset
with dense interior
and $f\colon U\to F$ be a map.
We say that $f$ is $MC^0_\K$ if
it is continuous.
If~$f$
is~$C^1_\K$,
$f'(U)\sub \cL_{d,d'}(E,F)$
and the map $f'\colon U\to \cL_{d,d'}(E,F)$
is continuous,
then~$f$ will be called
an \emph{$MC^1_\K$-map}.\footnote{These are M\"{u}ller's
``bounded differentiable''
maps. We avoid his terminology because
of the risk of confusion with Colombeau's
venerable ``bounded differential calculus,''
and also because not boundedness is the main feature
of the approach, but Lipschitz continuity
with respect to a given choice of metrics.}
We also write $f^{(0)}:=f$ and $f^{(1)}:=f'$.
If~$f$ is $MC^1_\K$,
$x_0\in U$ and~$V\sub U$ a connected
open neighborhood of~$x_0$
(e.g., an open convex neighborhood),
then $f'(V)$ is connected and hence
contained in the connected component
$f'(x_0)+\cL_{d,d'}(E,F)_0$ of~$f'(x_0)$
in $\cL_{d,d'}(E,F)$
(cf.\ Proposition~\ref{niceelements}).
Thus
$f'|_V-f'(x_0)\colon V\to \cL_{d,d'}(E,F)_0$
is again a map between subsets of Fr\'{e}chet
spaces. This enables a recursive definition:\\[2.5mm]
If $f$ is $MC^1_\K$ and
$V$ (as before) can be chosen
for each $x_0\in U$ such that\linebreak
$f'|_V-f'(x_0)\colon V\to \cL_{d,d'}(E,F)_0$
is $MC^{k-1}_\K$,
then~$f$ is called an \emph{$MC^k_\K$-map},
and we make a piecewise
definition of~$f^{(k)}$
via $f^{(k)}|_V:=(f'|_V-f'(x_0))^{(k-1)}$
for $x_0$ and~$V$ as before.
The map $f$ is $MC^\infty_\K$ if it is
$MC^k_\K$ for each $k\in \N_0$.
\end{defn}
We mention that a suitable version
of the Chain Rule holds: Compositions
of composable $MC^k_\K$-maps
are~$MC^k_\K$ (see Lemma~\ref{invMCextended}\,(f)
in Appendix~\ref{appB}).
%
%\ma{proverpr}
\begin{rem}\label{proverpr}
In the setting
of Corollary~\ref{evenworse},
we have $\cL_d(F)_0=\{0\}$,
whence $f'$ has to be locally constant
for any map $f\colon F \supseteq U\to F$
which is $MC^1$.
Therefore, locally around a given point~$x_0$,
$f$ merely
is an affine linear map (with linear part
in $\cL_d(F)$),
since the derivative of $f-f'(x_0)$
vanishes close to~$x_0$.
%(cf.\ also \cite[Proposition~4]{BaS1}).
This observation shows that
the class of $MC^1$-maps
(used as the basis of~\cite[Theorem~4.7]{Mul})
can be quite small in some cases.
By contrast, our Theorem~B does not require continuity of
$x\mto f'(x_0)^{-1}f'(x)$ as a map into $\cL_d(F)$.
\end{rem}
%
%
%\ma{commt4}
\begin{rem}\label{commt4}
We mention that the topology on $\cL_{d,d'}(E,F)$
arising from the metric is not the only
useful one:
Occasionally, it might be convenient
to equip $\cL_{d,d'}(E,F)$
with the translation-invariant
manifold structure which makes\linebreak
$\cL_{d,d'}(E,F)_0$
an open $MC^\infty$-submanifold
(and the corresponding finer topology).\\[2.5mm]
Besides the preceding
definition of $MC^k$-maps
between metric Fr\'{e}chet
spaces,
it might be interesting
to explore the possibility
of a metric differential calculus
of~$MC^k$-maps in arbitrary
metric locally convex vector groups.\\[2.5mm]
Presumably, to obtain a meaningful
differential calculus for
mappings
between
metric locally convex vector groups
$(E,d)$ and $(F,d')$,
one should differentiate a map
$f\colon E\supseteq U \to F$
only along directions in~$E_0$;
thus $df\colon U\times E_0\to F_0$.\\[2.5mm]
Motivated by the fact that
inverses in $\cL_d(E)^\times$
close to~$\id_E$ can be expressed
in terms of the Neumann series,
it would also be natural
to consider a certain
(restrictive) class of analytic
functions between metric locally convex vector
groups, which are locally given by series of
(metrically) Lipschitz continuous, homogeneous
polynomials, with sufficiently
strong convergence.
\end{rem}
\section{{\boldmath $C^k$}-dependence
of fixed points on parameters}\label{secfp}
%\ma{secfp}
%
%
We now study the dependence
of fixed points of contractions on parameters.
In particular, we shall establish
$C^k$-dependence under
natural hypotheses.
These results form the technical backbone
of our generalizations of the inverse- and
implicit function theorems.
\begin{defn}
A mapping $f\colon X\to Y$
between metric spaces $(X,d_X)$ and $(Y,d_Y)$
is called a \emph{contraction}
if there exists $\theta\in [0,1[$
(a ``contraction constant'')
such that
$d_Y(f(x),f(y)) \leq
\theta\, d_X(x,y)$
for all $x,y\in X$.
\end{defn}
Banach's Contraction Theorem
is a paradigmatic fixed point
theorem for contractions. We recall
it as a model for the slight generalizations
which we actually need for our purposes:
%
%
%\ma{banachfix}
\begin{la}\label{banachfix}
Let $(X,d)$ be a $($non-empty$)$
complete metric
space and $f\colon X\to X$ be a contraction,
with contraction constant
$\theta\in [0,1[$.
Then $f(p)=p$ for a unique point $p\in X$.
Given any $x_0\in X$, we have $\lim_{n\to\infty}f^n(x_0)=p$.
Furthermore, the following
a priori estimate holds, for each $n\in \N_0$:
\[
d(f^n(x_0),p)\; \leq \; \frac{\theta^n}{1-\theta}\, d(f(x_0),x_0)\,.
\]
\end{la}
Unfortunately, we are not always
in the situation of this lemma. But
the simple variants compiled
in the next proposition
are flexible enough for our purposes.
%
%
%\ma{banfix2}
\begin{prop}\label{banfix2}
Let $(X,d)$ be a metric
space, $U\sub X$ be a subset and $f\colon U\to X$
be a contraction, with contraction constant~$\theta$.
Then the following holds:
\begin{itemize}
\item[\rm (a)]
$f$ has at most one fixed point.
\item[\rm (b)]
If $x_0\in U$ is a point and $n\in \N_0$
such that
$f^{n+1}(x_0)$ is defined, then
%
%\ma{iterctr}
\begin{equation}\label{iterctr}
d(f^{k+1}(x_0), f^k(x_0))\;\leq\; \theta^k\, d(f(x_0),x_0)
\end{equation}
for all $k\in \{0,\ldots, n\}$,
and $\,d(f^{n+1}(x_0), x_0)\leq \frac{1-\theta^{n+1}}{1-\theta}\,
d(f(x_0),x_0)$.
\item[\rm (c)]
If $x_0\in U$ is a point
such that $f^n(x_0)$ is defined
for all $n\in \N$, then $(f^n(x_0))_{n\in \N}$
is a Cauchy sequence in~$U$,
and
%
%\ma{makeCauex}
\begin{equation}\label{makeCauex}
d(f^{n+k}(x_0), f^n(x_0))\;\leq\; \frac{\theta^n(1-\theta^k)}{1-\theta}\,
d(f(x_0),x_0)\quad\mbox{for all $\, n,k\in \N_0$.}
\end{equation}
If $(f^n(x_0))_{n\in \N}$ converges
to some $x\in U$, then $x$ is a fixed point
of~$f$, and
%
%\ma{aprio2}
\begin{equation}\label{aprio2}
d(x, f^n(x_0))\;\leq\; \frac{\theta^n}{1-\theta}\,
d(f(x_0),x_0)\quad\mbox{for all $\,n \in \N_0$.}
\end{equation}
If $f^n(x_0)$ is defined for all $n\in \N$
and $f$ has a fixed point~$x$, then
$f^n(x_0)\to x$ as $n\to\infty$.
\item[\rm (d)]
Assume that $U=\wb{B}_r(x_0)$
is a closed ball of radius~$r$
around a point $x_0\in X$,
and
$d(f(x_0),x_0)\leq (1-\theta)r$.
Then $f^n(x_0)$ is defined for all
$n\in \N_0$.
Hence $f$ has a fixed point
inside~$\wb{B}_r(x_0)$, provided $X$
is complete.
Likewise, $f$ has a fixed point in the open
ball $B_r(x_0)$
if $X$ is complete,
$U=B_r(x_0)$,
and $d(f(x_0),x_0)<(1-\theta)r$.
\end{itemize}
\end{prop}
\begin{proof}
(a) If $x, y\in U$ are fixed points of~$f$,
then $d(x,y)=d(f(x),f(y))\leq \theta d(x,y)$,
entailing that $d(x,y)=0$ and thus $x=y$.

(b) For $k=0$, the formula (\ref{iterctr})
is trivial.
If $k<n$ and
$d(f^{k+1}(x_0),f^k(x_0))\leq \theta^k\, d(f(x_0),x_0)$,
then
$d(f^{k+2}(x_0),f^{k+1}(x_0))
=d(f(f^{k+1}(x_0)), f(f^k(x_0)))\leq $\linebreak
$\theta \, d(f^{k+1}(x_0), f^k(x_0))
\leq \theta^{k+1}\, d(f(x_0),x_0)$.
Thus (\ref{iterctr}) holds in general.

Using the triangle inequality and
the summation formula for the geometric series,
we obtain the estimates $d(f^{n+1}(x_0),x_0)
\leq\sum_{k=0}^n d(f^{k+1}(x_0),f^k(x_0))
\leq \sum_{k=0}^n \theta^k\, d(f(x_0),x_0)
=\frac{1-\theta^{n+1}}{1-\theta}\, d(f(x_0),x_0)$,
as asserted.

(c) Using both of the estimates from (b), we obtain
\[
d(f^{n+k}(x_0),f^n(x_0)) \;\leq\;
\frac{1-\theta^k}{1-\theta}\, d(f^{n+1}(x_0),f^n(x_0))
\;\leq\;
\frac{1-\theta^k}{1-\theta}\, \theta^n d(f(x_0),x_0)\,.
\]
Thus (\ref{makeCauex}) holds,
and thus $(f^n(x_0))_{n\in \N}$ is a Cauchy sequence.
If $f^n(x_0)\to x$ for some $x\in U$,
then $x=\lim_{n\to\infty}f^{n+1}(x_0)=f(\lim_{n\to\infty}f^n(x_0))=
f(x)$ by continuity of~$f$, whence indeed~$x$
is a fixed point of~$f$.
Letting now $k\to\infty$ in (\ref{makeCauex}),
we obtain (\ref{aprio2}).

To prove the final assertion, assume that $f$ has a fixed
point $x$ and that $f^n(x_0)$
is defined for all~$n$.
We choose a completion $\wb{X}$
of~$X$ (with $X\sub \wb{X}$)
and let $\wb{U}$ be the closure of~$U$ in~$\wb{X}$.
Then $f$ extends to a contraction
$\wb{U}\to \wb{X}$, which we also denote
by~$f$. Since $(f^n(x_0))_{n\in \N}$
is a Cauchy sequence in $\wb{U}$ and $\wb{U}$
is complete, we deduce that $f^n(x_0)\to y$
for some $y\in \wb{U}$.
Then both $y$ and $x$ are fixed points of~$f$
and hence $x=y$, by~(a).

(d) We show by induction that $f^n(x_0)$ is defined
for all $n\in \N$.
For $n=1$, this is trivial.
If $f^n(x_0)$ is defined,
then
\[
d(f^n(x_0),x_0)\, \leq\,
\frac{1-\theta^n}{1-\theta}\,d(f(x_0),x_0)
\,\leq\,
\frac{1-\theta^n}{1-\theta}\,(1-\theta)r\, \leq \, r
\]
by~(b) and thus $f^n(x_0)\in \wb{B}_r(x_0)$,
whence also $f^{n+1}(x_0)=f(f^n(x_0))$ is defined.
Then $f^n(x_0)\in \wb{B}_r(x_0)$ for each $n\in \N$.
By (c), $(f^n(x_0))_{n\in \N}$ is a Cauchy
sequence. If $X$ is complete, then
so is $\wb{B}_r(x_0)$ and
thus $(f^n(x_0))_{n\in \N}$
converges to some point $x\in \wb{B}_r(x_0)$,
which is a fixed point of $f$ by~(c).
Finally, if $U=B_r(x_0)$ and $d(f(x_0),x_0)<(1-\theta)r$,
there exists $s\in \;]0,r[$ such that
$d(f(x_0),x_0)\leq (1-\theta)s$.
By the preceding, $f^n(x_0)\in \wb{B}_s(x_0)\sub B_r(x_0)$
for all $n\in \N$,
and $f$ has a fixed
point in $\wb{B}_s(x_0)\sub B_r(x_0)$.
\end{proof}
A certain class of contractions
of Fr\'{e}chet spaces is of utmost importance
for our purposes.
\begin{defn}
Let $(F,d)$ be a metric Fr\'{e}chet space
over~$\K$
and $U\sub F$.
A map $f\colon U\to F$ is called
a \emph{special contraction}
if there exists $\theta\in [0,1[$
such that
\[
(\forall s\in \K)\,(\forall x,y\in U)\quad
d(sf(x),sf(y))\;\leq\; \theta\, d(sx,sy)\,.
\]
We then call $\theta$ a \emph{special
contraction constant} for~$f$.
\end{defn}
%
%
%\ma{charactspec}
\begin{la}\label{charactspec}
Let $(F,d)$ be a metric Fr\'{e}chet space
with absolutely convex balls,
$U\sub F$ be a locally convex subset
with dense interior and $f\colon U\to F$
be a $C^1$-map.
Consider the following conditions:
\begin{itemize}
\item[\rm(a)]
$f$ is a special contraction;
\item[\rm(b)]
$f'(U)\sub \cL_d(F)$
and
$\sup_{x\in U}\|f'(x)\|_{d,d}<1$.
\end{itemize}
Then {\rm(a)} implies {\rm(b)}.
If $U$ is convex, then {\rm(a)} and {\rm(b)}
are equivalent.
\end{la}
\begin{proof}
Assume that $f$ is a special contraction with special
contraction constant $\theta\in [0,1[$.
Let $x\in U^0$, $y\in F$
and $t\in \K^\times$ with $x+ty\in U$.
Since $f$ is a special contraction,
we have $\|\frac{1}{t}(f(x+ty)-f(x))\|_d
\leq \theta\|\frac{1}{t}((x+ty)-x)\|_d=\theta \|y\|_d$.
Letting $t\to 0$, we deduce
that $\|f'(x).y\|_d\leq \theta \|y\|_d$.
Since $U^0$ is dense in~$U$, it follows
by continuity that
$\|f'(x).y\|_d\leq\theta\|y\|_d$
for all $x\in U$ and $y\in F$.
For each $x\in U$, this gives
$\|f'(x)\|_{d,d}\leq\theta$,
since $y$ was arbitrary.
Thus
\,$\sup\, \|f'(U)\|_{d,d}\leq\theta$.\\[2.5mm]
Conversely, suppose that
$U$ is convex and $\theta:=\sup\,\|f'(U)\|_{d,d}<1$.
Given $x,y\in U$ and $s\in \K$,
the map $\gamma\colon [0,1]\to F$,
$\gamma(t):=sf(x+t(y-x))$ is $C^1_\R$
and
\[
sf(y)-sf(x)\;=\; \gamma(1)-\gamma(0)\;=\;
\int_0^1 \gamma'(t)\, dt\,,
\]
where $\gamma'(t)=sf'(x+t(y-x)).(y-x)=
f'(x+t(y-x)).(sy-sx)$
with $\|\gamma'(t)\|_d\leq \|f'(x+t(y-x))\|_{d,d}.\|sy-sx\|_d
\leq\theta \|sy-sx\|_d$.
Hence $d(sf(y),sf(x))\leq \theta d(sy,sx)$,
by Lemma~\ref{estinteg},
and so~$f$ is a special contraction.
\end{proof}
We are interested in uniform
families of contractions.
%
%
%\ma{defunicon}
\begin{defn}\label{defunicon}
Let $(F,d)$ be a metric Fr\'{e}chet space over~$\K$,
$U\sub F$ and~$P$ be a set.
A family $(f_p)_{p\in P}$ of mappings
$f_p\colon U\to F$
is called a
\emph{uniform family of contractions}
if there exists $\theta\in [0,1[$
(a ``uniform contraction constant'')
such that
\[
\|f_p(x)-f_p(y)\|_d\,\leq\, \theta\|x-y\|_d\quad
\mbox{for all $x,y\in U$ and $p\in P$.}
\]
If
$\|s(f_p(x)-f_p(y))\|_d\leq
\theta\|s(x-y)\|_d$
for all $s\in \K$, $x,y\in U$ and $p\in P$,
we call $(f_p)_{p\in P}$
a \emph{uniform family of special contractions}
and $\theta$ a \emph{uniform special contraction constant}.
\end{defn}
If~$U$ is closed
and each $f_p$ is a self-map
of~$U$ here, then
Banach's Contraction Theorem
ensures
that, for each $p\in P$, the map
$f_p$ has a unique
fixed point $x_p$.
Our goal is to understand
the dependence of $x_p$ on the parameter~$p$.
In particular, for $P$ a subset
of a topological $\K$-vector space,
we want to find conditions
ensuring that the map $P\to F$,
$p\mto x_p$ is continuously differentiable.
Let us discuss continuous dependence of fixed points
on parameters first.
%
%
%\ma{preparadep}
\begin{la}\label{preparadep}
Let $P$ be a topological space and
$(F,d)$ be a metric Fr\'{e}chet space.
Let $U\sub F$ and
$f\colon P\times U\to F$ be a continuous
map such that $(f_p)_{p\in P}$
is a uniform family of contractions,
where $f_p:=f(p,\sbull)\colon U\to F$.
We assume that $f_p$ has a fixed point
$x_p$, for each $p\in P$.
Furthermore,
we assume that $U$ is open or
$f(P\times U)\sub U$
$($whence every $f_p$ is a self-map
of $U)$.
Then the map
$\phi\colon P\to F$, $\phi(p):=x_p$
is continuous.
\end{la}
\begin{proof}
Let $\theta\in [0,1[$ be a uniform contraction constant
for $(f_p)_{p\in P}$.
Given $p\in P$ and $\ve>0$,
we find a neighborhood $Q\sub P$ of~$p$
such that $\|x_p-f_q(x_p)\|_d=\|f(p,x_p)-f(q,x_p)\|_d
\leq (1-\theta)\ve$
for all $q\in Q$.
If $f(P\times U)\sub U$,
then $\|x_p-x_q\|_d\leq \frac{1}{1-\theta}\|x_p-f_q(x_p)\|_d
\leq \ve$, by (\ref{aprio2})
in Proposition~\ref{banfix2}\,(c).
If $U$ is open, we may
assume that $\wb{B}_\ve^d(x_p)\sub U$
after shrinking~$\ve$ and~$Q$.
Then Proposition~\ref{banfix2}\,(d)
applies to $f_q$
as a map $\wb{B}_\ve^d(x_p)\to F$
for each $q\in Q$,
showing that $f_q^n(x_p)$ is defined
for each $n\in \N$ and
$x_q=\lim_{n\to\infty}f_q^n(x_p)\in \wb{B}_\ve^d(x_p)$,
that is, $\|x_p-x_q\|_d\leq\ve$.
\end{proof}
%
%
%\ma{prepdep2}
\begin{la}\label{prepdep2}
Let~$E$ be a topological
$\K$-vector space
and $P\sub E$ be a subset with dense
interior.
Let $(F,d)$ be a metric Fr\'{e}chet space over~$\K$
and
$U\sub F$ be a subset with dense interior.
Furthermore,
let $f\colon P\times U\to F$ be a
$C^1_\K$-map such that $(f_p)_{p\in P}$
is a uniform family of special contractions,
where $f_p:=f(p,\sbull)\colon U\to F$.
We assume that $f_p$ has a fixed point
$x_p$, for each $p\in P$.
Finally,
we assume that $U$ is open or
$f(P\times U)\sub U$
$($whence every $f_p$ is a self-map
of $U)$.
Then the map
$\phi\colon P\to F$, $\phi(p):=x_p$
is~$C^1_\K$.
\end{la}
\begin{proof}
Since~$f$ is continuous,
Lemma~\ref{preparadep}
shows that~$\phi$ is continuous.
Thus $\phi^{]1[}$ (as in (\ref{defopf})
in Definition~\ref{C12}) is continuous.
To see that $\phi$ is~$C^1$,
it only remains to show that,
for all $p_0\in P$ and $q_0\in E$,
there exists an open neighborhood
$W\sub P^{[1]}$ of $(p_0,q_0,0)$
and a continuous map $g\colon W\to F$
which extends the difference quotient map
$\phi^{]1[}|_{W\cap P^{]1[}}\colon W\cap P^{]1[}\to F$.
Then $\phi^{]1[}$ has a continuous
extension $\phi^{[1]}$ to all
of $P^{[1]}$
(cf.\ \cite[Exercise~3.2\,A\,(b)]{Eng})
and thus~$\phi$ will be~$C^1$.
Our strategy is the following: We write
\begin{equation}
(f^{n+1}_{p+tq}(x_p)-x_p)/t\;=\;
\sum_{k=0}^n(f_{p+tq}^{k+1}(x_p)-f^k_{p+tq}(x_p))/t
\end{equation}
for
$(p,q,t)$ in a suitable neighborhood
$W$ of $(p_0,q_0,0)$, with $t\not=0$.
For $W$ sufficiently small,
the left hand side converges to $\frac{x_{p+tq}-x_p}{t}
=\frac{\phi(p+tq)-\phi(p)}{t}$
as $n\to\infty$.
Furthermore, we can achieve that
each term on the right hand side extends continuously
to all of $W$, and
that the series
converges uniformly to a continuous
function on~$W$.
This will be our desired continuous extension~$g$.\\[1.4mm]
Let us carry this out in detail.\\[1.4mm]
Case~1. If $f(P\times U)\sub U$,
we set $W_0:=P^{[1]}$.\\[1.4mm]
Case~2. Otherwise, $U$ is open,
whence there exists $\ve>0$ such that
$\wb{B}_{2\ve}^d(x_{p_0})\sub U$.
Since $f$ and $\phi$ are continuous
and $f_{p_0}(x_{p_0})=x_{p_0}$,
we find an open neighborhood $Q\sub P$ of $p_0$
such that $\|x_p-x_{p_0}\|_d\leq\ve$
and $\|x_p-f_q(x_p)\|_d\leq (1-\theta)\ve$
for all $p,q\in Q$.
Then $f_q^k(x_p)$ is defined for all $k\in \N_0$,
$f_q^k(x_p)\in B_\ve^d(x_p)$,
and $x_q=\lim_{k\to\infty}f^k_q(x_p) \in \wb{B}_\ve^d(x_p)$,
by Proposition~\ref{banfix2}\,(d).
We now set $W_0:=Q^{[1]}$
and note that, if $(p,q,t)\in Q^{[1]}$,
then $p,p+tq\in Q$,
whence $f_{p+tq}^k(x_p)$
is defined
for all $k\in \N_0$
and $\lim_{k\to\infty}f_{p+tq}^k(x_p)=x_{p+tq}$
(by the preceding considerations).

In either case, we define
\[
h_0\colon W_0 \to F\,,\qquad
h_0(p,q,t)\,=\, f^{[1]}(p,x_p,q,0,t)
\]
and note that $h_0$ is a continuous map such that
%
%\ma{dexter}
%
\begin{equation}\label{dexter}
h_0(p,q,t)=(f_{p+tq}(x_p)-f_p(x_p))/t
=(f_{p+tq}(x_p)-x_p)/t\quad\mbox{if $t\not=0$.}
\end{equation}
For all $k\in \N$ and $(p,q,t)\in W_0$ with $t\not=0$,
we have
%
%\ma{dext2}
\begin{eqnarray}
\lefteqn{\frac{f^{k+1}_{p+tq}(x_p)-f^k_{p+tq}(x_p)}{t}}\qquad\notag \\
&=&
\frac{f \bigl( p+tq, f^{k-1}_{p+tq} (x_p)+t
\, \frac{f^k_{p+tq}(x_p)-f^{k-1}_{p+tq}(x_p)}{t} \bigr)
- f(p+tq, f^{k-1}_{p+tq} (x_p))}{t} \notag \\
&=& f^{[1]} \Big(
p+tq, f^{k-1}_{p+tq} (x_p), 0, \frac{f^k_{p+tq}(x_p)-f^{k-1}_{p+tq} (x_p)}{t},
t \Big)\,. \label{dext2}
\end{eqnarray}
Recursively, we define
\[
h_k \colon W_0 \to F\,,\quad
h_k (p,q,t)\, :=\,
f^{[1]}\bigl(p+tq,f^{k-1}_{p+tq}(x_p),0, h_{k-1}(p,q,t),t\bigr)
\]
for $k\in \N$.
A simple induction based on (\ref{dexter})
and (\ref{dext2}) shows that the definition of
$h_k$ makes sense for each $k\in \N_0$,
and that
%
%\ma{intph}
\begin{equation}\label{intph}
h_k(p,q,t)\,=\, \frac{f^{k+1}_{p+tq}(x_p)-f^k_{p+tq}(x_p)}{t}
\quad \mbox{for all $(p,q,t)\in W_0$ with $t\not=0$.}
\end{equation}
The function $h_0\colon W_0 \to F$,
$(p,q,t)\mto f^{[1]}(p,x_p,q,0,t)$
being continuous,
we find an open neighborhood
$W \sub W_0$ of $(p_0,q_0,0)$
and $C\in [0,\infty[$
such that
\[
\|f^{[1]}(p,x_p,q,0,t)\|_d\;\leq\; C\qquad \mbox{for all $(p,q,t)\in W$.}
\]
For all $(p,q,t)\in W$ such that $t\not=0$,
we have
\begin{eqnarray*}
\|t^{-1}f_{p+tq}^{k+1}(x_p)-t^{-1}f_{p+tq}^k(x_p)\|_d
&\leq & \theta^k\|t^{-1}f_{p+tq}(x_p)-t^{-1}x_p\|_d\\
&= & \theta^k \|f^{[1]}(p,x_p,q,0,t)\|_d
\; \leq \; \theta^k C\,;
\end{eqnarray*}
to obtain the inequality, we used repeatedly
that $f_{p+tq}$
is a special contraction with special contraction constant~$\theta$.
Combining the preceding estimates
with (\ref{intph}), for each $k\in \N_0$
we see that
$\|h_k(p,q,t)\|_d\leq \theta^k C$
for all $(p,q,t)\in W$ such that $t\not=0$,
and thus
%
%\ma{estih}
\begin{equation}\label{estih}
\|h_k(p,q,t)\|_d\;\leq \; \theta^k C
\quad\mbox{for all $(p,q,t)\in W$,}
\end{equation}
because $h_k$ is continuous and $W\cap P^{]1[}$
is dense in~$W$.
As a consequence, $\sum_{k=0}^\infty \sup\|h_k(W)\|_d \leq
\sum_{k=0}^\infty\theta^k C=\frac{1}{1-\theta}C<\infty$,
whence the series
$\sum_{k=0}^\infty h_k|_W$ of continuous
functions into $(F,d)$ converges
uniformly. Thus
\[
g(p,q,t)\; :=\; \sum_{k=0}^\infty \, h_k (p,q,t)
\]
exists for all $(p,q,t)\in W$,
and $g\colon W\to F$ is continuous.
It only remains to observe that
\[
\frac{f_{p+tq}^{n+1}(x_p)-x_p}{t}
\; =\;  \sum_{k=0}^n \, \frac{f_{p+tq}^{k+1}(x_p)
-f_{p+tq}^k(x_p)}{t}
\; =\;
 \sum_{k=0}^n \, h_k(p,q,t)
\]
for all $(p,q,t)\in W$ such that $t\not=0$.
Since the left hand side converges to $\frac{x_{p+tq}-x_p}{t}$
as $n\to\infty$
and the right hand side converges to $g(p,q,t)$,
we obtain
\[
\frac{\phi(p+tq)-\phi(p)}{t}\;=\;\sum_{k=0}^\infty \, h_k(p,q,t)
\;=\; g(p,q,t)\,.
\]
Thus $g\colon W\to F$ is a continuous map
which extends $\phi^{]1[}|_{W\cap P^{]1[}}$,
as desired. This completes the proof.
\end{proof}
We are now in the position to prove
Theorem~\hspace*{-.2mm}D \hspace*{.4mm}from the introduction
(the main result of
this section).\\[3mm]
{\bf Proof of Theorem~D.}
(a) For the proof of~(a),
we assume only that~$P$ is a
topological space
(the hypotheses that $P\sub E$ for a topological
vector space~$E$ is not required).
Let $\theta\in [0,1[$ be a uniform contraction constant
for $(f_p)_{p\in P}$.
If $p\in Q$, there is $r>0$
such that $\wb{B}_r^d(x_p)\sub U$.
There is a neighborhood
$S\sub Q$ of~$p$ such that
$\|f_q(x_p)-x_p\|_d=\|f_q(x_p)-f_p(x_p)\|_d\leq (1-\theta)r$
for all $q\in S$.
Now Proposition~\ref{banfix2}\,(d)
shows that $f_q$ has a fixed point $x_q$ in
$\wb{B}_r^d(x_p)$,
for each $q\in S$. Thus $S\sub Q$ and we deduce that
$Q$ is open. By Lemma~\ref{preparadep}, the map
$\phi\colon Q\to U$, $\phi(p):=x_p$ is continuous.

(b) Let $\theta\in [0,1[$ be a uniform special
contraction constant
for $(f_p)_{p\in P}$.
We may assume that $k<\infty$.
The proof is by induction on $k\in \N$.
Since~$\phi$ is continuous by~(a),
and since we only need
to check that~$\phi$ is~$C^k$ on an open neighborhood
of a given point $p_0\in Q$, after replacing~$U$ by some
open ball around~$\phi(p_0)$ and replacing~$Q$
by a smaller open neighborhood of~$p_0$ in~$Q$
we may assume henceforth that~$U$ is convex.\\[2.5mm]
The case~$k=1$ is covered by
Lemma~\ref{prepdep2}.
Since $\phi(p)=f(p,\phi(p))$ for all
$p\in Q$, the Chain Rule shows that
\begin{equation}\label{fixford}
\phi'(p).q\;=\; f'(p,\phi(p)).(q,\phi'(p).q)\quad
\mbox{for all $p\in Q$ and $q\in E$.}
\end{equation}
If $k\geq 2$ and $\phi$ is $C^{k-1}$ by induction,
then the map
\[
g\colon (Q\times E)\times F\to F\,,\quad
g(p,q,y)\; :=\; f'(p,\phi(p)).(q,y)
\]
is $C^{k-1}$. By linearity in $(q,y)$,
for the partial differential
with respect to~$y$ we obtain
$g_{(p,q)}'(y).z= f'(p,\phi(p)).(0,z)
=f_p'(\phi(p)).z$,
whence $g_{(p,q)}'(y)=f_p'(\phi(p))$
and thus $\|g_{(p,q)}'(y)\|_{d,d}=\|f_p'(\phi(p))\|_{d,d}\leq\theta$,
using Lemma~\ref{charactspec}.
Hence, by Lemma~\ref{charactspec},
$(g_{(p,q)})_{(p,q)\in Q\times E}$
is a uniform family of special
contractions.
Since $d\phi(p,q)=\phi'(p).q$ is the fixed point
of $g_{(p,q)}$ by~(\ref{fixford}),
applying the inductive hypothesis
we see that $d\phi\colon Q\times E\to F$
is $C^{k-1}$. Hence $\phi$ is $C^k$,
which completes the inductive proof.\,\Punkt
%
%
%
%
%
%
%
%
%
%
%
%
%
%
%
%
%
%
%
%
%
%
%
%
%\ma{snultra}
\section{Preparatory results concerning local inverses}\label{snultra}
We now prove an Inverse
Function Theorem for self-maps
of a Fr\'{e}chet space,
which provides
local inverses
which are Lipschitz continuous
with respect to a given metric.
We also provide a variant
dealing with families of
local inverses,
and use it to construct
continuous implicit functions.\\[3mm]
Our first theorem, and its proof,
is the direct analogue of \cite[Theorem~5.3]{IM2}
(dealing with Banach spaces)
for Fr\'{e}chet spaces.
It is a variant of \cite[Theorem~4.5]{Mul}.
Parts of the proof will be re-used later.
%
%
%\ma{newton}
\begin{thm}[Lipschitz Inverse Function Theorem]\label{newton}
Let $(E,d)$ be a metric Fr\'{e}chet space over~$\K$,
with absolutely convex balls.
Let $r>0$, $x\in E$,
and $f\colon  B_r^d(x)\to E$ be a map.
We suppose that there exists $A\in \cL_d(E)^\times$
such that
%
%\ma{shrink}
\begin{equation}\label{shrink}
\sigma:=\sup\left\{
\frac{\|f(z)-f(y)-A.(z-y)\|_d}{\|z-y\|_d}\colon 
\;\mbox{$y,z\in B_r(x)$, $y\not=z$}\right\}<\frac{1}{\|A^{-1}\|_d}\,.
\end{equation}
Then the following holds:
\begin{itemize}
\item[\rm (a)]
$f$ has open image and is a homeomorphism
onto its image.
\item[\rm (b)]
The inverse map $f^{-1}\colon f(B_r(x))\to B_r(x)$
is Lipschitz continuous with\linebreak
respect to the metric~$d$, with
%
%\ma{lipofinv}
\begin{equation}\label{lipofinv}
\Lip(f^{-1})\; \leq\; \frac{1}{\|A^{-1}\|^{-1}_d-\sigma}\,.
\end{equation}
\item[\rm (c)]
Abbreviating $a :=\|A^{-1}\|^{-1}_d-\sigma>0 $ and
$b :=\|A\|_d+\sigma$, we have
%
%\ma{quasiisom5}
\begin{equation}\label{quasiisom5}
a \|z-y\|_d\;\leq\; \|f(z)-f(y)\|_d\; \leq\;
b \|z-y\|_d\quad
\mbox{for all $y,z\in B_r(x)$.}
\end{equation}
\item[\rm (d)]
The following estimates for the size of images
of balls are available:
For every $y\in B_r(x)$ and $s\in \;]0,r-\|y-x\|_d]$,
%
%\ma{schachtel5}
\begin{equation}\label{schachtel5}
B_{a s}(f(y))\; \sub \; f(B_s(y))\; \sub \; B_{b s}(f(y))
\end{equation}
holds.
In particular,
$B_{a r}(f(x)) \sub f(B_r(x)) \sub B_{b r}(f(x))$.
\end{itemize}
\end{thm}
%
%
%\ma{againlp}
\begin{rem}\label{againlp}
Note that the condition (\ref{shrink})
means that the
remainder term
\[
\tilde{f} \colon B_r(x)\to E,\quad \tilde{f}(y):=
f(y)-f(x)-A.(y-x)
\]
in the affine-linear approximation $f(y)=f(x)+A.(y-x)+\tilde{f}(y)$
is Lipschitz continuous with respect to the metric~$d$,
with $\Lip(\tilde{f})=\sigma<\|A^{-1}\|^{-1}_d$.
\end{rem}
\begin{rem}
To understand the constants in Theorem~\ref{newton} better,
we recall that $\|A^{-1}\|^{-1}_d$
can be interpreted as a minimal
distortion factor, in the following sense:
For each $u\in E$, we have
$\|u\|_d=\|A^{-1}.(A.u)\|_d\leq \|A^{-1}\|_d\cdot
\|A.u\|_d$ and thus
%
%\ma{leastdisto}
\begin{equation}\label{leastdisto}
\|A.u\|_d\;\geq\; \|A^{-1}\|^{-1}_d\|u\|_d\quad
\mbox{for all $u\in E$.}
\end{equation}
Thus $A$ increases the distance of each given vector
from~$0$
by a factor of at least $\|A^{-1}\|^{-1}_d$.
Furthermore, $\|A^{-1}\|^{-1}_d$
is maximal among such factors,
as one verifies by going backwards
through the preceding lines.
Similarly, since $A^{-1}B_s(0)\sub B_{\|A^{-1}\|_ds}(0)$
and thus $B_s(0)\sub A.B_{\|A^{-1}\|s}(0)$
for each $s>0$, we find that
%
%\ma{leastinfla}
\begin{equation}\label{leastinfla}
A.B_s(0)\;\supseteq \; B_{\|A^{-1}\|^{-1}s}(0)\quad
\mbox{for each $s>0$.}
\end{equation}
\end{rem}
{\bf Proof of Theorem~\ref{newton}.}
(c) Given $y,z\in B_r(x)$, we have
\begin{eqnarray*}
\|f(z)- f(y)\|_d &= & \|f(z)-f(y)-A.(z-y)+A.(z-y)\|_d\\
&\leq & \|f(z)-f(y)-A.(z-y)\|_d+\|A.(z-y)\|_d\\
&\leq & (\sigma+\|A\|_d)\|z-y\|_d\; =\; b\|z-y\|_d
\end{eqnarray*}
and
\begin{eqnarray*}
\|z-y\|_d&= &\|A^{-1}.(f(z)-f(y)-A.(z-y))-(A^{-1}.f(z)-A^{-1}.f(y))\|_d\\
&\leq&
\|A^{-1}\|_d\cdot \|f(z)-f(y)-A.(z-y)\|_d+\|A^{-1}.f(z)-A^{-1}.f(y)\|_d\\
&\leq & \sigma\|A^{-1}\|_d\cdot \|z-y\|_d+\|A^{-1}\|_d\cdot\|f(z)-f(y)\|_d\,,
\end{eqnarray*}
whence $\|f(z)-f(y)\|_d\geq (\|A^{-1}\|^{-1}_d-\sigma)\|z-y\|_d
=a\|z-y\|_d$. Thus (\ref{quasiisom5}) holds.

(b) As a consequence of (\ref{quasiisom5}),
$f$ is injective,
a homeomorphism onto its image,
and $\Lip(f^{-1})\leq a^{-1}=(\|A^{-1}\|_d^{-1}-\sigma)^{-1}$.

(d) Suppose that $y\in B_r(x)$
and $s\in \;]0,r-\|y-x\|_d]$.
By (\ref{quasiisom5}),
we have $f(B_s(y))\sub B_{bs}(f(y))$,
proving the second half of (\ref{schachtel5}).
We now show that
%
%\ma{sfhf}
\begin{equation}\label{sfhf}
f(y)+A.B_{\alpha s}(0)\; \sub \; f(B_s(y))\,,
\end{equation}
where $\alpha:=a\|A^{-1}\|_d=1-\sigma\|A^{-1}\|_d$.
Then also the first half of (\ref{schachtel5})
will hold, as
$A.B_{\alpha s}(0)\supseteq B_{\|A^{-1}\|^{-1}_d\alpha s}(0)
=B_{as}(0)$
by (\ref{leastinfla}).
To prove (\ref{sfhf}),
let $c \in f(y)+A.B_{\alpha s}(0)$.
There exists $t\in \;]0,1[$ such that $c \in f(y)+A.\wb{B}_{t\alpha s}(0)$.
For $v\in \wb{B}_{st}(y)$,
we define
\[
g(v)\, :=\, v-A^{-1}.(f(v)-c)\,.
\]
Then $g(v)\in \wb{B}_{st}(y)$, because
\begin{eqnarray*}
\|g(v)-y\|_d &\leq &
\underbrace{\|v-y-A^{-1}.f(v)+A^{-1}.f(y)\|_d}_{\leq \|A^{-1}\|_d\,\sigma
\|v-y\|_d
\leq \|A^{-1}\|_d\, \sigma st}
+ \underbrace{\|A^{-1}.c -A^{-1}.f(y)\|_d}_{\leq t\alpha s}\\
&\leq & (\|A^{-1}\|_d\, \sigma + \alpha )st=st\,.
\end{eqnarray*}
Thus $g(\wb{B}_{st}(y))\sub \wb{B}_{st}(y)$.
The map $g\colon  \wb{B}_{st}(y)\to \wb{B}_{st}(y)$
is a contraction, since
\begin{eqnarray}
\|g(v)-g(w)\|_d & = & \|v-w-A^{-1}.(f(v)-f(w))\|_d\nonumber\\
&\leq & \|A^{-1}\|_d \cdot \|f(v)-f(w)-A.(v-w)\|_d \nonumber\\
& \leq & \sigma \cdot \|A^{-1}\|_d\cdot \|v-w\|_d \label{hencecontr}
\end{eqnarray}
for all $v,w\in \wb{B}_{st}(y)$,
where $\sigma\|A^{-1}\|_d<1$.
By Banach's Contraction Theorem (Lemma~\ref{banachfix}),
there exists a unique element $v_0\in
\wb{B}_{st}(y)$ such that $g(v_0)=v_0$
and hence $f(v_0)= c$.

(a) We have already seen that $f$ is a homeomorphism onto
its image. As a consequence of (d), the image of $f$
is open.\vspace{3mm}\Punkt

\noindent
We are now in the position to formulate
the first version of an inverse function
theorem with parameters.
The result, and its proof,
can be re-used later to prove
the corresponding results for $C^k$-maps.
%
%
%
%
%\ma{newtpara}
\begin{thm}[Continuous families of local inverses]\label{newtpara}
Let $(F,d)$ be a metric Fr\'{e}chet space
with absolutely convex balls,
and~$P$ be a topological space.
Let $r>0$, $x\in F$,
and $f\colon  P\times B \to F$ be a continuous
mapping, where $B:=B_r^d(x)$.
Given $p\in P$, we abbreviate
$f_p:=f(p,\sbull)\colon  B\to F$.
We suppose that there exists $A\in \cL_d(F)^\times$
such that
%
%\ma{shrinkpara}
\begin{equation}\label{shrinkpara}
\sigma:=\sup\left\{
\frac{\|f_p(z)-f_p(y)-A.(z-y)\|_d}{\|z-y\|_d}\colon 
\;\mbox{$p\in P$, $y,z\in B$, $y\not=z$}\right\}<\frac{1}{\|A^{-1}\|_d}\,.
\end{equation}
Then the following holds:
\begin{itemize}
\item[\rm (a)]
$f_p(B)$ is open in~$F$ and $f_p|_B$ is a homeomorphism
onto its image, for each $p\in P$.
\item[\rm (b)]
The set
$W:=\bigcup_{p\in P}\, \{p\}\times f_p(B)$
is open in $P\times F$, and the map
$\psi\colon  W\to F$,
$\psi(p,z):=(f_p|_B^{f_p(B)})^{-1}(z)$ is continuous.
\item[\rm (c)]
The map
$\xi \colon  P\times B\to W$,
$\xi (p,y):=(p,f(p,y))$
is a homeomorphism, with inverse given by
$\xi^{-1}(p,z)=(p,\psi(p,z))$.
\end{itemize}
\end{thm}
\begin{proof}
By Theorem~\ref{newton}, applied to~$f_p$,
the set $f_p(B)$ is open in~$F$ and $f_p|_B$
a homeomorphism onto its image.
Define $a :=\|A^{-1}\|_d^{-1}-\sigma$.
Let us show openness of $W$ and
continuity of $h$.
If $(p,z)\in W$, there exists $y\in B$
such that $f_p(y)=z$.
Let $\ve \in\,]0,r-\|y-x \|_d]$ be given.
There is an open neighborhood~$Q$ of~$p$ in~$P$
such that $d(f_q(y),f_p(y))< \frac{a\ve}{2}$
for all $q\in Q$, by continuity of~$f$.
Then, by (\ref{schachtel5}) in Theorem~\ref{newton}\,(d),
\[
f_q(B_\ve(y))\supseteq
B_{a \ve}(f_q(y))\supseteq
B_{\frac{a \ve}{2}}(f_p(y))
=  B_{\frac{a\ve}{2}}(z)\,.
\]
By the preceding, $Q\times B_{\frac{a\ve}{2}}(z))\sub W$,
whence $W$ is a neighborhood of $(p,z)$.
Furthermore, $\psi(q,z')=(f_q)^{-1}(z')\in B_\ve(y)=B_\ve((f_p)^{-1}(z))=
B_\ve(\psi(p,z))$
for all $(q,z')$ in the neighborhood
$Q\times B_{\frac{a\ve}{2}}(z)$ of $(p,z)$.
Thus $W$ is open and $\psi$
is continuous. The assertions concerning~$\xi$ follow immediately.
\end{proof}
As a direct consequence, we obtain an implicit function
theorem.
%
%\ma{lipimpl}
\begin{cor}[Continuous Implicit Functions]\label{lipimpl}
In the situation
of Theorem~{\rm\ref{newtpara}},
let $(p_0,y_0)\in P\times B$.
Then there exists an open neighborhood
$Q\sub P$ of~$p_0$ such that
$z_0:=f(p_0,y_0)\in f_p(B)$
for all $p\in Q$.
The mapping $\lambda \colon Q\to B$, $\lambda(p):=\psi(p,z_0)$
is continuous,
satisfies $\lambda (p_0)=y_0$, and
\[
\{(p,y)\in Q\times B \colon f(p,y)=z_0\}\; =\;
\graph\,(\lambda )\, .
\]
\end{cor}
\begin{proof}
Because $W$ is an open neighborhood
of $(p_0,z_0)$ in $P\times F$,
there exists an open neighborhood~$Q$
of $p_0$ in $P$
such that $Q\times \{z_0\}\sub W$.
Then $\lambda(p):=\psi(p,z_0)$ makes sense
for all $p\in Q$.
The rest is obvious from Theorem~\ref{newtpara}.
\end{proof}
\section{Inverse Function Theorem with Parameters}
We are now in the position
to formulate and prove
our main result,
an Inverse Function Theorem
with Parameters for Keller $C^k_c$-maps
in the presence of metric estimates on
partial differentials.
%
%
%
%\ma{advif}
\begin{thm}[Inverse Function Theorem with Parameters]\label{advif}
Let $E$ be a topological $\K$-vector
space, and $(F,d)$ be a metric Fr\'{e}chet
space over~$\K$, with absolutely convex balls.
Let $P_0 \sub E$ be an open subset,
or a locally convex subset with dense
interior if~$E$ is locally convex.
Let $U\sub F$
be open,
$k\in \N\cup\{\infty\}$
and $f\colon  P_0 \times U \to F$ be a $C^k_\K$-map.
Abbreviate $f_p:=f(p,\sbull)\colon  U \to F$ for $p\in P_0$.
Assume that $(p_0,x_0)\in P_0\times U$
and $f'_{p_0}(x_0)\colon F\to F$ is invertible.
Furthermore, assume that
\begin{equation}\label{nicestcase}
\sup_{(p,x)\in P_0\times U}\|\id_F-f'_{p_0}(x_0)^{-1}f'_p(x)\|_{d,d}
\;<\; 1\,;
\end{equation}
or, more generally, assume
that there exist isomorphisms
of topological vector spaces $S,A,T\colon F\to F$
such that $S\circ A \circ T\in \cL_d(F)^\times$ and
%
%\ma{vitcond}
\begin{equation}\label{vitcond}
\sup_{(p,x)\in P_0\times U}\|S\circ (A-f'_p(x))
\circ T)\|_{d,d}
\;<\; \frac{1}{\|(S\circ A \circ T)^{-1}\|_{d,d}}.
\end{equation}
Then there exists an open neighborhood $P\sub P_0$ of~$\,p_0$
and $r>0$ such that $B:=B_r(x_0)\sub U$ and the following holds:
\begin{itemize}
\item[\rm (a)]
$f_p(B)$ is open in~$F$, for each $p\in P$,
and $\phi_p\colon  B\to f_p(B)$,
$\phi_p(x):=f_p(x)=f(p,x)$ is a
$C^k_\K$-diffeomorphism.
\item[\rm (b)]
$W:=\bigcup_{p\in P} (\{p\}\times f_p(B))$ is open in $P_0\times F$,
and the map
\[
\psi\colon  W\to B, \quad
\psi(p,z)\, :=\, \phi_p^{-1}(z)
\]
is~$C^k_\K$.
Furthermore, the map
\[
\xi \colon  P \times B\to W, \quad
\xi (p,x):=(p, f(p,x))
\]
is a $C^k_\K$-diffeomorphism
with inverse $\xi^{-1}(p,z)=(p,\psi(p,z))$.
\item[\rm (c)]
$P\times B_\delta(f_{p_0}(x_0)) \sub W$ for some $\delta>0$.
\end{itemize}
In particular, for each $p\in P$ there is a unique element
$\lambda(p)\in B$ such that $f(p,\lambda(p)) =f(p_0,x_0)$,
and the map $\lambda\colon  P \to B$
so obtained is~$C^k_\K$.
\end{thm}
\begin{rem} Typical choices of $S,A,T$ are as follows:
\begin{itemize}
\item[(a)]
With $S=f_{p_0}(x_0)^{-1}$,
$A=f_{p_0}(x_0)$ and $T=\id_F$,
we recover~(\ref{nicestcase}).
\item[(b)]
If $f'_{p_0}(x_0)\in \cL_d(F)^\times$,
then a typical choice is
$A:=f'_{p_0}(x_0)$, $S:=T:=\id_F$.
In this case, we make a requirement
concerning $\sup_{x,p}\|f'_{p_0}(x_0)-f'_p(x)\|_{d,d}$.
\end{itemize}
\end{rem}
{\bf Proof of Theorem~\ref{advif}.}
Let $h:=T^{-1}\circ A^{-1}
\circ f\circ (\id_{P_0}\times T)\colon
P_0\times T^{-1}(U)\to F$. Then~$h$
satisfies
\begin{eqnarray*}
\sup_{p,x}\|\id_F -h'_p(x)\|_d
& = & \sup_{p,x}\|(SAT)^{-1}S(A-f'_p(x))T\|_d\\
& \leq &   \|(SAT)^{-1}\|_d
\,\sup_{p,x}\|S(A-f'_p(x))T\|_d
\;<\; 1\,,
\end{eqnarray*}
by~(\ref{vitcond}).
After replacing~$f$ with~$h$, we may
assume henceforth that $S=A=T=\id_F$ and
%
%\ma{vitcond2}
\begin{equation}\label{vitcond2}
\theta \; :=\; \sup_{(p,x)\in P_0\times U} \|\id_F-f'_p(x)\|_d\; <\; 1\,.
\end{equation}
Let $r>0$ such that
$B:=B_r(x_0)\sub U$.
Given $y,z\in B$ and $p\in P_0$, we have
$f_p(z)-f_p(y)=\int_0^1f'_p(y+t(z-y)).(z-y)\,dt$
and $z-y=\int_0^1 (z-y)\,dt$.
Then
\begin{eqnarray*}
\hspace*{-10mm}\lefteqn{\|f_p(z)-f_p(y)- (z-y)\|_d}\quad\\
& \leq & \sup_{t\in [0,1]}\|(f'_p(y+t(z-y))-\id_F).(z-y)\|_d\\
& \leq & \sup_{t\in [0,1]}\|f'_p(y+t(z-y))-\id_F\|_d
\| z-y\|_d \; \leq \; \theta \|z-y\|_d\, ,
\end{eqnarray*}
using
Corollary~\ref{corrfalse}
and~(\ref{vitcond2}).
Hence
\begin{equation}\label{thirty1}
\kappa \; := \; \sup\left\{
\frac{\|f_p(z)-f_p(y)-(z-y)\|_d}{\|z-y\|_d}\colon 
\,\mbox{$p\in P_0$, $z\not=y\in B$} \right\}
\; \leq \; \theta
\; < \; 1\,.
\end{equation}
Thus Theorem~\ref{newtpara} applies to $f|_{P_0\times B}$
with $A=\id_F$,
whence $f_p(B)$ is open in~$F$ and $\phi_p:=f_p|_B^{f_p(B)}$ a
homeomorphism onto its image, for each $p\in P_0$;
the set
$W:=\bigcup_{p\in P_0}\{p\}\times f_p(B)$
is open in $P_0\times F$;
the map $\psi\colon  W\to B$, $\psi(p,z):=\phi_p^{-1}(z)$ is continuous;
and the mapping
$\xi \colon  P_0\times B\to W$, $\xi(p,y):=(p,f(p,y))$
is a homeomorphism, with inverse given by
$\xi^{-1}(p,z)=(p,\psi(p,z))$.
Set $\alpha:=1-\kappa$ and $\beta:=1+\kappa$.
In view of (\ref{thirty1}),
Theorem~\ref{newton} applies to
$f_p|_B$, for each $p\in P_0$. Hence
%
%\ma{eqaddclaim}
\begin{equation}\label{eqaddclaim}
f_p(x)+ B_{\alpha s}(0) \; \sub \; f_p(B_s(x))\; \sub\;
f_p(x)+ B_{\beta s}(0)
\end{equation}
holds for all $p\in P_0$, $x\in B$
and $s\in \,]0,r-\|x-x_0\|_d]$.\\[3mm]
Also (c) is easily established:
we set $\delta:= \frac{\alpha r}{2}$.
There is an open neighborhood $P\sub P_0$ of~$p$
such that $\|f(p,x_0)-f(p_0,x_0)\|_d<\delta$
for all $p\in P$.
Then, using~(\ref{eqaddclaim})
with $x:=x_0$ and $s:=r$,
we get $f_p(B)\supseteq B_{\alpha r}(f_p(x_0))=
B_{2\delta}(f_p(x_0))
\supseteq B_\delta(f_{p_0}(x_0))$,
for all $p\in P$.
Thus (c) holds.\\[3mm]
(a) and (b): If we can show that
$\psi$ is~$C^k_\K$,
then clearly all of the maps
$\psi$, $\xi$, $\lambda$ and
$\phi_q$ will have the desired
properties. It suffices to show that~$\psi$
is~$C^k_\K$ on an open neighborhood
of each given element $(p,z)\in W$.
Given $(p,z)\in W$,
there exists $y\in B$
such that $f_p(y)=z$.
Let $\ve \in\,]0,r-\|y-x_0\|_d]$;
then $B_\ve(y)\sub B$.
There is an open neighborhood~$Q$ of~$p$ in~$P$
such that $\|f(q,y)-f(p,y)\|_d< \frac{\alpha \ve}{2}$
for all $q\in Q$, by continuity of~$f$.
Then, using (\ref{eqaddclaim}),
\[
f_q(B_\ve(y))\supseteq
B_{\alpha \ve}(f_q(y))\supseteq
B_{\frac{\alpha \ve}{2}}(f_p(y))
=B_{\frac{\alpha \ve}{2}}(z)\,.
\]
By the preceding, $Q\times B_{\frac{\alpha \ve}{2}}(z)\sub W$
and $\psi(Q\times B_{\frac{\alpha \ve}{2}}(z))\sub B_\ve(y)$. 
Now consider the $C^k_\K$-map
\[
g\colon Q\times B_{\frac{\alpha\ve}{2}}(z)
\times B_\ve (y)\to
F \,,\quad
g(q,c,v)\, :=\,
v- (f_q(v)-c)\,.
\]
For all $(q,c)\in
Q\times B_{\frac{\alpha\ve}{2}}(z)$,
the map $g_{(q,c)}:=g(q,c,\sbull)\colon B_\ve(y)\to F$
satisfies $g_{(q,c)}'(v)=\id_F-f_q'(v)$.
Thus $\sup_{(q,c)} \|g_{(q,c)}'(v)\|_d
\leq \theta<1$
using~(\ref{vitcond2}),
and so
$(g_{(q,c)})_{(q,c)}$
is a uniform family
of special contractions,
with uniform special contraction
constant~$\theta$
(see Lemma~\ref{charactspec}).
Note that $g_{(q,c)}(v)=v$ if and only if $f_q(v)=c$,
i.e., if and only if $v=\psi(q,c)$.
Thus $\psi(q,c)$
is a fixed point
of~$g_{(q,c)}$.
Since $g$ is~$C^k_\K$,
Theorem~D shows that
$\psi$ is~$C^k_\K$ on $Q\times B_{\frac{\alpha\ve}{2}}(z)$.\,\Punkt
%
%
%
%
%
%\ma{variaC1}
\begin{rem}\label{variaC1}
Theorem~\ref{advif}
remains valid in the $C^1_\K$-case
if $P_0\sub E$ is
any subset with non-empty
interior, no matter whether~$E$
and~$P_0$ are locally convex.
To achieve this, use Lemma~\ref{prepdep2}
instead of Theorem~D
at the end of the proof
of Theorem~\ref{advif}.
\end{rem}
%
%
%\ma{proofstanda}
\begin{rem}\label{proofstanda}
Note that Theorem~\ref{advif}
subsumes Theorem~A from the introduction
as its final assertion.
Using a singleton set of parameters,
we also obtain Theorem~B
as a special case.
\end{rem}
\begin{rem}
Let~$E$, $(F,d)$, $P_0$ and~$U$ be as in
Theorem~\ref{advif}.
If a $C^k_\K$-map\linebreak
$f\colon P_0\times U\to F$
satisfies
$f'_p(x)\in \cL_d(F)$
for all $(p,x)\in P_0\times U$
and also
\begin{equation}\label{toostrong}
P_0\times U \to \cL_d(F)\, ,
\quad
(p,x)\mto f'_p(x)
\end{equation}
is continuous at $(p_0,x_0)$
and $f'_{p_0}(x_0)\in \cL_d(F)^\times$,
then~(\ref{nicestcase})
is satisfied after shrinking
$P_0$ and~$U$ if necessary,
by continuity of composition in $\cL_d(F)$.
However, only (\ref{nicestcase})
(or (\ref{vitcond}))
is needed as an hypothesis
for the theorem,
not any continuity property
concerning the map in~(\ref{toostrong})
because this would be too restrictive
(at least continuity on an open set),
as we have seen
in Remark~\ref{proverpr}.
\end{rem}
\begin{rem}
We mention that $f'_p(x)\in \cL_d(F)^\times$
for all $p\in P_0$ and $x\in U$
if~(\ref{vitcond2}) holds,
because $A\in \cL_d(F)^\times $
with $A^{-1}=\sum_{n=0}^\infty (\id_F-A)^n$
for all $A\in \cL_d(F)$ such that
$\|\id_F-A\|_{d,d}<1$
(see \cite[Theorem~4.1]{Mul}).
\end{rem}
%
%
%
%
%
%
%%\ma{secBanach}
%\section{Implicit functions
%from topological vector spaces to Banach spaces}\label{secBanach}
%%
%INSERT
%%
%%
%Also inv fctn with param
%
%
%
%
%
%
%
%
\section{Application to families of linear operators}
This section describes
a simple application
of the inverse
function theorem with parameters
(Theorem~\ref{advif})
concerning the
inversion of linear operators
(which cannot be deduced
from the results in~\cite{Mul}).\\[2.5mm]
For the purposes of infinite-dimensional
Lie theory, it is useful to be able
to speak of smooth mappings
from a smooth manifold~$M$
to certain groups~$G$ which are not manifolds.
For example, $G$ might be the diffeomorphism
group of an infinite-dimensional
manifold,
or the group~$\cL(F)^\times$
of automorphisms of a (non-Banach)
topological vector space~(see~\cite{GaN}, \cite{NaV}).
We now discuss the following concept.
\begin{defn}
Let $F$ be a locally convex topological
$\K$-vector space, and
$k\in \N_0\cup\{\infty\}$.
If $k=0$, let $P$ be a topological
space. If $k\geq 1$,
let~$P$ be a locally convex subset with dense
interior
of a locally convex topological
$\K$-vector space~$E$.
Let
\[
\iota\colon \cL(F)^\times\to\cL(F)^\times\,,\quad
\iota(A):=A^{-1}
\]
be the inversion map.
We say that a map
$g\colon P \to \cL(F)^\times$
is \emph{$k$ times pseudo-differentiable}
(or~$PC^k_\K$, for short),
if both of the maps
\[
g^\wedge\colon P\times F\to F\,,
\quad g^{\wedge}(x,y):=g(x).y
\]
and $(\iota \circ g)^\wedge
\colon P\times F\to F$,
$(\iota \circ g)^\wedge(x,y):=g(x)^{-1}.y$
are~$C^k_\K$.
\end{defn}
Our application concerns a case
where the
condition on~$(\iota\circ g)^\wedge$ is superfluous.
\begin{prop}\label{parainves}
Let $F$ be a Fr\'{e}chet space
over~$\K$ and $k\in \N_0\cup\{\infty\}$.
If~$k=0$, let~$P$ be a topological space.
If $k\geq 1$,
let $E$ be a locally convex
space over~$\K$
and $P\sub E$ be a
locally convex subset with dense interior.
Let $g\colon P\to \cL(F)^\times$
be a map.
Assume that there exists
a translation invariant metric~$d$
on~$F$ defining its topology
and having absolutely convex balls,
and isomorphisms $S,A,T\colon F\to F$
of topological vector spaces with
\[
\sup_{p\in P}\, \|S(A-g(p))T\|_{d,d}\; <\; \frac{1}{\|(SAT)^{-1}\|_{d,d}}\,.
\]
Then $g$ is $PC^k_\K$
if and only if $g^\wedge\colon P \times F\to F$
is $C^k_\K$.
\end{prop}
\begin{proof}
If $k\geq 1$, then
the function $f:=g^\wedge$
satisfies the hypotheses of Theorem~\ref{advif},
with $x_0:=0$, $P_0:=P$ and
any $p_0\in P_0$.
Noting that only the proof of
Theorem~\ref{advif}\,(c) required
to shrink~$P_0$
(which is inessential for us here),
Part~(a) and~(b) of the theorem
show
that there is an open neighborhood
$W$ of $P\times \{0\}$ in $P\times F$
such that the map
\[
(\iota \circ g)^\wedge|_W\colon W\to F\,,\quad
(p,x)\mto g(p)^{-1}.x=(f_p)^{-1}(x)
\]
is~$C^k_\K$.
For $n\in \N$, set $W_n:=\{(p,nx)\colon (p,x)\in W\}$.
Then $W_n$ is open in $P\times F$ and
$\bigcup_{n\in \N}W_n=P\times F$.
Since $(\iota\circ g)^{\wedge}(p,x)=n
(\iota\circ g)^{\wedge}(p,\frac{1}{n}x)$
for each $(p,x)\in W_n$ by linearity
in the second argument,
we see that $(\iota\circ g)^\wedge|_{W_n}$
is $C^k_\K$ for each $n\in \N$.
Hence $(\iota\circ g)^\wedge$
is~$C^k_\K$.
If $k=0$, we use Theorem~\ref{newtpara}
instead Theorem~\ref{advif}
to reach the desired conclusion.
\end{proof}
Note that, if also~$E$
happens to be a Fr\'{e}chet space,
we need not assume
that $g^\wedge\colon P\times F\to F$
is $MC^1$, we only need the $C^1$-property.
This is essential, as the following example
shows.
\begin{ex}
In the situation of Example~\ref{exRN},
set $F:=\R^\N$ and consider the curve
\[
g\colon [0,1]\to \cL_d(F)\,,\quad
g(t):=\id_F-tS\,.
\]
Then $\|\id_F-g(t)\|_{d,d}=\|tS\|_{d,d}\leq
\|S\|_{d,d}=a<1$ for each $t\in [0,1]$,
and furthermore the map
\[
g^\wedge\colon [0,1]\times F\to  F\,,\quad
(t,x)\mto g(t)(y)=y  -tS.y
\]
is $C^\infty$.
By Proposition~\ref{parainves},
also the mapping
$h:=(\iota\circ g)^\wedge \colon [0,1]\times F\to F$,\linebreak
$(t,x)\mto (\id_F-tS)^{-1}(x)$
is a $C^\infty$-map.\footnote{This
implies that
also the map $[0,1]\to \cL(F)_b$, $t\mto dh((t,0),(0,\sbull))
=(\id_F-tS)^{-1}$
is~$C^\infty$, by general facts
of infinite-dimensional calculus
(see~\cite{HYP}).
Here $\cL(F)_b$ denotes $\cL(F)$, equipped
with the topology of uniform convergence
on bounded sets.}\\[2.5mm]
Since $g$ is discontinuous
as a map into $\cL_d(F)$ (because
$g(0)-\id_F=0$
but\linebreak
$\|g(t)-\id_F\|_{d,d}=a>0$
for each $t>0$), it follows that $g^\wedge$ is not~$MC^1$
(using the maximum metric on $\R\times F$).
Therefore,
we cannot get smoothness of~$h$,
say, by trying
to apply an inverse
function theorem for~$MC^1$-maps
(or $MC^\infty$-maps)
to $\R\times F\to \R\times F$, $(t,x)\mto (t,g(t).x)$.
\end{ex}
%
%
%
%
%
%
%
%
%
%
%
%
%
%
%\ma{secinvfun}
\section{Inverse and implicit {\boldmath $MC^k$}-maps}\label{secinvfun}
For later use, we record a variant
of M\"{u}ller's Inverse Function
Theorem,\footnote{Theorem~4.7 in~\cite{Mul}
and the inductive proof
of his Theorem~4.6 are our models here.}
which can do with slightly weaker
hypotheses.
%
%
%\ma{MCINV}
\begin{prop}[{\bf M\"{u}ller's Inverse Function Theorem for
{\boldmath$MC^k$}-maps}]\label{MCINV}
Let $(F,d)$ be a metric Fr\'{e}chet space
over $\K\in \{\R,\C\}$,
with absolutely convex balls.
Let $k\in \N\cup\{\infty\}$
and $f\colon U\to F$ be an $MC^k_\K$-map
on an open subset $U\sub F$.
Let $x_0\in U$.
If $f'(x_0)\in \cL(F)^\times$,
then there exists an open neighborhood
$V\sub U$ of~$x$ such that $f(V)$ is open
in~$F$ and $f|_V\colon V\to f(V)$
is an $MC^k_\K$-diffeomorphism.
\end{prop}
Before we prove the proposition,
let us record a standard consequence.
%
%
%\ma{MCIMP}
\begin{cor}[{\bf Implicit Function Theorem for
{\boldmath $MC^k$}-maps.}]\label{MCIMP}
\,Let~$(E,d_E)$\linebreak
and $(F,d_F)$ be metric
Fr\'{e}chet spaces over~$\K$,
with absolutely convex balls.
Equip $E\times F$ with the metric~$d$ given by
\[
d((x_1,y_1),(x_2,y_2))\;:=\;
\max\{d_E(x_1,x_2),d_F(y_1,y_2))\,.
\]
Let $f\colon U\times V\to F$ be an $MC^k_\K$-map,
where $U\sub E$
and $V\sub F$ are open sets.
Given $x\in U$, abbreviate $f_x:=f(x,\sbull)
\colon V\to F$.
If $f(x_0,y_0)=0$ for some $(x_0,y_0)\in U\times V$
and $f_{x_0}'(y_0)\in \cL_d(F)^\times$,
then there exist open neighborhoods
$U_0\sub U$ of~$x_0$ and $V_0\sub V$ of~$y_0$
such that
\[
\{(x,y)\in U_0\times V_0\colon f(x,y)=0\}
\; =\; \graph\lambda
\]
for an $MC^k_\K$-map $\lambda \colon U_0\to V_0$.\,\Punkt
\end{cor}
\begin{proof}
Apply Proposition~\ref{MCINV}
to the $MC^k$-map
$U\times V\to E\times F$, $g(x,y):=(x,f(x,y))$.
\end{proof}
The following lemma will help
us to deduce Proposition~\ref{MCINV}
from Theorem~B.
Its proof (recorded in Appendix~\ref{appB})
is simple but requires longish preparations,
because one has to struggle
with the fact that~$\cL_d(F)$
only is a locally convex vector group.
%
%
%\ma{invMC}
\begin{la}\label{invMC}
Let $(F,d)$
be a metric Fr\'{e}chet
space with absolutely convex balls,
and $A\in \cL_d(F)^\times$.
Let $\Omega$ be the connected component
of~$0$ of the set
\[
((A+\cL_d(F)_0)\cap \cL_d(F)^\times)-A\, .
\]
Then $(A+B)^{-1}-A^{-1}\in \cL_d(F)_0$
for each $B\in \Omega$,
and the map
$\iota_A\colon \Omega\to \cL_d(F)_0$,
$\iota_A(B):= (A+B)^{-1}-A^{-1}$
is $MC^\infty_\K$.\,\Punkt
\end{la}
{\bf Proof of Proposition~\ref{MCINV}.}
By continuity of $f'\colon U\to \cL_d(F)$,
the point~$x_0$ has an open connected neighborhood
$V\sub U$ with \mbox{$\sup_{x\in V}\|\id_F-f'(x_0)^{-1}f'(x)\|_{d,d}
<1$.} Set $y_0:=f(x_0)$.
Then $f(V)$ is an open neighborhood of~$y_0$ and $f|_V$ is
a $C^k_\K$-diffeomorphism, by
Theorem~B. Set $g:=(f|_V)^{-1}\colon f(V)\to V$.
Then $g'(y)=f'(g(y))^{-1}\in \cL_d(F)^\times$
for each $y\in f(V)$,
and the formula shows that~$g'$ is continuous.
Hence~$g$ is~$MC^1$.
By connectedness of~$V$,
$g'(f(V))\sub g'(y_0)+\cL_d(F)_0$.
Hence, setting $A:=g'(y_0)$,
we have
%\ma{borbor}
\begin{equation}\label{borbor}
g'\; =\; \iota_A\circ (\tau_{-A} \circ f')\circ g
\end{equation}
with $\iota_A$ as in Lemma~\ref{invMC}
and the translation map $\tau_{-A}\colon A+\cL_d(F)_0\to
\cL_d(F)_0$, $B\mto B-A$.
Now assume that~$g$ is $MC^{k-1}$,
by induction. Since~$\iota_A$ is
$MC^\infty$,
$\tau_{-A}\circ f'\colon U\to \cL_d(F)_0$
is $MC^{k-1}$,
it follows that $g'$ is~$MC^{k-1}$.
Hence~$g$ is~$MC^k$, which completes
the inductive proof.\Punkt
%
%
%Also $C^k$-version for impl
%into complete metric lcx vector
%groups ?.
%Does parameter-dependence etc
%still work\,?
%
%
%
%
%
%
%
%
%
%
%
%
%
%
%
%
%\ma{secglob}
\section{Global Inverse Function Theorems}\label{secglob}
In this section, we generalize
Hadamard's global inverse function theorem
from the classical Banach case
to the case of Fr\'{e}chet spaces.
%
%
%
%
%\ma{globCk}
\begin{thm}[{\bf Global\hspace*{.8mm}
Inverse\hspace*{.8mm} Function\hspace*{.8mm} Theorem\hspace*{.8mm}
for\hspace*{.8mm} {\boldmath$C^k$}-Maps}]\label{globCk}
\hfill Let\linebreak
$(E,d)$
and~$(F,d')$ be metric
Fr\'{e}chet spaces over~$\K$ with absolutely convex
balls and $f\colon E\to F$ be a $C^k_\K$-map
which is a local $C^k_\K$-diffeomorphism
around each point.
We assume that there exist isomorphisms
of topological vector spaces
$S \colon F\to F$ and $T\colon E\to E$
such that $Sf'(x)T\colon E\to F$ is
invertible with inverse\linebreak
$(Sf'(x)T)^{-1}\in \cL_{d',d}(F,E)$
for each $x\in E$,
and
%\ma{essntl}
\begin{equation}\label{essntl}
M\; :=\; \sup_{x\in E}\, \|(Sf'(x)T)^{-1}\|_{d',d}\; <\; \infty\,.
\end{equation}
Then~$f$ is a $C^k_\K$-diffeomorphism from~$E$
onto~$F$.
\end{thm}
\begin{rem}
Conditions ensuring
that $f$ is a local diffeomorphism
can be deduced from Theorem~B.
If $E=F$ and $d=d'$,
$f\colon F\to F$ is~$C^k_\K$
and
%
%\ma{ncstcs}
\begin{equation}\label{ncstcs}
\sup_{x\in F}\|\id_F-f'(x)\|_{d,d}\; <\; 1\,,
\end{equation}
then all hypotheses of Theorem~\ref{globCk}
are satisfied with $S=T=\id_F$
(noting that the estimate~(\ref{essntl})
is ensured by \cite[Theorem~4.1]{Mul}).
This case is useful for the
construction of Lie groups
of diffeomorphisms (see \S\,\ref{ligpdiff}).
The injectivity
of~$f$ is quite obvious in
the special case when~(\ref{ncstcs}) holds
(cf.\ \cite[Lemma~5.1]{DRN}).
\end{rem}
We also have a version for $MC^k$-maps.
%
%
%
%\ma{globMC}
\begin{thm}[{\bf Global Inverse Function Theorem
for {\boldmath$MC^k$}-Maps}]\label{globMC}
\hfill Let\linebreak
$(E,d)$ and~$(F,d')$ be metric
Fr\'{e}chet spaces over $\K\in \{\R,\C\}$,
with absolutely convex
balls.
Let $k\in \N\cup\{\infty\}$
and $f\colon E\to F$ be an $MC^k_\K$-map
such that $f'(x)\colon E\to F$
is invertible for each $x\in E$
and $\,\sup_{x\in E}\, \|f'(x)^{-1}\|_{d',d}<\infty$.
Then~$f$ is an $MC^k_\K$-diffeomorphism from~$E$
onto~$F$.\vspace{.7mm}
\end{thm}
{\bf Proof of Theorems~\ref{globCk} and~\ref{globMC}.}\\[1.3mm]
Step~1.
In the situation of Theorem~\ref{globCk},
after replacing $f$ with $S\circ f\circ T$, 
we may assume that $S=\id_F$ and $T=\id_E$.
In the situation of Theorem~\ref{globMC},
$f$ is a local $MC^k$-diffeomorphism,
by Theorem~\ref{MCINV}.
We therefore only need to show that $f$
is a bijection.\\[2.5mm]
Step~2.
For each continuous map $\gamma\colon Z\to F$
on a connected topological space~$Z$
and elements $z_0\in Z$
and $x_0 \in f^{-1}(\gamma(z_0))$,
there exists at most one continuous
map (``lift'') $\eta \colon Z\to E$
such that $\eta(z_0)=x_0$,
because~$E$ and~$F$ are Hausdorff
spaces and~$f$ is a local homeomorphism
(see \cite[Theorem~4.8]{For}).\\[2.5mm]
Step~3.
We now show that for each $C^1$-curve
$\gamma\colon \R\to F$
and each $x_0\in f^{-1}(\gamma(0))$,
there exists a
$C^1$-curve $\eta \colon \R\to E$
such that $f\circ \eta =\gamma$, and $\eta(0)=x_0$.\\[2mm]
Lifts being unique by Step~2,
there exists a largest
interval $I\sub \R$
such that $\gamma|_I$ admits a lift
$\eta$ with $\eta(0)=x_0$.
Since $\{0\}\to E$, $0\mto x_0$ is a lift,
$I$ is non-empty.
Because~$f$ is a local homeomorphism,
$I$~is an open interval.
Furthermore, $\eta$ is~$C^1$,
because for each $t_0\in I$,
there exists an open neighborhood
$U\sub E$ of~$\eta(t_0)$
on which~$f$ is injective,
and thus
%
%\ma{usefulfore}
\begin{equation}\label{usefulfore}
\eta|_J=(f|_U)^{-1}\circ \gamma|_J
\end{equation}
by uniqueness of lifts,
where $J\sub I$
is a connected open neighborhood
of~$t_0$ such that $\gamma(J)\sub f(U)$.\\[2.5mm]
From (\ref{usefulfore}), we also deduce that
%
%\ma{estimdiff}
\begin{equation}\label{estimdiff}
\|\eta'(t)\|_d\;\leq\; M\|\gamma'(t)\|_{d'}
\end{equation}
holds for $\eta'(t)\in E$
and $\gamma'(t)\in F$.
Write
$I=\,]a,b[$ with \mbox{${-\infty}\leq a<0<b\leq {+\infty}$.}
We show that $b=\infty$.
If not, then $J:=[0,b[$ has compact closure~$\wb{J}$
and hence
$L:=\max\{\|\gamma'(t)\|_{d'}\colon t\in \wb{J}\}<\infty$.
Using~(\ref{estimdiff}), we deduce that
\begin{equation}\label{enabcauchy}
\|\eta(t)-\eta(s)\|_d
\,\leq\,  \sup_{\tau \in J}\|\eta'(\tau)\|_d
\cdot |t-s|\,\leq\, ML\cdot |t-s|
\;\;\mbox{for all $t,s\in J$.}
\end{equation}
By (\ref{enabcauchy}), $(\eta(t))_{t\in J}$
is a Cauchy net indexed by~$J$,
which is a directed set with respect to the
order on~$J$ induced by~$\R$.
Since~$E$ is complete,
the limit $\eta(b):=\lim_{t\uparrow b}\eta(t)$
exists and provides a continuous
extension of~$\eta$ to a map
on the interval $I\cup\{b\}$.
By continuity of~$f$,
the extended function is a lift.
As~$I$ is a proper subset
of $I\cup \{b\}$, this contradicts
the maximality of~$I$.
Hence $b=\infty$, and an analogous argument
shows that $a=-\infty$.\\[2.5mm]
Step~4. $f$ is surjective.
To see this, let $z_0\in F$ be given.
Pick any $x_0\in E$ and set $y_0:=f(x_0)$.
Then $\gamma\colon \R\to F$,
$\gamma(t):=y_0+t(z_0-y_0)$
is a $C^1$-curve such that $\gamma(0)=y_0$ and
$\gamma(1)=z_0$.
By Step~3, there exists
a $C^1$-curve $\eta\colon \R\to E$
such that $\eta(0)=x_0$ and
$f\circ \eta=\gamma$.
Thus $z_0=\gamma(1)=f(\eta(1))$ in particular.\\[2.5mm]
Step~5. To see that $f$ is injective,
let $x_0, y_0\in E$
such that $f(x_0)=f(y_0)$.
Then $\eta\colon\R\to E$, $\eta(t):=x_0+t(y_0-x_0)$
is a $C^1$-curve in~$E$,
and $\gamma:=f\circ \eta\colon \R\to F$
is a $C^1$-curve such that $z_0:=\gamma(0)=\gamma(1)$.
Now consider the
$C^1$-map
\[
\Gamma\colon \R\times \R\to F\,,\quad
\Gamma(t,s):= z_0+s(\gamma(t)-z_0)\,.
\]
Then $\Gamma(0,s)=\Gamma(1,s)=z_0$
for each $s\in \R$ and $\Gamma(\sbull,0)\equiv z_0$,
$\Gamma(\sbull,1)=\gamma$.
By Step~3, for each $s\in \R$
the curve $\Gamma(\sbull,s)\colon\R\to F$
can be lifted to a curve $\zeta_s$
such that $\zeta_s(0)=x_0$.
Then $\eta=\zeta_1$,
since $f\circ \eta=\gamma=\Gamma(\sbull,1)
=f\circ \zeta_1$
and lifts are unique.
Likewise, $\zeta_0\equiv x_0$
since $\Gamma(\sbull,0)\equiv z_0$.
Now \cite[Theorem~4.10]{For}
shows that $x_0=\zeta_0(1)=\zeta_1(1)=\eta(1)=y_0$.\,\Punkt
%
%
%
%
%
%
%
%
%
%
%
%
%
%
%
%
%
%
%\ma{sec-push-it}
\section{Function spaces and mappings between them}\label{sec-push-it}
As a preliminary for our studies
of ODEs in Fr\'{e}chet spaces,
we now study differentiability properties
of certain types of mappings between
spaces of continuous
vector-valued functions on compact
topological spaces.
\begin{numba}\label{defncsp}
If $E$ is a locally convex topological $\K$-vector space
and $K$ a compact topological space,
we equip the space $C(K,E)$
of continuous $E$-valued maps
in~$K$ with the topology of uniform convergence.
This topology makes $C(K,E)$
a locally convex topological $\K$-vector
space;
the sets $C(K,U)$ with $U\sub E$
an open $0$-neighborhood
form a basis of open $0$-neighborhoods
in $C(K,E)$. If $(E,d)$ is a metric
Fr\'{e}chet space, we set
\[
\|\gamma\|_{d,\infty}\; :=\;
\max\{\|\gamma(x)\|_d\colon x\in K\}
\]
for $\gamma\in C(K,E)$ and note that
$(\gamma,\eta)\mto \|\gamma-\eta\|_{d,\infty}$
is a metric on $C(K,E)$
which defines the given topology
and has absolutely convex balls
if this is the case of~$d$.
We shall always equip $C(K,E)$
with the latter metric, and shall refer
to it as the ``maximum metric with respect to~$d$.''
\end{numba}
%
%
%\ma{relsetop}
\begin{numba}\label{relsetop}
If $U\sub E$ is an open subset,
then $C(K,U)$ is open in $C(K,E)$.
In fact, given\linebreak
$\gamma\in C(K,U)$,
the image $\gamma(K)\sub U$ is compact
and hence has a uniform neighborhood
of the form $\gamma(K)+V\sub U$
for some open $0$-neighborhood
$V\sub E$. Then $C(K,V)$ is an open $0$-neighborhood
in $C(K,E)$ and $\gamma+C(K,V)\sub C(K,U)$.
\end{numba}
%
%
%\ma{prebackb}
\begin{prop}\label{prebackb}
Let $E$ and $F$ be locally convex topological
$\K$-vector spaces,
$U\sub E$ be open,
$P$ be a topological space and~$K$
a compact topological space.
Let
\[
f\colon K\times U\times P \to F
\]
be a continuous map.
Given $p\in P$,
abbreviate $f^p:=f(\sbull,p)\colon K\times U\to F$.
Define $(f^p)_*(\gamma)\in C(K,F)$
for $\gamma\in C(K,U)$
via
\[
(f^p)_*(\gamma)(x)\; :=\; f^p(x,\gamma(x))
\;=\; f(x,\gamma(x),p)\quad \mbox{for $x\in K$.}
\]
Then the map
$\phi\colon C(K,U)\times P  \to C(K,F)$,
$\phi(\gamma,p):=(f^p)_*(\gamma)$
is continuous.
\end{prop}
\begin{proof}
Let $\gamma \in C(K,U)$,
$p\in P$, and $V\sub F$ be an open $0$-neighborhood.
Let $W\sub F$ be an open $0$-neighborhood
such that $W-W\sub V$.
For each $x\in K$,
we find an open neighborhood $A_x\sub K$
of~$x$,
an open neighborhood
$C_x\sub P$ of~$p$
and an open $0$-neighborhood
$B_x\sub E$
such that $\gamma(A_x)+B_x\sub U$ and
\[
f(y,u,q)-f(x,\gamma(x),p)\in W
\]
for all $y\in A_x$, $u\in \gamma(A_x)+B_x$,
and $q\in C_x$.
By compactness, $K\sub \bigcup_{x\in I}A_x$
for some finite subset $I\sub K$.
Then $B:=\bigcap_{x\in I}B_x\sub E$
is an open $0$-neighborhood
and $C:=\bigcap_{x\in I}C_x\sub P$
an open neighborhood of~$p$.
Let $\eta\in \gamma +C(K,B)$
and $q\in C$.
Given $y\in K$,
there is $x\in I$ with
$y\in A_x$. Then
$f(y,\eta(y),q)-f(y,\gamma(y),p) =
f(y,\eta(y),q)-f(x,\gamma(x),p)
-(f(y,\gamma(y),p)-f(x,\gamma(x),p))
\in W-W\sub V$.
We have shown that
$\phi(\eta,q)-\phi(\gamma,p)\in C(K,V)$
for all $(\eta,q)$
in the open neighborhood
$(\gamma+C(K,B))\times C$
of $(\gamma ,p)$.
Thus $\phi$ is continuous.
\end{proof}
If $k\in \N$ and $f\colon E\supseteq U\to F$
is a $C^k$-map, we can associate iterated
differentials with~$f$
via $d^0f:=f\colon U\to F$, $d^1f:=df\colon U\times E\to F$,
\[
d^2f\; :=\; d(df)\colon (U\times E)\times (E\times E)\to F
\]
(if $k\geq 2$),
and recursively $d^kf:=d^{k-1}(df)\colon U\times E^{2^k-1}\to F$.
%
%
%\ma{backbone}
%
\begin{prop}\label{backbone}
Let $K$ be a compact topological
space, $E$ and $F$ be locally
convex topological $\K$-vector spaces,
$Z$ be a topological $\K$-vector
space, and
$k\in \N_0\cup\{\infty\}$.
Let $U\sub E$ be an open subset
and $P\sub Z$ be a subset with dense
interior. If $k\geq 2$,
we assume that $P$ is open or
that both $Z$ and~$P$ are locally convex.
Let
$f\colon K \times (U\times P)\to F$
be a map such that
\begin{itemize}
\item[\rm (a)]
$f_x:=f(x,\sbull)\colon U\times P\to F$
is $C^k_\K$ for each $x\in K$; and
\item[\rm (b)]
For each $j\in \N_0$
such that $j\leq k$, the map
$K\times (U\times P)\times (E\times Z)^{2^j-1} \to F$,
\[
(x,u,p,y)\mto (d_2^jf)(x,u,p,y):=(d^jf_x)(u,p; y)
\]
for $x\in K$, $u\in U$, $p\in P$,
$y\in (E\times Z)^{2^j-1}$ 
is continuous.
\end{itemize}
Then
\[
\phi\colon
C(K,U)\times P \to C(K,F),\qquad
\phi(\gamma,p)\, :=\, (f^p)_*(\gamma)
\]
is a $C^k_\K$-map,
where $f^p:=f(\sbull,p)\colon K\times U\to F$
for $p\in P$
and
$(f^p)_*(\gamma)(x):=f(x,\gamma(x),p)$
for $x\in K$.
Furthermore, the differentials of~$\phi$ are given by
%
%\ma{diffparpush}
\begin{equation}\label{diffparpush}
\phi'(\gamma,p).(\eta,q)\;=\; (g^{(p,q)})_*(\gamma,\eta)\,,
\end{equation}
where $g^{(p,q)}:=g(\sbull,(p,q))$ with
%
%\ma{deffg}
\begin{equation}\label{deffg}
g\colon
K\times (U\times E)\times (P\times Z)\to F\,,\quad
g(x,(u,v),(p,q))\,:=\,
d_2f(x,u,p, q,v)\,.
\end{equation}
\end{prop}
\begin{proof}
We may assume that $k<\infty$.
The case $k=0$ having been settled
in Proposition~\ref{prebackb},
we may assume that $k\geq 1$.
The proof is by induction.\\[2mm]
\emph{The case $k=1$.}
Let $\gamma\in C(K,U)$, $\eta\in C(K,E)$,
$q\in Z$ and $p\in P^0$, the interior
of~$P$. Since $C(K,U)$ and~$P^0$ are open,
there exists $r>0$ such that $\gamma+B_r^\K(0)\eta
\sub C(K,U)$ and $p+B_r^\K(0)q\sub P^0$.
For each $x\in K$ and
$t\in B_r^\K(0)$
such that $t\not=0$, we have
%
%\ma{agn}
\begin{eqnarray}
\Delta_t(x) \; & := &  \frac{\phi(\gamma+t\eta,p+tq)-\phi(\gamma,p)}{t}(x)\notag\\
&=& \frac{f(x,\gamma(x)+t\eta(x),p+tq)-f(x,\gamma(x),p)}{t}\notag\\
&=&
\int_0^1 d_2f(x,(\gamma(x),p)+st(\eta(x),q); (q,\eta(x)))\,ds\,.\label{agn}
\end{eqnarray}
The map $h\colon B_r^\K(0)\times K\times [0,1]\to F$,
\[
h(t,x,s)\; :=\;
d_2f(x,(\gamma(x),p)+st(\eta(x),q); (q,\eta(x)))
\]
is continuous.
By~(\ref{agn}), the weak integral
$H(t,x):=\int_0^1h(t,x,s)\,ds$
exists in~$F$ for all $x\in K$ and $0\not=t\in B^\K_r(0)$.
But it also exists for $t=0$ because
the integrand is constant in this case.
Now the continuity of~$h$ implies
continuity of the parameter-dependent weak
integral $H\colon B_r^\K(0)\times K\to F$
(see, e.g., \cite[Chapter~1]{GaN}).
By the first half of the
exponential law (\cite[Theorem~3.4.1]{Eng}),
continuity of~$H$ implies continuity of
\[
H^\vee\colon B^\K_r(0)\to C(K,F)\,,\quad H^\vee(t)\, :=\, H(t,\sbull)\,.
\]
Since $H^\vee(t)=\Delta_t$ for $t\not=0$
by (\ref{agn}) and~$H^\vee$ is continuous,
we deduce that
$\Delta_t\to H^\vee(0)$ as $t\to 0$,
where $(H^\vee(0))(x)=H(0,x)=d_2f(x,(\gamma(x),p), (q,\eta(x)))$
for all $x\in K$.
Hence $d\phi((\gamma,p),(\eta,q))$ exists
for $(\gamma,p,\eta,q)$ as before,
and is given by~(\ref{diffparpush}).
Since~$g$ from~(\ref{deffg})
is continuous
by hypothesis~(b),
Proposition~\ref{prebackb}
shows that
the map described in~(\ref{diffparpush})
is continuous.
As the map in~(\ref{diffparpush})
extends $d\phi$ (defined so far only
on $C(K,U)\times P^0$),
we see that~$\phi$ is~$C^1_\K$ with~$d\phi$
given by~(\ref{diffparpush}).\\[2.5mm]
\emph{Induction step.}
Let $k\geq 2$ and assume that
the proposition holds when
$k$ is replaced with $k-1$.
We already know that~$\phi$ is~$C^1_\K$
and that $d\phi$ is given
by~(\ref{diffparpush}).
Since, by hypothesis~(b),
$g$ satisfies
a condition analogous to
hypothesis~(b)
with $k-1$ in place of~$k$,
the parameter-dependent pushforward in~(\ref{diffparpush})
is $C^{k-1}_\K$ by induction.
Thus $\phi$ is $C^1_\K$ with $d\phi$ a $C^{k-1}_\K$-map
and hence~$\phi$ is~$C^k_\K$.
\end{proof}
%
%
%
%\ma{deflocSCC}
\begin{defn}\label{deflocSCC}
Let $(E,d)$ and $(F,d')$ be metric Fr\'{e}chet spaces
over~$\K$,
$U\sub E$ and $X$ be a topological space.
We say that a function $f\colon X\times U\to F$
satisfies the \emph{local contraction condition}
(or ``local CC'') in its second argument,
if for $x_0\in X$ and $y_0\in U$,
there exist neighborhoods $X' \sub X$ of~$x_0$
and $U'\sub U$ of $y_0$ such that
$(f_x|_{U'})_{x\in X'}$ is a family
of special contractions,
where $f_x|_{U'}\colon U'\to F$, $y\mto f(x,y)$.
If we can always find~$X'$ and~$U'$
as before such that $(f_x|_{U'})_{x\in X'}$ is a uniform
family of special contractions, we say that~$f$
satisfies the \emph{local special contraction condition}
(or ``local SCC'')
in its second argument.
Likewise, we speak of a local SCC (resp., a local CC)
in the second argument
if $f\colon X \times U\times Z\to F$
with topological spaces~$X$ and~$Z$
and $(x_0,y_0,z_0)$ always has a box
neighborhood $X'\times U'\times Z'$
such that the maps $U'\to F$, $y\mto f(x,y,z)$
form a uniform family of special contractions
for $(x,z)\in X'\times Z'$
(resp., a uniform family of contractions).
\end{defn}
The following simple observation is useful:
%
%\ma{inheritlip}
\begin{la}\label{inheritlip}
Let $K$ be a compact topological space,
$(X,d)$ and $(Y,d')$ be
metric spaces and
$f\colon K\times X\to Y$ be
a map such that $f_x:=f(x,\sbull)\colon X\to Y$
is Lipschitz continuous for each $x\in K$
and $\theta:=\sup_{x\in K}\Lip(f_x)<\infty$.
Equip $C(K,X)$ and $C(K,Y)$
with the maximum metrics.
Then also the map
\[
f_*\colon C(K,X)\to C(K,Y)\,,\quad
f_*(\gamma)(x):=f(x,\gamma(x))\quad
\mbox{for $\gamma\in C(K,X)$, $x\in K$}
\]
is Lipschitz continuous,
with minimal Lipschitz constant
$\Lip(f_*)\leq \theta$.
In particular, if
$g\colon X\to Y$
is a Lipschitz continuous map,
then also
\[
C(K,g)\colon C(K,X)\to C(K,Y)\,,\quad
\gamma\mto g\circ \gamma
\]
is Lipschitz continuous,
with $\Lip(C(K,g))\leq \Lip(g)$.
\end{la}
\begin{proof}
Let $\gamma,\eta\in C(K,X)$.
Then
\begin{eqnarray*}
d'(f_*(\gamma)(x),f_*(\eta)(x))
& = & d'(f(x,\gamma(x)),f(x,\eta(x)))
\; \leq \; \Lip(f_x)d(\gamma(x),\eta(x))\\
& \leq & \theta \max_{y\in K}d(\gamma(y),\eta(y))
\end{eqnarray*}
for each $x\in K$ and thus
$\max_{x\in K}d'(f_*(\gamma)(x),f_*(\eta)(x))
\leq \theta \max_{y\in K}d(\gamma(y),\eta(y))$,
from which the assertions follow.
\end{proof}
\begin{rem}
Proposition~\ref{backbone}
is a variant of \cite[Proposition~3.3]{ZOO};
%in the unpublished manuscript~\cite{ZOO}.
corresponding results without parameters
are well-known
(see, e.g., \cite[Proposition~3.10]{GCX}).
Certain pushforwards (without parameters)
between certain spaces of sections
in finite-dimensional fibre bundles
(with a different type of metric)
have also been discussed in
\cite[Theorem~3.31]{Mul}.
\end{rem}
%
%
%
%
%
%
%
%
%
%
%
%
%
%
%
%
%
%
%
%
%\ma{secODE}
\section{ODEs \,in Fr\'{e}chet spaces}\label{secODE}
This section is devoted to applications.
We use our preceding results
to discuss existence and uniqueness
for solutions
to ordinary differential equations
in Fr\'{e}chet spaces,
as well
as their dependence on parameters
and initial conditions.
To this end, we adapt a classical idea
by Chow and Hale
concerning ordinary differential equations in Banach spaces
(see \cite[Chapter~3, proof of Theorem~1.1]{CaH}),
who reduced the problems in contention
to the implicit function theorem in\linebreak
Banach spaces.\\[2.5mm]
Besides the real case spelled out in Theorem~E,
also local solutions to complex differential
equations are of interest
(which are suitable complex
differentiable vector-valued maps on
a connected, locally convex subset of~$\C$
with dense interior),
but also mixed cases where we look
for ordinary solutions (on intervals
in~$\R$) with values in a complex
Fr\'{e}chet space and would like
to establish complex differentiable dependence
on initial values and parameters.
Such mixed situations
are of interest for infinite-dimensional
Lie theory, where they can simplify
the proof of regularity for a given Lie group
(cf.\ \cite[Theorem~8.1]{DL2}).\\[2.5mm]
We begin with a simple uniqueness result.
%
%
%\ma{uniqns}
\begin{prop}\label{uniqns}
Let $(F,d)$ be a metric Fr\'{e}chet space
with absolutely convex balls,
$U\sub F$ be a subset,
$J\sub \K$ be a locally convex, connected subset with dense
interior,
$f\colon J\times U\to F$ be a continuous
function
and $\gamma,\eta\colon J\to U$ be $C^1_\K$-solutions
to the differential equation
$x'(t)=f(t,x(t))$ such that $\gamma(t_0)=\eta(t_0)$
for some $t_0\in J$.
If~$f$ satisfies a local contraction condition in its
second argument,
then $\gamma=\eta$.
\end{prop}
\begin{proof}
Local uniqueness:
We show first that~$\gamma$ and~$\eta$
coincide on some neighborhood of~$t_0$.
To this end, after shrinking~$J$ and~$U$,
we may assume that~$J$
is convex,
of diameter $\leq 1$,
and that $f(t,\sbull)\colon U\to F$
is a uniform family of contractions
for $t\in J$, with some uniform
contraction constant $\theta\in \;]0,1[$.
We may also assume that $M:=\sup \|f(J\times U)\|_d<\infty$.
For each $t\in J$,
we have $\|\gamma'(t)-\eta'(t)\|_d=\|f(t,\gamma(t))-f(t,\eta(t))\|_d
\leq \min\{2M, \theta \|\gamma(t)-\eta(t)\|_d\}$.
Hence
%
%\ma{laterbett}
\begin{eqnarray}
\|\gamma(t)-\eta(t)\|_d & = & \left\|\int_0^1
(t-t_0)\cdot (\gamma'-\eta')(t_0+s(t-t_0)) \, ds\right\|_d\notag\\
& \leq & \sup_{s\in [0,1]} \|(\gamma'-\eta')(t_0+s(t-t_0))\|_d\label{47}\\
& \leq &
\theta \sup_{s\in [0,1]} \|(\gamma-\eta)(t_0+s(t-t_0))\|_d\, ,
\label{laterbett}
\end{eqnarray}
where~(\ref{47}) is also~$\leq 2M$.
Hence $\Delta:=\sup_{t\in J}\|\gamma(t)-\eta(t)\|_d<\infty$.
If $\Delta>0$, we pick $t\in J$ such that
$\Delta< \theta^{-1}\|\gamma(t)-\eta(t)\|_d$.
Since the right hand side of~(\ref{laterbett})
is $\leq \Delta$,
we obtain the contradiction
$\|\gamma(t)-\gamma(t)\|_d< \|\gamma(t)-\gamma(t)\|_d$.

(b) The set $E:=\{t\in J\colon \gamma(t)=\eta(t)\}$
is closed in~$J$ by continuity of~$\gamma$ and~$\eta$.
By~(a), $E$ is also a neighborhood in~$J$
of any of its points and hence open
in~$J$. Since $E\not=\emptyset$ (as $t_0\in E$)
and~$J$ is connected, it follows that $E=J$.
\end{proof}
%
%
%\ma{settT}
\begin{numba}\label{settT}
Our general setting is as follows.
We let
$(F,d)$ be a metric Fr\'{e}chet space
over $\K\in \{\R,\C\}$,
with absolutely convex balls.
Also, we let
$k\in \N_0\cup\{\infty\}$,
$J \sub \K$ be
a locally convex subset with dense interior,
and $\bL \in \{\R,\K\}$.
If $k=0$, we let
$P$ be a topological space
and assume that $\K=\bL=\R$.
If $k\geq 1$, we let~$E$ be a topological
$\K$-vector space
and $P\sub E$ be a subset with
dense interior.
If $k\geq 2$, we assume that
both $P\sub E$
and~$J\sub \K$ are open or that
%$P\sub E$ is open or
$E$ and~$P$ are locally convex.
We let $f\colon J \times U\times P \to F$ be a $C^k_\K$-map,
$t_0\in J \cap \bL$, $x_0\in U$ and $p_0\in P$.
\end{numba}
%
%
%\ma{exandun}
\begin{thm}[{\bf Solutions to ODEs in Fr\'{e}chet Spaces}]\label{exandun}
Let $f$ be as in~{\rm \S\,\ref{settT}}.\linebreak
If $k=0$,
assume that~$f$
satisfies a local
contraction condition in its second argument.
If $k\geq 1$,
assume that~$f$
satisfies a special local
contraction condition in its second argument.
Then there exists a convex open neighborhood
$J_1\sub J$ of~$t_0$
and open neighborhoods
$U_1\sub U$ of~$x_0$ and
$P_1\sub P$ of~$p_0$
such that
for all $(t_1,x_1,p_1)\in (J_1\cap\bL)
\times U_1\times P_1$,
the initial value problem
%
%\ma{initval2}
\begin{equation}\label{initval2}
x'(t)\;=\; f(t,x(t),p_1)\,,\qquad
x'(t_1)\;=\; x_1
\end{equation}
has a $C^k_\bL$-solution
$\phi_{t_1,x_1,p_1}\colon J_1\cap \bL  \to U$
with the following properties:
\begin{itemize}
\item[\rm(a)]
The map $\Psi \colon (J_1\cap\bL)\times (J_1\cap\bL)
\times U_1\times P_1\to U$,
$\Psi(t_1,t,x_1,p_1):=\phi_{t_1,x_1,p_1}(t)$
is~$C^k_\bL$.
\item[\rm(b)]
For fixed $(t_1,t) \in (J_1\cap\bL)\times (J_1\cap\bL)$,
the map $\Psi(t_1,t,\sbull) \colon U_1\times P_1\to F$
is~$C^k_\K$.
\item[\rm(c)]
If $(t_1,x_1,p_1)\in (J_1\cap \bL) \times U_1\times P_1$
and $\psi\colon W\to F$ is a~$C^1_\bL$-solution
to~{\rm(\ref{initval2})} on a convex neighborhood
$W\sub J_1\cap \bL$ of~$t_1$,
then $\psi=\phi_{t_1,x_1,p_1}|_W$.
\end{itemize}
\end{thm}
\begin{proof}
Once~$\Psi$ exists,
(c) is a special case of Proposition~\ref{uniqns}.
To construct solutions,
we assume first that $\K=\bL$.
By the local CC (resp.,
SCC), we may assume that $f(t,\sbull,p)\colon U\to F$,
for $(t,p)\in J\times P$,
is a uniform family
of contractions (resp., special
contractions) with uniform (special) contraction
constant $\theta\in \;]0,1[$,
after replacing~$J$, $U$ and~$P$
with smaller
neighborhoods of~$t_0$, $x_0$ and~$p_0$,
respectively
(with properties as described in the hypotheses).
We may also assume that~$J$ is convex.\\[2.5mm]
Let $V\sub U$ be an open neighborhood
of~$x_0$ and $W\sub F$ be an open
$0$-neighborhood
such that $V+W\sub U$.
Define
$g\colon [0,1]\times W\times J\times J\times V\times P \to F$,
%
%\ma{str1}
\begin{equation}\label{str1}
g(\tau,w, t_1,t_2,x_1, p_1)\; :=\;
(t_2-t_1)f(t_1+\tau(t_2-t_1),w+x_1, p_1)\,.
\end{equation}
Given $(t_1,t_2,x_1,p_1)
\in J\times J\times V\times P$,
a continuous map
$\eta \colon [0,1]\to W$
is $C^1_\R$ and satisfies
%
%\ma{str2}
\begin{equation}\label{str2}
\eta'(\tau)\; =\;
g(\tau,\eta(\tau), t_1,t_2,x_1, p_1) \;\;\;
\mbox{for all $\tau\in [0,1]$, and $\eta(0)=0$}
\end{equation}
if and only if
\[
(\forall \tau\in [0,1])\quad
\eta(\tau)\;=\; \int_0^\tau
g(\sigma ,\eta(\sigma), t_1,t_2,x_1, p_1)\,d\sigma\,,
\]
if and only if
%
%\ma{str3}
\begin{equation}\label{str3}
(\forall \tau\in [0,1])\quad
\eta(\tau)\;=\;
\int_0^1
\tau g(\sigma \tau,\eta(\sigma\tau), t_1,t_2,x_1, p_1)\,d\sigma\,.
\end{equation}
Using notation as in Proposition~\ref{backbone},
the preceding equation can be rewritten as
%
%\ma{str4}
\begin{equation}\label{str4}
(\forall \tau\in [0,1])\quad
\eta(\tau)\;=\;
\tau\cdot \int_0^1
(g^{t_1,t_2,x_1,p_1})_*(\eta)(\sigma\tau) \,d\sigma\,.
\end{equation}
To obtain a more transparent formula,
we introduce the continuous mapping\linebreak
$m\colon [0,1]\times [0,1]\to [0,1]$,
$m(\tau,\sigma):=\sigma\tau$
and the pullback
\[
C(m,F)\colon C([0,1],F)\to C([0,1]\times [0,1],F)\,,
\quad \zeta\mto \zeta\circ m
\]
which is $\K$-linear and Lipschitz continuous
with $\Lip(C(m,F))\leq 1$,
as we are using maximum metrics
on the function spaces.
Given $\zeta\in C([0,1]\times [0,1],F)$,
the map
\[
\zeta^\vee\colon [0,1]\to C([0,1],F)\,,\quad
\zeta^\vee(\tau)(\sigma):=\zeta(\tau,\sigma)
\]
is continuous and the $\K$-linear map
\[
\Phi\colon C([0,1]\times [0,1],F)\to C([0,1], C([0,1],F))\,,
\quad \Phi(\zeta):=\zeta^\vee
\]
is continuous (see \cite[Theorem~3.4.7]{Eng}).
Since maximum metrics are used
on the function spaces, it is
obvious that~$\Phi$ is isometric
and hence Lipschitz continuous
with $\Lip(\Phi)\leq 1$.
We also need the integration operator
\[
I\colon C([0,1],F)\to F\,,\quad \zeta\mto \int_0^1\zeta(\sigma)\,d\sigma
\]
which is $\K$-linear and Lipschitz continuous
with $\Lip(I)\leq 1$ (see Lemma~\ref{estinteg}).
Finally, we need the map
\[
C([0,1],I)
\colon C([0,1],C([0,1],F))\to C([0,1],F)\,,\quad
\zeta\mto I\circ \zeta
\]
which is $\K$-linear (as is clear)
and Lipschitz continuous
with $\Lip(C([0,1],I))\leq 1$
(by Lemma~\ref{inheritlip});
and the multiplication operator
\[
\mu\colon C([0,1],F)\to C([0,1],F)\,,\quad
\mu(\zeta)(\tau):=\tau\zeta(\tau)
\]
which is $\K$-linear, and Lipschitz continuous
with $\Lip(\mu)\leq 1$ (again by Lemma~\ref{inheritlip}).
We can now rewrite~(\ref{str4}) as
\[
h(t_1,t_2,x_1,p_1, \eta)\;=\; 0
\]
where
$h\colon J\times J\times V\times P\times C([0,1],W)\to C([0,1],F)$
is given by
\[
h(t_1,t_2,x_1, p_1, \eta)\; =\;
\eta-\tilde{h}(t_1,t_2,x_1,p_1,\eta)
\]
with
$\tilde{h}\colon J\times J\times V\times P\times C([0,1],W)\to C([0,1],F)$
defined via
%
%\ma{longcomp}
\begin{equation}\label{longcomp}
\tilde{h}(t_1,t_2,x_1,p_1,\eta)\; :=\;
(\mu\circ C([0,1],I)\circ \Phi\circ C(m,F)
\circ (g^{t_1,t,x_1,p_1})_*)(\eta)\,.
\end{equation}
By Lemma~\ref{inheritlip},
$\Lip((g^{t_1,t_2,x_1,p_1})_*)\leq \theta$
for all $(t_1,t_2,x_1,p_1)\in J\times J\times V\times P$.
If $k\geq 1$, for each $s\in \K^\times$
we can apply Lemma~\ref{inheritlip}
also with the metric given by~$d_s(x,y):=d(sx,sy)$
(instead of~$d$),
from which we conclude that
$(g^{t_1,t_2,x_1,p_1})_*$,
for $(t_1,t_2,x_1,p_1)\in J\times J\times V\times P$,
is a uniform family of \emph{special}
contractions with constant~$\theta$.
Since
all other maps involved in~(\ref{longcomp})
are $\K$-linear and Lipschitz continuous
with constant~$\leq 1$
(as explained before),
we deduce that $\tilde{h}(t_1,t_2,x_1,p_1,\sbull)
\colon C([0,1],W)\to C([0,1],F)$,
for $(t_1,t_2,x_1,p_1)\in J\times J\times V\times P$,
is a uniform family of contractions
(resp., of special contractions if $k\geq 1$),
with constant~$\theta$.
Furthermore, $h$ is
$C^k_\K$ as a composition
of continuous $\K$-linear maps
and a map which is~$C^k_\K$
by Proposition~\ref{backbone}.
Also, $h(t_0,t_0,x_0,p_0,0)=0$.
Hence Corollary~\ref{lipimpl}
can be applied with $A:=\id\colon C([0,1],F)\to C([0,1],F)$,
if $k=0$.
Furthermore,
Theorem~A can be applied (if $k\geq 1$)
with $A=S=T=\id$ in~(\ref{nowmoregen2}),
because the supremum on the left hand
side of~(\ref{nowmoregen2})
is $\leq \theta<1=\frac{1}{\|A\|_{D,D}}$,
by Lemma~\ref{charactspec}
and its proof (where~$D$
is the maximum metric on $C([0,1],F)$).
Now the corollary or theorem
provides open neighborhoods
$J_1\sub J$, $V_1\sub V$ and $P_1\sub P$
of $t_0$, $x_0$, resp.\ $p_0$,
and a $C^k_\K$-map
$\lambda\colon J_1\times J_1\times V_1\times P_1 \to C([0,1],F)$
such that
\[
h(t_1,t_2,x_1,p_1, \lambda(t_1,t_2,x_1,p_1))\, =\, 0\quad
\mbox{for all $(t_1,t_2,x_1,p_1)\in J_1\times J_1\times V_1\times P_1$.}
\]
\emph{The case $\K=\bL=\R$.}
Given
$(t_1,t_2,x_1,p_1)\in J_1\times J_1\times V_1\times P_1$,
consider the map $\eta:=\lambda(t_1,t_2,x_1,p_1)
\colon [0,1]\to W$.
Since~$\eta$ is continuous and satisfies~(\ref{str3}),
it is~$C^1_\R$ and satisfies~(\ref{str2}).
If $t_2\not=t_1$,
we define
$\gamma\colon [t_1,t_2]\to F$,
$t\mto x_1+\eta(\frac{t-t_1}{t_2-t_1})$
on the line segment $[t_1,t_2]$
joining~$t_1$ and~$t_2$.
Then $\gamma$ is~$C^1_\R$,
$\gamma(t_1)=x_1$,
and $\gamma'(t)=
\frac{1}{t_2-t_1}\eta(\frac{t-t_1}{t_2-t_1})
=f(t,\gamma(t),p_1)$,
using~(\ref{str2}) and expressing~$g$ in terms
of~$f$ as in~(\ref{str1}).
Hence $\gamma$ is a solution to~(\ref{initval2})
on $[t_1,t_2]$,
and its value at~$t_2$ is
$x_1+\eta(1)=x_1+\lambda(t_1,t_2,x_1,p_1)(1)$.
By uniqueness of solutions
(Proposition~\ref{uniqns}),
the former solutions
on smaller intervals combine to a solution
$\phi_{x_1,t_1,p_1}\colon J_1\to W$,
given by
%
%\ma{nowdfn}
\begin{equation}\label{nowdfn}
\phi_{t_1,x_1, p_1}(t)
\; =\;
x_1+\lambda(t_1,t,x_1,p_1)(1)\,.
\end{equation}
Since~$\lambda$ is~$C^k_\R$
and the evaluation map
%\ma{immedia}
\begin{equation}\label{immedia}
\ev_1\colon C([0,1],F)\to F\, ,\quad
\zeta\mto\zeta(1)
\end{equation}
is continuous and linear,
we deduce that~$\Psi$ (and hence also
the map in~(b)) is~$C^k_\R$.\\[2.5mm]
\emph{The case $\K=\bL=\C$}.
Then $k\geq 1$.
In this case, we simply use~(\ref{nowdfn})
to \emph{define}
$\phi_{x_1,t_1,p_1}(t)$.
Next, we define the map~$\Psi$
as in part~(a) of the theorem.
Since~$\ev_1$ (as in~(\ref{immedia}))
is continuous and complex linear,
and~$\lambda$ is~$C^k_\C$,
our definitions ensure that~$\Psi$
(and hence also the map in~(b))
is~$C^k_\C$.
By definition, $\phi_{x_1,t_1,p_1}(t_1)=x_1$,
and since~$\Psi$ is~$C^k_\C$
and hence~$C^1_\C$,
also $\phi_{x_1,t_1,p_1}$
is~$C^1_\C$.
In order that~$\phi_{x_1,t_1,p_1}$
solves~(\ref{initval2}),
it only remains to show
that $\phi_{x_1,t_1,p_1}'(t_2)=f(t_2,\phi_{x_1,t_1,p_1}(t_2),p_1)$
holds for each $t_2\in J_1$.
By continuity,
it suffices to check this for $t_2\not= t_1$.
But then we can define~$\eta$
and a $C^1_\R$-map
$\gamma\colon [t_1,t_2]\to U\sub F$
by the same formulas as in the proof
of the case~$\K=\bL=\R$,
considering now the line segment $[t_1,t_2]$
joining $t_1,t_2$
as a $1$-dimensional
real manifold with boundary immersed into~$\C$.
Again, $\gamma$ solves (\ref{initval2})
(considered now as an ODE on the manifold~$[t_1,t_2]$),
and we deduce as above that
$\gamma(\tau)=x_1+\lambda(t_1,\tau,x_1,p_1)(1)=
\phi_{t_1,x_1,p_1}(\tau)$
for each $\tau\in [t_1,t_2]$.
Calculating the complex
derivative as a suitable real directional
derivative, we find that
$\phi_{t_1,x_1,p_1}'(t_2)=
\gamma'(t_2)
=f(t_2, \phi_{t_1,x_1,p_1}(t_2) ,p_1)$,
as desired.\\[2.5mm]
\emph{The case $\K=\C$, $\bL=\R$}.
Then $k\geq 1$, and the case $\K=\bL=\C$
provides $J_1$, $U_1$ and~$P_1$
as described in the theorem
such that~(\ref{initval2})
admits a $C^k_\C$-solution $\xi_{t_1,x_1,p_1}\colon J_1\to U$
for all $t_1\in J_1$, $x_1\in U_1$ and $p_1\in P_1$,
and such that $\Theta \colon J_1\times J_1\times U_1\times P_1\to U$,
$\Theta(t_1,t,x_1,p_1):=\xi_{t_1,x_1,p_1}(t)$
is~$C^k_\C$.
Then $\phi_{t_1,x_1,p_1}:=
\xi_{t_1,x_1,p_1}|_{J_1\cap\R}
\colon J_1\cap \R\to U$
is a $C^k_\R$-solution to~(\ref{initval2})
whenever $t_1\in J_1\cap\R$,
and the map~$\Psi$ (defined in~(a))
is~$C^k_\R$, being the restriction
of the~$C^k_\C$-map~$\Theta$
to $(J_1\cap\R)\times (J_1\times \R)\times U_1\times P_1$.
Since $\Psi(t_1,t,\sbull)=\Theta(t_1,t,\sbull)$
is~$C^k_\C$, also~(b) is verified.
\end{proof}
\begin{rem}
If $F$ is a Banach space and $f\colon J\times U\times P\to F$
satisfies a Lipschitz
condition in its second argument,
then $sf$ satisfies a local SCC
in its second argument (even a global
such condition),
for $s\in \R^\times$ sufficiently
small. If $\gamma$ solves~(\ref{initval2}),
then $\eta\colon s^{-1}J_1\to U$,
$\eta(t):=\gamma(st)$
solves $\eta'(t)=sf(st,\eta(t),p_1)$,
$\eta(t)=x_1$.
Similarly, we can pass
from~$\eta$ back to~$\gamma$.
As a consequence,
all conclusions
of Theorem~\ref{exandun}
remain valid if no local CC or local SCC
is assumed, but
$F$ is a Banach space
and $f$ satisfies a local Lipschitz
condition in its second argument.
\end{rem}
\begin{rem}
Of course, we can prove existence,
uniqueness and $C^k$-dependence
just as well for higher order equations
under appropriate analogous conditions,
by rewriting them as first-order
systems.
\end{rem}
\begin{rem}
It is possible to extract
quantitative information
from the proof of Theorem~\ref{exandun}
because Theorem~A
and Corollary~\ref{lipimpl}
can be traced back
to Theorem~\ref{newton},
which provides
quantitative information on the
size of the images of balls.
For example, if~$U$ is a ball
and~$U_1$ a ball with same center
of half the radius of~$U$,
it is possible
to describe explicit conditions
on the size of the differentials
and a condition
on the diameter
of the image of~$f$
which ensure that
$J_1=J$ and $P_1=P$ can be chosen
in Theorem~\ref{exandun}.
\end{rem}
\begin{rem}\label{ODEMCk}
If~$E$ is a metric Fr\'{e}chet
space and $k\geq 1$, it is possible
to prove an analogue
of Theorem~\ref{exandun}
for~$f$ an~$MC^k_\K$-map,
in which case~$\Psi$ will
be $MC^k_\bL$ and the map in~(b)
will be $MC^k_\K$.
For the proof, note that
Proposition~\ref{backbone}
has an analogue
for~$MC^k_\K$-maps,
and use
% Proposition
Corollary~\ref{MCIMP}
instead of Theorem~A.\vspace{2mm}
\end{rem}
\begin{numba}\label{ligpdiff}
{\bf Prospect: A new class of
infinite-dimensional Lie groups.}\\[1.6mm]
Using results from this article,
it is possible to construct
certain Lie groups of rapidly decreasing
diffeomorphisms of Fr\'{e}chet spaces.\\[2.5mm]
Let $(F,d)$ be a metric Fr\'{e}chet space
with absolutely convex balls,
and $\cW$ be a set of functions
$w\colon F\to \R\cup\{\infty\}$
containing the constant function~$1$.
Let $C^\infty_\cW(F,F)$ be the ``weighted
function space''
of all $MC^\infty$-maps
$\gamma\colon F\to F$ such that
$\,\sup_{x\in F}\,|w(x)|\cdot \|\gamma(x)\|_d<\infty$
and $\,\sup_{x\in F} \, |w(x)|\cdot \|\gamma^{(k)}(x)\|_{d,d}<\infty$
for all $w\in \cW$ and $k\in \N$.
For example, take $F=\R$ and
let~$\cW$ be the set of all polynomial
functions $\R\to\R$;
then $C^\infty_\cW(\R,\R)=\cS(\R)$
is the Schwartz space of rapidly
decreasing smooth functions on~$\R$.
Returning to the general case,
let $\Diff_\cW(F)$ be the set of
all diffeomorphisms $\gamma\colon F\to F$
such that $\gamma-\id_F$,\linebreak
$\gamma^{-1}-\id_F\in C^\infty_\cW(F,F)$.
It was shown recently that
$\Diff_\cW(F)$ can be made
a Lie group modelled on $C^\infty_\cW(F,F)$,
for each Banach space~$F$
(see~\cite{Wal});
this Lie group has a smooth exponential map
and is regular (in Milnor's sense,
as in~\cite{Mil}).
Using results provided
in this article (notably, Theorem~\ref{globMC})
instead of the standard facts
of Banach differential calculus used
in~\cite{Wal},
it is possible to turn $\Diff_\cW(F)$ into
a Lie group along the lines of~\cite{Wal}.
Using the results on ODEs in Fr\'{e}chet spaces
sketched
in Remark~\ref{ODEMCk},
one also sees similarly as in the Banach case
that $\Diff_\cW(F)$ is regular.
\end{numba}
\appendix
\section{Proof of Proposition~\ref{essLspa1}}\label{appA}
\hspace*{3.2mm}(a)
We know from Remark~\ref{simplestpr}
that $\cL_{d,d'}(E,F)$ is an additive subgroup
of $F^E$. To see closedness under scalar multiplication, let
$A\in \cL_{d,d'}(E,F)$ and $t\in \K$.
Then $\|tA.x\|_{d'}\leq \max\{1,2|t|\}\,\|A.x\|_{d'}\leq
\max\{1,2|t|\}\, \|A\|_{d,d'}\|x\|_d$
for all $x\in E$ (see Lemma~\ref{splobsv}
and (\ref{standest1})).
Hence $\|tA\|_{d,d'}\leq \max\{1,2|t|\}\, \|A\|_{d,d'}<\infty$
and thus $tA\in \cL_{d,d'}(E,F)$.

(b) The evaluation map~$\ve$ is continuous
at $(0,0)$ by (\ref{standest1}) in
Remark~\ref{simplestpr}\,(a),
and furthermore $\ve(A,\sbull)$ and $\ve(\sbull,x)$
are continuous at~$0$ for all
$A\in \cL_{d,d'}(E,F)$ and $x\in E$,
by~(\ref{standest1}).
Hence~$\ve$ is
continuous, being bilinear
(see Lemma~\ref{bilinea} below).

(c)
The composition mapping
$\Gamma\colon \cL_{d',d''}(F,G)\times\cL_{d,d'}(E,F)\to
\cL_{d,d''}(E,G)$,\linebreak
$(A,B)\mto A\circ B$
is continuous at $(0,0)$
by (\ref{standest2})
in Remark~\ref{simplestpr}\,(b),
and the maps $\Gamma(A,\sbull)$ and
$\Gamma(\sbull,B)$ are continuous at~$0$,
as a consequence of~(\ref{standest2}).
Since~$\Gamma$ is bilinear,
this implies continuity of~$\Gamma$.

(d) We already know from Remark~\ref{simplestpr}\,(d)
that $D:=D_{d,d'}$ is a metric.
Now $\|t.A\|_{d,d'}\leq \|A\|_{d,d'}$
for all $A\in \cL_{d,d'}(E,F)$ and
$t\in \K$ such that $|t|\leq 1$,
by the case $|t|\leq 1$ of Lemma~\ref{splobsv}.
Hence $D$ has absolutely convex balls.

\noindent
% might need to remove again
To see that $(\cL_{d,d'}(E,F), D)$
is complete, let $(A_n)_{n\in \N}$
be a Cauchy sequence in $\cL_{d,d'}(E,F)$.
Given $x\in E$,
the point evaluation $\cL_{d,d'}(E,F)\to F$,
$B\mto B.x$ is continuous linear.
Hence $(A_n.x)_{n\in \N}$
is a Cauchy sequence in~$F$
and hence convergent,
to $A.x$ say. It is clear that the map
$A\colon E\to F$ so obtained is linear.
Given $\ve>0$, there is $N\in \N$ such that
$\|A_n-A_m\|_{d,d'}\leq\ve$ for all
$n,m\geq N$.
Given $x\in E$, this implies
that
$\|A_n.x\|_{d'}\leq
\|(A_n-A_N).x\|_{d'}+\|A_N.x\|_{d'}
\leq$\linebreak
$\|A_n-A_N\|_{d,d'}\|x\|_d+\|A_N\|_{d,d'}\|x\|_d
\leq (\ve+\|A_N\|_{d,d'})\|x\|_d$
for all $n\geq N$
and hence also
$\|A.x\|_{d'}\leq (\ve+\|A_N\|_{d,d'})\|x\|_d$,
letting $n\to \infty$.
Thus $\|A\|_{d,d'}\leq (\ve+\|A_N\|_{d,d'})\|x\|_d<\infty$
and hence $A\in \cL_{d,d'}(E,F)$.
Let $\ve$ and $N$ be as before.
Given $m \geq N$
and $x\in E$,
we have
$\|(A_n-A_m).x\|_{d'}\leq \|A_n-A_m\|_{d,d'}\|x\|_d
\leq \ve\|x\|_d$
for all $n\geq N$ and hence
$\|(A-A_m).x\|_{d'}\leq \ve\|x\|_d$,
letting $n\to\infty$.
Since $x$ was arbitrary,
we deduce that $\|A-A_m\|_{d,d'}\leq\ve$
for all $n\geq N$. Thus $A=\lim_{n\to\infty}A_n$
in $\cL_{d,d'}(E,F)$.

(e) This is \cite[Theorem~4.2]{Mul}.\,\vspace{2mm}\Punkt

\noindent
We used the following simple fact.
%
%\ma{bilinea}
\begin{la}\label{bilinea}
Let $A,B,C$ be abelian topological
groups and $\beta \colon A\times B\to C$
be a bi-additive map $($viz., a $\Z$-bilinear map$)$.
If $\beta $ is continuous at~$(0,0)$
and all of the maps
$\beta(a,\sbull)\colon B\to C$ for
$a\in A$ and $\beta(\sbull,b)\colon A\to C$
for $b\in B$
are continuous at~$0$,
then $\beta$ is continuous.
\end{la}
\begin{proof}
Let $(a_j,b_j)_{j\in J}$
be a convergent net
in $A\times B$, with limit $(a,b)$.
Since
\[
\beta(a_j,b_j)-\beta(a,b)\;=\;
\beta(a_j-a,b_j-b)+\beta(a,b_j-b)+
\beta(a_j-a,b)\,\to\, 0
\]
by the hypotheses, we see that
$\beta(a_j,b_j)\to \beta(a,b)$.
\end{proof}
%
%
%
%
%
%
%
%
%
%
%
%\ma{appB}
\section{Basic facts concerning {\boldmath $MC^k$}-maps}\label{appB}
In this appendix, we prove
compile various basic facts
concerning~$MC^k$-maps,
and deduce Lemma~\ref{invMC}
from them.\\[2.5mm]
On a product $E\times F$
of metric Fr\'{e}chet
spaces $(E,d)$ and $(F,d')$,
we shall always use the maximum metric
%
%\ma{maxmetr}
\begin{equation}\label{maxmetr}
(E\times F)^2\to [0,\infty[\,,\quad
((x_1,y_1),(x_2,y_2))\mto\max\{d(x_1,x_2),
d'(y_1,y_2)\}\,.
\end{equation}
%
%
%\ma{invMCextended}
\begin{la}\label{invMCextended}
For all metric Fr\'{e}chet
spaces
$(E,d)$, $(F,d')$ and $(G,d'')$
with absolutely convex balls,
the following holds:
\begin{itemize}
\item[\rm(a)]
Each $A\in \cL_{d,d'}(E,F)$
is an $MC^\infty_\K$-map $E\to F$.
Furthermore, the translation
$\tau_x\colon E\to E$, $y\mto x+y$
is $MC^\infty_\K$ for each $x\in E$.
\item[\rm(b)]
For each $A\in \cL_d(E)$,
both the left multiplication map
\[
\lambda_A\colon \cL_d(E)_0\to\cL_d(E)_0\,,
\quad B\mto AB
\]
and
the right multiplication map
$\rho_A\colon \cL_d(E)_0\to\cL_d(E)_0$,
$B\mto BA$ are~$MC^\infty_\K$.
More generally,
for each $A\in \cL_{d',d''}(F,G)$
and $C\in \cL_{d,d'}(E,F)$, the maps
\[
\lambda_A\colon \cL_{d,d'}(E,F)_0\to\cL_{d,d''}(E,G)_0\,,
\quad B\mto AB
\]
and $\rho_A\colon \cL_{d',d''}(F,G)_0\to\cL_{d,d''}(E,G)_0$,
$B\mto BA$ are~$MC^\infty_\K$.
\item[\rm(c)]
The map
$L\colon \cL_d(E)_0\to \cL_D(\cL_d(E)_0)_0$,
$L(A):=\lambda_A$
is $MC^\infty_\K$
and also the map
$R\colon \cL_d(E)_0\to \cL_D(\cL_d(E)_0)_0$,
$R(A):=\rho_A$ is $MC^\infty_\K$,
where $D:=D_{d,d}$ is the natural
metric on $\cL_d(E)$.
\item[\rm(d)]
Let $\beta\colon E\times F\to G$
be a bilinear map such that
$\|\beta(x,y)\|_{d''}\leq \|x\|_d\|y\|_{d'}$
for all $x\in E$, $y\in F$.
Then~$\beta$ is~$MC^\infty_\K$.
In particular, the
composition map~$\Gamma$ is $MC^\infty_\K$,
where
\[
\Gamma\colon \cL_d(E)_0\times\cL_d(E)_0\to\cL_d(E)_0\,,
\quad (A,B)\mto A\circ B\,.
\]
\item[\rm(e)]
If $U\sub E$ is a locally convex subset
with dense interior
and both $f\colon U\to F$ and $g\colon U\to G$
are~$MC^k_\K$, then also $(f,g)\colon U\to F\times G$
is~$MC^k_\K$.
\item[\rm(f)]
If $U\sub E$ and $V\sub F$ are locally convex subsets
with dense interior
and $f\colon U\to V\sub F$, $g\colon V\to G$
are~$MC^k_\K$,
then also $g\circ f\colon U\to G$
is~$MC^k_\K$.
\item[\rm(g)]
The quasi-inversion map
$q\colon Q(\cL_d(E)_0)\to \cL_d(E)_0$
is~$MC^\infty_\K$.
\end{itemize}
\end{la}
\begin{proof}
(a) Being continuous, $A$ is~$MC^0$.
Furthermore, being continuous linear,
$A$ is $C^1$ with $A'\colon E\to \cL(E,F)$,
$x\mto A$ a constant (and hence continuous)
map into $\cL_d(E,F)$.
Hence~$A$ is $MC^1$ and it follows by a trivial
induction that $A$ is $MC^k$
for each $k\geq 2$ with
$A^{(k)}=0$.

The translation $\tau_x$ is~$C^1$
with $(\tau_x)'(y)=\id_E$ for each
$y\in E$. Thus $(\tau_x)'$ is a constant map
to $\cL_d(E)$ and hence continuous.
As before, we deduce that $\tau_x$ is~$MC^\infty$
with $(\tau_x)^{(k)}=0$ for each $k\geq 2$.

(b) By (\ref{standest2})
in Remark~\ref{simplestpr}\,(b),
the map $\lambda_A$
is Lipschitz continuous.
Since $\lambda_A$ is
a linear map, it follows with~(a)
that~$\lambda_A$ is~$MC^\infty$.
The maps $\rho_A$ (and~$\rho_C$)
can be discussed analogously.

(c) Since $\|AB\|_d\leq \|A\|_d\|B\|_d$
for all $B\in \cL_d(E)$,
it follows that $\|\lambda_A\|_D\leq\|A\|_d$
and $\lambda_A\in \cL_D(\cL_d(E)_0)$.
If $A\in \cL_d(E)_0$,
then $tA\to 0$
as $t\to 0$, whence
$\|tA\|_d\to 0$ and thus
$\|t\lambda_A\|_D=\|\lambda_{tA}\|_D\leq\|tA\|_d\to 0$.
Hence $\lambda_A\in \cL_D(\cL_d(E)_0)_0$.
Summing up,
$L\colon \cL_d(E)_0\to \cL_D(\cL_d(E)_0)_0$
is a Lipschitz continuous
linear map and thus~$MC^\infty$, by~(a).
The map~$R$ can be discussed
along the same lines.

(d) Being continuous bilinear, $\beta$
is~$C^1$ with $\beta'(x,y).(u,v)=\beta(x,v)+\beta(u,y)$.
Since $\|\beta(x,v)+\beta(u,y)\|_{d''}\leq
2\max\{\|x\|_d,\|y\|_{d'}\}\cdot
\max\{\|u\|_d,\|v\|_{d'}\}$,
we deduce that $\beta'(x,y)\in \cL_{D,d''}(E\times F,G)$,
where~$D$ is the maximum metric
(as in~(\ref{maxmetr})).
Furthermore,
$\|\beta'(x,y)\|_{D,d''}\leq 2 \|(x,y)\|_D$.
Thus $\beta'\colon E\times F\to
\cL_{D,d''}(E\times F,G)$ is Lipschitz
continuous and linear.
The space $E\times F$ and hence also its
image being connected,
we deduce that $\beta'(E\times F)\sub
\cL_{D,d''}(E\times F,G)_0$.
Now $\beta'$ is $MC^\infty$ by~(a).
Hence also~$\beta$ is~$MC^\infty$.

(e) and~(f): We may assume that $k\in \N_0$.
We now prove~(e) and~(f)
in parallel for $k\in \N_0$,
by induction. In both cases, the case $k=0$ is trivial.
Thus, let~$k$ be an integer~$\geq 1$
now and assume that~(e) and~(f) hold with $k-1$
in place of~$k$.

\emph{Induction step for}~(f):
After shrinking~$U$ and~$V$,
we may assume that both sets are connected.
Pick $x_0\in U$ and set $y_0:=f(x_0)$.

We know that $g\circ f$ is~$C^1$,
with
\begin{eqnarray*}
(g\circ f)'(x) \!\! & \!\! =\!\!  & \! g'(f(x))\circ f'(x)\\
\!\! & \!\! = \!\! & \! (g'(f(x))-g'(y_0))\circ (f'(x)-f'(x_0))
+ (g'(f(x))-g'(y_0))\circ f'(x_0)\\
\!\! &\!\! \!\!  & \,  +\; g'(y_0)\circ (f'(x)-f'(x_0))
+g'(y_0)\circ f'(x_0)\, .
\end{eqnarray*}
Thus
\begin{eqnarray}
(g\circ f)' -(g\circ f)'(x_0)&=& \Gamma\circ
((g'-g'(y_0))\circ f,f'-f'(x_0))\notag\\
& & \;\; +\, \rho_{f'(x_0)}\circ ((g'-g'(y_0))\circ f)\notag\\
& & \;\; +\, \lambda_{g'(y_0)}\circ (f'-f'(x_0))\,,\label{wuergg}
\end{eqnarray}
using suitable left and right translations
(which are $MC^\infty$)
and the composition map
$\Gamma\colon \cL_{d',d''}(F,G)_0\times
\cL_{d,d'}(E,F)_0\to\cL_{d,d''}(E,G)_0$,
which is~$MC^\infty$ by~(d).
All maps involved being continuous,
we infer from~(\ref{wuergg})
that \mbox{$(g\circ f)'\colon U\to \cL_{d,d''}(E,G)$} is continuous.
Assume now that compositions of~$MC^{k-1}$-maps
are~$MC^{k-1}$.
Using the $MC^{k-1}$-case of~(e),
we then deduce from~(\ref{wuergg})
that the mapping
$(g\circ f)'-(g\circ f)'(x_0)\colon
U\to \cL_{d,d''}(E,G)_0$ is~$MC^{k-1}$,
whence $g\circ f$ is~$MC^k$.

\emph{Induction step for}~(e):
We may assume that~$U$ is connected
and pick~$x_0\in U$.
We let $\pr_1\colon F\times G\to F$
and $\pr_2\colon F\times G\to G$
be the projections onto the first
and second component, respectively.
These maps are Lipschitz continuous
and linear.
Also, we let
$\alpha\colon
\cL_{d,d'}(E,F)_0\times \cL_{d,d''}(E,G)_0
\to \cL_{d,D}(E,F\times G)_0$,
$(A,B)\mto (x\mto (Ax,Bx))$
be the natural isomorphism
of vector spaces, which is
a linear contraction.
Then
\[
(f,g)'-(f,g)'(x_0)\; =\;
\alpha \circ
(\rho_{\pr_1}\times \rho_{\pr_2})\circ (f'-f'(x_0),g'-g'(x_0))\, .
\]
Using the inductive hypotheses
(both for~(e) and~(f))
and~(b),
the preceding formula
shows that $(f,g)'-(f,g)'(x_0)$ is~$MC^{k-1}$
and thus $(f,g)$ is~$MC^k$.

(g) We already know from Proposition~\ref{isCIA}
that $Q(\cL_d(E)_0)$ is open
in $\cL_d(E)_0$ and
that~$q$ is~$C^\infty$
(and hence continuous).
Since $q(A)=\id_E-(\id_E-A)^{-1}$,
the well-known formula
$b^{-1}-a^{-1}=b^{-1}(a-b)a^{-1}$
for invertible elements
in a unital algebra
implies that
%
%\ma{ssoonn}
\begin{eqnarray}
q(B)-q(A) \! &  \! =  \! & \!(\id_E-A)^{-1}
-(\id_E-B)^{-1}
\; =\; 
(\id_E-A)^{-1}(A-B)(\id_E-B)^{-1}\notag\\
\! &\!  =  \! &\!
(q(A)-\id_E)(A-B)(q(B)-\id_E)\notag\\
\! &\! =  \! &\!
q(A)(A\!-\!B)q(B)-(A\!-\!B)q(B)-q(A)(A\!-\!B)+(A\!-\!B) \label{ssoonn}
\end{eqnarray}
for all $A,B\in Q(\cL_d(E))$.
Let $A\in Q(\cL_d(E)_0)$
and $B\in \cL_d(E)_0$ now.
For $0\not=t\in \K$
sufficiently small,
using~(\ref{ssoonn}) we see that
\[
\frac{q(A+tB)-q(A)}{t}
= 
-q(A)Bq(A+tB)+Bq(A+tB)+q(A)B-B\,,
\]
which tends to
$dq(A,B)=-q(A)Bq(A)+Bq(A)+q(A)B-B$
as $t\to 0$.
Thus, writing ${\bf 1}:=\id_{\cL_d(E)_0}$,
we have
$q'(A)+{\bf 1} =-\lambda_{q(A)}\circ\rho_{q(A)}+\rho_{q(A)}
+\lambda_{q(A)}\in\cL_D(\cL_d(E)_0)_0$
by~(b) and~(c), and 
%\ma{eninducti}
\begin{equation}\label{eninducti}
q'+{\bf 1} \, =\, -\Gamma\circ (L,R)\circ q + R\circ q+L\circ q
\end{equation}
where~$L$ and~$R$ are the $MC^\infty$-maps
from~(c) and the composition map
\[
\Gamma\colon \cL_D(\cL_d(E)_0)_0\times
\cL_D(\cL_d(E)_0)_0\to \cL_D(\cL_d(E)_0)_0
\]
is $MC^\infty$ by~(d).
Since~$q$ is continuous,~(\ref{eninducti})
shows that also the mapping\linebreak
$q'+{\bf 1} \colon Q(\cL_d(E)_0)\to \cL_D(\cL_d(E)_0)_0$
is continuous, whence~$q$ is~$MC^1$.
If~$q$ is~$MC^k$ by induction, then
(\ref{eninducti})
shows that also $q'+{\bf 1}$
is~$MC^k$ and so~$q$ is~$MC^{k+1}$.
\end{proof}
{\bf Proof of Lemma~\ref{invMC}.}
Since $\cL_d(F)^\times$ is open,
$M:=((A+\cL_d(F)_0)\cap \cL_d(F)^\times)$
is open in the affine space
$A+\cL_d(F)_0$ which is
homeomorphic to $\cL_d(F)_0$
and hence locally connected.
Thus, the connected component of~$M$ containing~$A$
is open in $A+\cL_d(F)_0$.
Since $\cL_d(F)_0\to \cL_d(F)_0$,
$C\mto C-A$ is a homeomorphism,
openness of~$\Omega$~follows.\\[2.5mm]
Since $\iota_A(B)=(\id_F+A^{-1}B)^{-1}A^{-1}-A^{-1}
=((\id_F+A^{-1}B)^{-1}-\id_F)A^{-1}
=-q(-A^{-1}B)A^{-1}$
for $B\in \Omega$
using the quasi-inversion map~$q$
of~$\cL_d(F)_0$,
we see that
$\iota_A
=-\rho_{A^{-1}}\circ q\circ (-\lambda_{A^{-1}})|_\Omega$.
Hence~$\iota_A$ is an~$MC^\infty$-map,
by Part~(b) and~(g) of Lemma~\ref{invMCextended}.\,\Punkt
{\footnotesize
\noindent
{\bf Helge Gl\"{o}ckner}, TU~Darmstadt, FB~Mathematik~AG~5,
Schlossgartenstr.\,7,\\
64289~Darmstadt, Germany.
\,E-Mail: gloeckner@mathematik.tu-darmstadt.de}
\end{document}